\documentclass{article}
\usepackage{graphicx}
\usepackage{amsmath}
\usepackage{amsfonts}
\usepackage{amssymb}
\newtheorem{theorem}{Theorem}

\newtheorem{corollary}[theorem]{Corollary}

\newtheorem{definition}[theorem]{Definition}

\newtheorem{lemma}[theorem]{Lemma}

\newtheorem{proposition}[theorem]{Proposition}

\newenvironment{proof}[1][Proof]{\textbf{#1.} }{\ \rule{0.5em}{0.5em}}
\input diagrams

\begin{document}

\title{The categorical theory of relations and quantization}
\maketitle
\author{Per K. Jakobsen, Valentin V. Lychagin\\Faculty of science, University of Troms\o , Troms\o \ 9037, Norway}
\begin{abstract}
In this paper we develops a categorical theory of relations and
use this formulation to define the notion of quantization for
relations. Categories of relations are defined in the context of
symmetric monoidal categories. They are shown to be symmetric
monoidal categories in their own right and are found to be
isomorphic to certain categories of $A-A$ bicomodules. Properties
of relations are defined in terms of the symmetric monoidal
structure. Equivalence relations are shown to be commutative
monoids in the category of relations. Quantization in our view is
a property of functors between monoidal categories. This notion of
quantization induce a deformation of all algebraic structures in
the category, in particular the ones defining properties of
relations like transitivity and symmetry.
\end{abstract}

\tableofcontents

\section{\ Introduction}

The concept of quantization is somewhat mysterious and rather ill
defined. It first appeared in a rudimentary form in the work of
Max Planck \cite{Planck1} .. Its role there was as a purely
technical device to solve a problem central to the physics of
radiation at the time, the so called ultraviolet catastrophe for
the blackbody radiation spectrum. Planck's original idea was
shortly thereafter used by Einstein to explain the photoelectric
effect \cite{Einstein1} and was further developed by N. Bohr into
what we today call the Old Quantum Theory. This theory explained
with greater precision than ever before the position of the
spectral lines for the hydrogen atom. The theory was however
rather ad hoc and it was difficult to generalize the theory to
more complicated atomic systems. The next step forward was
introduced by Louise De Broglie \cite{debroglie0},
\cite{debroglie1},\cite{debroglie2}, . He generalized the already
well known wave-particle duality for light to matter and
postulated that electrons confined to an atom would display
wavelike properties. The idea of wave-particle duality inspired E.
Schr\o dinger in 1926 to write down a wave equation for matter
waves. A different view on the notion of quantization was
introduced by Heisenberg \cite{Heis1}\cite{Heis2} in 1925 through
his matrix mechanics. These two approaches was soon shown to be
equivalent. From a modern point of view the difference in the two
approaches lies in Schr\o dingers use of the Hamiltonian
formulation of classical mechanics and of Heisenbergs use of a
formulation of classical mechanics in terms of Poissont brackets.
Schr\o dinger's approach gave rise to the canonical quantization
procedure. This procedure has been applied successfully to many
systems but contain ambiguities, like variable ordering, and has
invariance problems. The method of Geometric Quantization
\cite{Geom1} was introduced in order to resolve these problems.
Heisenbergs approach to quantization although equivalent to Schr\o
dingers approach at an elementary level, has a distinctly more
algebraic flavor than the wave mechanics of Schr\o dinger. Here
the structure of a physical system is represented in terms of an
algebra of observables. Representations of this algebra of
observables are possible models of the system in question. Whereas
algebras derived from a classical description of the system are
commutative, the algebras representing quantized systems are in
general noncommutative although still associative. Deformation
quantization \cite{Deform1},\cite{deform2} is a collection of
tools and methods that have been developed in order to find
quantized version of classical systems by deforming the algebraic
description of the system within some class of algebras. What is
clear from the existence of all these different approaches is that
the notion of quantization is not well defined. The various
approaches agree for simple systems, but they have different
domains of applicability and even for a single approach several
possible quantizations are possible for a given system. What are
the properties, or constraints, a system need in order for the
notion of quantization to be applicable? Is quantization one thing
or several different things? What is the relation between
constraints and quantizations? These are just some of the
questions that comes to mind. This paper will not give a definite
answer to any of these questions but will introduce a mathematical
framework that emphasize the idea that quantization is something
that depends on constraints and that these constraints may not
belong to the domain of mechanics or not even to physics. In fact
we believe that quantization has its natural description in terms
of a theory of representation for constraints. We also believe
that at the present time the only mathematical framework with the
right kind of generality for the formulation of a representation
theory of constraints is Category Theory \cite{MacLane}.
Constraints will in this framework take the form of relations
between natural transformations and \ a representation of the
constraints will be a category that supports all given functors
and natural transformation with the assumed relations.
Quantizations will be related to morphisms in the category of
possible representations of a given set of constraints. What we
describe here is of course a lot of bones with very little flesh.
The goal of this paper is to put a little more flesh on the bones.
This we will do by developing a theory for the quantization of
relations along the lines described above. This theory illustrate
our view of quantization, but is also of independent interest
since it gives a framework for the quantization of logic and
machines as described in the classical theory of computing. In
these days when the whole domain of classical computing is in the
process of being quantized a wider point of view on the process of
quantization is certainly needed. The categorical approach to
quantization has been introduced by one of the authors in \cite{Lych1}%
,\cite{Lych2},\cite{Lych3}.

\section{Categorical framework}

In this first chapter we formulate the basic categorical machinery that we
need in order to categorize the notion of relation. In the first subsection we
introduce the notion of a semimonoidal and a monoidal category. In line with
our general ideas of constraints and representations both notions are defined
entirely in terms of functors and natural transformations. This leads to a
slightly more general notion of monoidal category than the usual one although
we does not pursue this here. Symmetries for monoidal categories is introduced
as a further set of constraints on monoidal categories. A certain derived
relation for the natural transformations defining a symmetric monoidal
category is described and shown to be equivalent to the usual Yang-Baxter
equation. This new formulation of the Yang-Baxter equation is essential when
we later in this paper introduce a generalization of the usual notion of
symmetry that we need in order to formulate commutativity in the context of
relations. We lay the groundwork for this generalization by showing how the
Yang-Baxter equation is intimately connected to a action by a certain $S_{2}%
$-graded group. In the last subsection in this part of the paper
we introduce the notion of M-categories and C-categories. These
categories have exactly the constraints needed in order to
formulate and develop a theory of relations.

\subsection{Symmetric monoidal categories}

A semimonoidal category is a category that has a product that is
associative up to a natural isomorphism. A semimonoidal category
is a monoidal category if there is a object that is a unit for the
product up to a natural isomorphism. Properties of categories are
most clearly expressed in terms of functors and natural
transformations. We now review this formulation. On any category
we have defined the identity functor $1_{C}$. Let us assume that
there also is a bifunctor $\otimes:$ $C\times C\longrightarrow C$
defined on $C$.

\begin{definition}
A semimonoidal category is a triple $\langle$ $C,\otimes,\alpha\rangle$ where
$C$ is a category, $\otimes:$ $C\times$ $C\longrightarrow$ $C$ $\ $is a
bifunctors,
\[
\alpha:\otimes\circ(1_{C}\times\otimes)\longrightarrow\otimes\circ
(\otimes\times1_{C})
\]

is a natural isomorphism and where the following relation holds
\begin{align*}
(\alpha\circ1_{\otimes\times1_{C}\times1_{C}})\cdot(\alpha\circ1_{1_{C}%
\times1_{C}\times\otimes})  & =(1_{\otimes}\circ(\alpha\times1_{1_{C}}))\\
& \cdot(\alpha\circ1_{1_{C}\times\otimes\times1_{C}})\cdot(1_{\otimes}%
\circ(1_{1_{C}}\times\alpha))
\end{align*}
\end{definition}

A semimonoidal category is strict if $\otimes\circ(1_{C}\times\otimes
)=\otimes\circ(\otimes\times1_{C})$ and $\alpha=1_{\otimes\circ(1_{C}%
\times\otimes)}$. The relation on $\alpha$ given in the previous definition is
the object-free formulation of the usual MacLane coherence condition for the
associativity constraint $\alpha$.

For any category $C$ we have defined two bifunctors $P:C$ $\times C$
$\longrightarrow C$ and $Q:C$ $\times C$ $\longrightarrow C$. These are the
projection on the first and second factor, $P(X,Y)=X$ and $Q(X,Y)=Y$ with
obvious extension to arrows. Let $e$ be a fixed object in the category $C$ and
define a constant functor $K_{e}:C$ $\longrightarrow C$ \ by $K_{e}(X)=e$ and
$K_{e}(f)=1_{e}$. Using these functors we can give a definition of a monoidal
category entirely in terms of functors and natural transformations.

\begin{definition}
\label{moncat}A monoidal category is a 6-tuple $\langle C$ $,\otimes
,K_{e},\alpha,\beta,\gamma\rangle$ such that $\langle C$ $,\otimes
,\alpha\rangle$ is a semimonoidal category and where
\begin{align*}
\beta & :\otimes\circ(K_{e}\times1_{C})\longrightarrow Q\\
\gamma & :\otimes\circ(1_{C}\times K_{e})\longrightarrow P
\end{align*}

are natural isomorphisms such that the following relations holds
\begin{align*}
(1_{\otimes}\circ(\gamma\times1_{1_{C}}))\cdot(\alpha\circ1_{1_{C}\times
K_{e}\times1_{C}})  & =(1_{\otimes}\circ(1_{1_{C}}\times\beta))\\
\beta\circ1_{1_{C}\times K_{e}}  & =\gamma\circ1_{K_{e}\times1_{C}}%
\end{align*}
\end{definition}

A monoidal category is strict if $\langle C,\otimes,\alpha\rangle$ is a strict
semimonoidal category and if $\otimes\circ(K_{e}\times1_{C})=Q$, $\otimes
\circ(1_{C}\times K_{e})=P$ and $\beta=1_{Q}$,$\gamma=1_{P}$.

Note that $\langle C,P,1_{P\circ(1_{C,}\times P)}\rangle$ and $\langle
C,Q,1_{Q\circ(1_{C,}\times Q)}\rangle$ both are strict semimonoidal
categories. None of them can be made into a monoidal category by selecting a
unit $e$. However if $\otimes$ is part of a monoidal structure on $C$ then we
can reduce the product to projections by fixing the first and second argument
to be the unit object.

Our definition in fact deviate somewhat from the standard formulation in terms
of objects. Recall that a monoidal category in the usual sense is a 6-tuple
$\langle C,\otimes,e,\alpha^{\prime},\beta^{\prime},\gamma^{\prime}\rangle$
where $\alpha_{X,Y,Z}^{\prime}:X\otimes(Y\otimes Z)\longrightarrow(X\otimes
Y)\otimes Z$ , $\beta_{X}^{\prime}:e\otimes X\longrightarrow X$ and
$\gamma_{X}^{\prime}:X\otimes e\longrightarrow X$ are isomorphisms in $C$ that
are natural in $X,Y$, and $Z$ and where the following MacLane Coherence
\cite{MacLane} conditions are satisfied

{\small
{\tiny
\begin{diagram}
X\otimes(Y\otimes(Z\otimes T))        &      \rTo^{\alpha'_{X,Y,Z\otimes T}}%
&      (X\otimes Y)\otimes(Z\otimes T)  &   \rTo^{\alpha'_{X\otimes Y,Z,T}}%
&  ((X\otimes Y)\otimes Z)\otimes T    \\
\dTo_{1_X\otimes\alpha'_{Y,Z,T}}%
&                                                               &                                                           &                                                          &   \uTo
{\alpha'_{X,Y,Z} \otimes1_T}%
\\
X\otimes((Y\otimes Z)\otimes
T)        &                                                               &     \rTo
_{\alpha'_{X,{Y\otimes Z},T}}%
&                                                          &  (X\otimes
(Y\otimes Z))\otimes T  \\
\end{diagram}
}
}%

{\tiny
\begin{diagram}
X\otimes(e\otimes
Y)        &                                                          &  \rTo
^{\alpha'{X,e,Y}}%
&                                                         &       (X\otimes
e)\otimes Y    \\
&      \rdTo_{1_X\otimes\beta'_{Y}}%
&                                       &  \ldTo_{\gamma'_X\otimes1_Y}%
&                                               \\
&                                                           &  A\otimes
B                   &                                                         &                                               \\                                                                          &                                                \\
\end{diagram}
}%

{\tiny
\begin{diagram}
e\otimes e & \pile{\rTo^{\beta'_e} \\ \rTo_{\gamma'_e}} &   e \\
\end{diagram}
}%

It is easy to see that if we define
\begin{align*}
\beta_{X,Y}  & =\beta_{Y}^{\prime}\\
\gamma_{X,Y}  & =\gamma_{X}^{\prime}\\
\alpha_{X,Y,Z}  & =\alpha_{X,Y,Z}^{\prime}%
\end{align*}

for all objects $X$ and $Y$ in $C$, then $\langle\otimes,K_{e},\alpha
,\beta,\gamma\rangle$ is a monoidal category as defined in \ref{moncat}. If we
assume that $C$ is a category such that for all pairs of objects there exists
at least one arrow $f:X\longrightarrow X^{\prime}$. Then $K_{e}(f)=1_{e}$ and
naturality of $\beta$ implies the commutativity of the following diagram%

{\tiny
\begin{diagram}
e\otimes Y                  &  \rTo^{\beta_{X,Y}}&   Y    \\
\dTo^{1e\otimes1_Y}&                               & \dTo_{1_Y} \\
e\otimes Y                   & \rTo_{\beta_{X',Y}}&  Y     \\
\end{diagram}
}%

We thus get $\beta_{X,Y}=\beta_{X^{\prime},Y}$. In a similar way
we find $\gamma_{X,Y}=\gamma_{X,Y^{\prime}}$. This gives us a
monoidal category in the usual sense if we define
$\beta_{X}^{\prime}=\beta_{Y,X}$ and $\gamma
_{X}^{\prime}=\gamma_{X,Y}$. Our aim in this paper is not to
investigate generalizations of the notion of a monoidal category
and we will therefore assume that solutions to the relations in
\ref{moncat} satisfy $\beta _{X,Y}=\beta_{X^{\prime},Y}$ and
$\gamma_{X,Y}=\gamma_{X,Y^{\prime}}$.

We will need to express categorically the process of changing order in a
product with several factors. For any category $C$ we have the transposition
functor $\tau:C\times C\longrightarrow C\times C$ defined by $\tau(X,Y)=(Y,X)$
and $\tau(f,g)=(g,f)$. A symmetry for a monoidal category is expressed using
the functor $\tau$.

\begin{definition}
\label{symcat}A symmetric monoidal category is a 7-tuple $\langle C$
$,\otimes,K_{e},\alpha,\beta,\gamma,\sigma\rangle$ such that $\langle C$
$,\otimes,K_{e},\alpha,\beta,\gamma\rangle$ is a monoidal category and where
\[
\sigma:\otimes\longrightarrow\otimes\circ\tau
\]

is a natural isomorphism such that the following relations holds
\begin{align*}
\sigma\circ1_{\otimes\times1_{C}}  & =(\alpha^{-1}\circ1_{\tau\times1_{C}%
}\circ1_{1_{C}\times\tau})\cdot(1_{\otimes}\circ(\sigma\times1_{1_{C}})\\
& \cdot(\alpha\circ1_{1_{C}\times\tau})\cdot(1_{\otimes}\circ(1_{1_{C}}%
\times\sigma))\cdot\alpha^{-1}\\
\sigma\circ1_{1_{C}\times\otimes}  & =(\alpha\circ1_{1_{C}\times\tau}%
\circ1_{\tau\times1_{C}})\cdot(1_{\otimes}\circ(1_{1_{C}}\times\sigma
)\circ1_{\tau\times1_{C}})\\
& \cdot(\alpha^{-1}\circ1_{\tau\times1_{C}})\cdot(1_{\otimes}\circ
(\sigma\times1_{1_{C}}))\cdot\alpha\\
\beta & =(\gamma\circ1_{\tau})\cdot(\sigma\circ1_{K_{e}\times1_{C}})\\
\gamma & =(\beta\circ1_{\tau})\cdot(\sigma\circ1_{1_{C}\times K_{e}})\\
\sigma\circ1_{\tau}  & =\sigma^{-1}%
\end{align*}
\end{definition}

A symmetric monoidal category is strict if the underlying monoidal category
$\langle C$ $,\otimes,K_{e},\alpha,\beta,\gamma\rangle$ is strict.

The conditions in the definition are not independent.

\begin{proposition}
Let $\langle C$ $,\otimes,K_{e},\alpha,\beta,\gamma\rangle$ be a monoidal
category and let $\sigma:\otimes\longrightarrow\otimes\circ\tau$ be a natural
isomorphism such that $\sigma\circ1_{\tau}=\sigma^{-1}$. Then the following
two conditions are equivalent
\begin{align*}
\sigma\circ1_{\otimes\times1_{C}}  & =(\alpha^{-1}\circ1_{\tau\times1_{C}%
}\circ1_{1_{C}\times\tau})\cdot(1_{\otimes}\circ(\sigma\times1_{1_{C}}))\\
& \cdot(\alpha\circ1_{1_{C}\times\tau})\cdot(1_{\otimes}\circ(1_{1_{C}}%
\times\sigma))\cdot\alpha^{-1}\\
\sigma\circ1_{1_{C}\times\otimes}  & =(\alpha\circ1_{1_{C}\times\tau}%
\circ1_{\tau\times1_{C}})\cdot(1_{\otimes}\circ(1_{1_{C}}\times\sigma
)\circ1_{\tau\times1_{C}})\\
& \cdot(\alpha^{-1}\circ1_{\tau\times1_{C}})\cdot(1_{\otimes}\circ
(\sigma\times1_{1_{C}}))\cdot\alpha
\end{align*}
\end{proposition}

\begin{proof}
We have the following relations $\tau\circ\tau=1_{C\times C}$ and $\tau
\circ(1_{C}\times\otimes)=(\otimes\times1_{C})\circ(1_{C}\times\tau)\circ
(\tau\times1_{C})$. Using these functorial relations we have
\begin{align*}
& \sigma\circ1_{1_{C}\times\otimes}\\
& =\sigma\circ1_{\tau}\circ1_{\tau}\circ1_{1_{C}\times\otimes}\\
& =\sigma\circ1_{\tau}\circ1_{\otimes\times1_{C}}\circ1_{1_{C}\times\tau}%
\circ1_{\tau\times1_{C}}\\
& =(\sigma\circ1_{\otimes\times1_{C}})^{-1}\circ1_{1_{C}\times\tau}%
\circ1_{\tau\times1_{C}}%
\end{align*}

We thus have a relations between $\sigma\circ1_{1_{C}\times\otimes}$ and
$\sigma\circ1_{\otimes\times1_{C}}$. The equivalence of the two conditions
stated in the proposition follows directly from this relation.
\end{proof}

The third and fourth relations are also equivalent

\begin{proposition}
Let $\langle C$ $,\otimes,K_{e},\alpha,\beta,\gamma\rangle$ be a monoidal
category and let $\sigma:\otimes\longrightarrow\otimes\circ\tau$ be a natural
isomorphism such that $\sigma\circ1_{\tau}=\sigma^{-1}$. Then the following
two conditions are equivalent
\begin{align*}
\beta & =(\gamma\circ1_{\tau})\cdot(\sigma\circ1_{K_{e}\times1_{C}})\\
\gamma & =(\beta\circ1_{\tau})\cdot(\sigma\circ1_{1_{C}\times K_{e}})
\end{align*}
\end{proposition}

\begin{description}
\item
\begin{proof}
Let the first condition be given. Then we have
\begin{align*}
& \beta\circ1_{\tau}\\
& =((\gamma\circ1_{\tau})\cdot(\sigma\circ1_{K_{e}\times1_{C}}))\circ1_{\tau
}\\
& =(\gamma\circ1_{\tau}\circ1_{\tau})\cdot(\sigma\circ1_{K_{e}\times1_{C}%
}\circ1_{\tau})\\
& =\gamma\cdot(\sigma\circ1_{\tau}\circ1_{1_{C}\times K_{e}})\\
& =\gamma\cdot(\sigma^{-1}\circ1_{1_{C}\times K_{e}})
\end{align*}

and this is equivalent to the last condition.
\end{proof}
\end{description}

The symmetry conditions have a consequence that will play an important role.

\begin{proposition}
Let $\langle C$ $,\otimes,K_{e},\alpha,\beta,\gamma,\sigma\rangle$ be a
symmetric monoidal category. Then the following equation holds
\[
(\alpha\circ1_{1_{C}\times\tau}\circ1_{\tau\times1_{C}}\circ1_{1_{C}\times
\tau})\cdot(\sigma\circ(\sigma\times1_{1_{C}}))\cdot\alpha=\sigma
\circ(1_{1_{C}}\times\sigma)
\]
\end{proposition}

\begin{proof}
We have
\begin{align*}
& \sigma\circ(1_{1_{C}}\times\sigma)\\
& =(\sigma\circ(1_{1_{C}}\times1_{\otimes\circ\tau}))\cdot(1_{\otimes}%
\circ(1_{1_{C}}\times\sigma))\\
& =(\sigma\circ1_{1_{C}\times\otimes}\circ1_{1_{C}\times\tau})\cdot
(1_{\otimes}\circ(1_{1_{C}}\times\sigma))\\
& =(((\alpha\circ1_{1_{C}\times\tau}\circ1_{\tau\times1_{C}})\cdot(1_{\otimes
}\circ(1_{1_{C}}\times\sigma)\circ1_{1_{C}\times\tau})\\
& \cdot(\alpha^{-1}\circ1_{\tau\times1_{C}})\cdot(1_{\otimes}\circ
(\sigma\times1_{1_{C}}))\cdot\alpha)\circ1_{1_{C}\times\tau})\cdot(1_{\otimes
}\circ(1_{1_{C}}\times\sigma))\\
& =(\alpha\circ1_{1_{C}\times\tau}\circ1_{\tau\times1_{C}}\circ1_{1_{C}%
\times\tau})\cdot(1_{\otimes}\circ(1_{1_{C}}\times\sigma)\circ1_{\tau
\times1_{C}}\circ1_{1_{C}\times\tau})\\
& \cdot(\alpha^{-1}\circ1_{\tau\times1_{C}}\circ1_{1_{C}\times\tau}%
)\cdot(1_{\otimes}\circ(\sigma\times1_{1_{C}})\circ1_{1_{C}\times\tau}%
)\cdot(\alpha\circ1_{1_{C}\times\tau})\\
& \cdot(1_{\otimes}\circ(1_{1_{C}}\times\sigma))\\
& =(\alpha\circ1_{1_{C}\times\tau}\circ1_{\tau\times1_{C}}\circ1_{1_{C}%
\times\tau})\cdot(1_{\otimes}\circ(1_{1_{C}}\times\sigma)\circ1_{\tau
\times1_{C}}\circ1_{1_{C}\times\tau})\\
& \cdot(\sigma\circ1_{\otimes\times1_{C}})\cdot\alpha\\
& =(\alpha\circ1_{1_{C}\times\tau}\circ1_{\tau\times1_{C}}\circ1_{1_{C}%
\times\tau})\cdot(1_{\otimes\circ\tau}\circ(\sigma\times1_{1_{C}}%
))\cdot(\sigma\circ1_{\otimes\times1_{C}})\cdot\alpha\\
& =(\alpha\circ1_{1_{C}\times\tau}\circ1_{\tau\times1_{C}}\circ1_{1_{C}%
\times\tau})\cdot(\sigma\circ(\sigma\times1_{1_{C}}))\cdot\alpha
\end{align*}
\end{proof}

If we introduce the expressions for $\sigma\circ1_{\otimes\times1_{C}}$ and
$\sigma\circ1_{1_{C}\times\otimes}$ into the equation from the previous
proposition we get a equation that is cubic in $\sigma$. This equation is the
well known Yang-Baxter equation. In terms of object it takes in the strict
case the following form
\begin{align*}
& (1_{Z}\otimes\sigma_{X,Y})\circ(\sigma_{X,Z}\otimes1_{Y})\circ(1_{X}%
\otimes\sigma_{Y,Z})\\
& =(\sigma_{Y,Z}\otimes1_{X})\circ(1_{Y}\otimes\sigma_{X,Z})\circ(\sigma
_{X,Y}\otimes1_{Z})
\end{align*}

The equation from the previous proposition is clearly equivalent to the
Yang-Baxter equation in a symmetric monoidal category. We will call this
equation also for the Yang-Baxter equation. A certain generalization of this
equation will play a fundamental role in our theory of relations. This
generalization is based on characterization of symmetries in terms of a group action.

\subsection{Symmetries and group action\label{symgroupaction}}

Let $S_{2}$ be the group of permutation of two elements with the single
generator given by $t$. Let $\tau:C\times C\longrightarrow C\times C$ be the
transposition bifunctor. The functors $T_{1}=1_{C}$, $T_{2}=\tau$ and
$T_{3}=(1_{C}\times\tau)\circ(\tau\times1_{C})\circ(1_{C}\times\tau)$ defines
action of the group $S_{2}$ on the categories $C,C^{2}=C\times C$ and
$C^{3}=C\times C\times C$. Let $[C^{2},C]$ and $[C^{3},C]$ be the category of
bifunctors and trifunctors on $C$ with natural transformations as arrows. We
can induce an action of $S_{2}$ on the functor categories $[C^{2},C]$ and
$[C^{3},C]$ in the usual way by defining for objects $F$ and arrows $\alpha$
in $[C^{i},C],i=2,3$
\begin{align*}
tF  & =F\circ T_{i}\\
ta  & =\alpha\circ1_{T_{i}}%
\end{align*}

It is easy to see that this really defines an action of $S_{2}$. Let us first
consider the case when $C$ is a semimonoidal category with product $\otimes$
and associativity constraint $\alpha$. Note that
\begin{align*}
& t(\otimes\circ(1_{C}\times\otimes))\\
& =\otimes\circ(1_{C}\times\otimes)\circ(\tau\times1_{C})\circ(1_{C}\times
\tau)\circ(\tau\times1_{C})\\
& =t\otimes\circ(\otimes\times1_{C})\circ(\tau\times1_{C})\\
& =t\otimes\circ(t\otimes\times1_{C})
\end{align*}

In a similar way we find that $t(\otimes\circ(\otimes\times1_{C}%
))=t\otimes\circ(1_{C}\times t\otimes)$. We have here used the fact that
$(1_{C}\times\tau)\circ(\tau\times1_{C})\circ(1_{C}\times\tau)=(\tau
\times1_{C})\circ(1_{C}\times\tau)\circ(\tau\times1_{C})$. We therefore have a
natural isomorphism
\[
t\alpha^{-1}:t\otimes\circ(1_{C}\times t\otimes)\longrightarrow t\otimes
\circ(t\otimes\times1_{C})
\]

This is in fact an associativity constraint as the next proposition show

\begin{proposition}
$\langle C,t\otimes,t\alpha^{-1}\rangle$ is a semimonoidal category
\end{proposition}

\begin{proof}
Let $g=(1_{C}\times\tau\times1_{C})\circ(\tau\times\tau)\circ(1_{C}\times
\tau\times1_{C})\circ(\tau\times\tau)$. Then we have
\begin{align*}
& (t\alpha^{-1}\circ1_{t\otimes\times1_{C}\times1_{C}})\cdot(t\alpha^{-1}%
\circ1_{1_{C}\times1_{C}\times t\otimes})\\
& =(\alpha^{-1}\circ1_{T_{3}}\circ1_{t\otimes\times1_{C}\times1_{C}}%
)\cdot(\alpha^{-1}\circ1_{T_{3}}\circ1_{1_{C}\times1_{C}\times t\otimes})\\
& =(\alpha^{-1}\circ1_{t\otimes\times1_{C}\times1_{C}}\circ1_{g})\cdot
(\alpha^{-1}\circ1_{1_{C}\times1_{C}\times t\otimes}\circ1_{g})\\
& =[(\alpha\circ1_{1_{C}\times1_{C}\times t\otimes}\circ1_{g})\cdot
(\alpha\circ1_{t\otimes\times1_{C}\times1_{C}}\circ1_{g})]^{-1}\\
& =[((\alpha\circ1_{1_{C}\times1_{C}\times t\otimes})\cdot(\alpha
\circ1_{t\otimes\times1_{C}\times1_{C}}))\circ1_{g}]^{-1}\\
& =[((1_{\otimes}\circ(\alpha\times1_{1_{C}}))\cdot(\alpha\circ1_{1_{C}%
\times\otimes\times1_{C}})\cdot(1_{\otimes}\circ(1_{1_{C}}\times\alpha
)))\circ1_{g}]^{-1}\\
& =((1_{\otimes}\circ(1_{1_{C}}\times\alpha^{-1}))\cdot(\alpha^{-1}%
\circ1_{1_{C}\times\otimes\times1_{C}})\cdot(1_{\otimes}\circ(\alpha
^{-1}\times1_{1_{C}})))\circ1_{g}\\
& =(1_{\otimes}\circ(1_{1_{C}}\times\alpha^{-1})\circ1_{g})\cdot(\alpha
^{-1}\circ1_{1_{C}\times\otimes\times1_{C}}\circ1_{g})\cdot(1_{\otimes}%
\circ(\alpha^{-1}\times1_{1_{C}})\circ1_{g})\\
& =(1_{t\otimes}\circ(t\alpha^{-1}\times1_{1_{C}}))\cdot(t\alpha^{-1}%
\circ1_{1_{C}\times t\otimes\times1_{C}})\cdot(1_{t\otimes}\circ(1_{1_{C}%
}\times t\alpha^{-1}))
\end{align*}
\end{proof}

Let us assume that there exists a natural isomorphism $\sigma:\otimes
\longrightarrow t\otimes$ and let $\alpha$ be a associativity constraint for a
semimonoidal category $\langle C,\otimes,\alpha\rangle$. Then $t\alpha
^{-1}:t\otimes\circ(1_{C}\times t\otimes)\longrightarrow t\otimes
\circ(t\otimes\times1_{C})$ is a associativity constraint for a semimonoidal
category $\langle C,\otimes,t\alpha^{-1}\rangle$. On the other hand we have
natural isomorphisms
\begin{align*}
\sigma\circ(1_{1_{C}}\times\sigma)  & :t\otimes\circ(1_{C}\times
t\otimes)\longrightarrow\otimes\circ(1_{C}\times\otimes)\\
\sigma\circ(\sigma\times1_{1_{C}})  & :t\otimes\circ(t\otimes\times
1_{C})\longrightarrow\otimes\circ(\otimes\times1_{C})
\end{align*}

We therefore have a natural isomorphism $\widehat{\alpha}:\otimes\circ
(1_{C}\times\otimes)\longrightarrow\otimes\circ(\otimes\times1_{C})$ where we
have defined
\[
\widehat{\alpha}=(\sigma^{-1}\circ(\sigma^{-1}\times1_{1_{C}}))\cdot
t\alpha^{-1}\cdot(\sigma\circ(1_{1_{C}}\times\sigma))
\]

This new isomorphism also an associativity constraint.

\begin{proposition}
$\langle C,\otimes,\widehat{\alpha}\rangle$ is a semimonoidal category.
\end{proposition}

\begin{proof}
We only need to show that the MacLane coherence condition hold for
$\widehat{\alpha}$. Let us first observe that
\begin{align*}
& (\sigma\circ(1_{1_{C}}\times\sigma)\circ1_{\otimes\times1_{C}\times1_{C}%
})\cdot(\sigma^{-1}\circ(\sigma^{-1}\times1_{1_{C}})\circ1_{1_{C}\times
1_{C}\times\otimes})\\
& =(\sigma\circ(1_{1_{C}}\times\sigma)\circ(1_{\otimes}\times1_{1_{C}}%
\times1_{1_{C}}))\\
& \cdot(\sigma^{-1}\circ(\sigma^{-1}\times1_{1_{C}})\circ(1_{1_{C}}%
\times1_{1_{C}}\times1_{\otimes}))\\
& =(\sigma\circ((1_{1_{C}}\circ1_{\otimes})\times(\sigma\circ(1_{1_{C}}%
\times1_{1_{C}}))))\\
& \cdot(\sigma^{-1}\circ((\sigma^{-1}\circ(1_{1_{C}}\times1_{1_{C}}%
))\times(1_{1_{C}}\circ1_{\otimes})))\\
& =(\sigma\circ(1_{\otimes}\times\sigma))\cdot(\sigma^{-1}\circ(\sigma
^{-1}\times1_{\otimes}))\\
& =(1_{t\otimes}\circ(\sigma^{-1}\times\sigma))\\
& =(1_{t\otimes}\circ(\sigma^{-1}\times1_{t\otimes}))\cdot(1_{t\otimes}%
\circ(1_{t\otimes}\times\sigma))\\
& =(1_{t\otimes}\circ((1_{1_{C}}\circ\sigma^{-1})\times(1_{t\otimes}%
\circ(1_{1_{C}}\times1_{1_{C}}))))\\
& \cdot(1_{t\otimes}\circ((1_{t\otimes}\circ(1_{1_{C}}\times1_{1_{C}}%
))\times(1_{1_{C}}\circ\sigma)))\\
& =(1_{t\otimes}\circ(1_{1_{C}}\times1_{t\otimes})\circ(\sigma^{-1}%
\times1_{1_{C}}\times1_{1_{C}}))\\
& \cdot(1_{t\otimes}\circ(1_{t\otimes}\times1_{1_{C}})\circ(1_{1_{C}}%
\times1_{1_{C}}\times\sigma))
\end{align*}

Let $g=(1_{C}\times\tau\times1_{C})\circ(\tau\times\tau)\circ(1_{C}\times
\tau\times1_{C})\circ(\tau\times\tau)$.

Using the previous identity we have for the left hand side of the coherence
condition
\begin{align*}
& (\widehat{\alpha}\circ1_{\otimes\times1_{C}\times1_{C}})\circ(\widehat
{\alpha}\circ1_{1_{C}\times1_{C}\times\otimes})\\
& =(\sigma^{-1}\circ(\sigma^{-1}\times1_{1_{C}})\circ1_{\otimes\times
1_{C}\times1_{C}})\cdot(t\alpha^{-1}\circ1_{\otimes\times1_{C}\times1_{C}})\\
& \cdot(\sigma\circ(1_{1_{C}}\times\sigma)\circ1_{\otimes\times1_{C}%
\times1_{C}})\cdot(\sigma^{-1}\circ(\sigma^{-1}\times1_{1_{C}})\circ
1_{1_{C}\times1_{C}\times\otimes})\\
& (t\alpha^{-1}\circ1_{1_{C}\times1_{C}\times\otimes})\cdot(\sigma
\circ(1_{1_{C}}\times\sigma)\circ1_{1_{C}\times1_{C}\times\otimes})
\end{align*}%

\begin{align*}
& =(\sigma^{-1}\circ(\sigma^{-1}\times1_{1_{C}})\circ1_{\otimes\times
1_{C}\times1_{C}})\cdot(t\alpha^{-1}\circ(1_{\otimes}\times1_{1_{C}}%
\times1_{1_{C}}))\\
& \cdot(1_{t\otimes}\circ(1_{1_{C}}\times1_{t\otimes})\circ(\sigma^{-1}%
\times1_{1_{C}}\times1_{1_{C}}))\\
& \cdot(1_{t\otimes}\circ(1_{t\otimes}\times1_{1_{C}})\circ(1_{1_{C}}%
\times1_{1_{C}}\times\sigma))\\
& \cdot(t\alpha^{-1}\circ(1_{1_{C}}\times1_{1_{C}}\times1_{\otimes}%
))\cdot(\sigma\circ(1_{1_{C}}\times\sigma)\circ1_{1_{C}\times1_{C}%
\times\otimes})\\
& =(\sigma^{-1}\circ(\sigma^{-1}\times1_{1_{C}})\circ(1_{\otimes}%
\times1_{1_{C}}\times1_{1_{C}}))\cdot(t\alpha^{-1}\circ(\sigma^{-1}%
\times1_{1_{C}}\times1_{1_{C}}))\\
& \cdot(t\alpha^{-1}\circ(1_{1_{C}}\times1_{1_{C}}\times\sigma))\cdot
(\sigma\circ(1_{1_{C}}\times\sigma)\circ(1_{1_{C}}\times1_{1_{C}}%
\times1_{\otimes}))\\
& =(\sigma^{-1}\circ(\sigma^{-1}\times1_{1_{C}})\circ(\sigma^{-1}%
\times1_{1_{C}}\times1_{1_{C}}))\cdot(\alpha^{-1}\circ1_{T_{3}}\circ
(1_{t\otimes}\times1_{1_{C}}\times1_{1_{C}}))\\
& \cdot(\alpha^{-1}\circ1_{T_{3}}\circ(1_{1_{C}}\times1_{1_{C}}\times
1_{t\otimes}))\cdot(\sigma\circ(1_{1_{C}}\times\sigma)\circ(1_{1_{C}}%
\times1_{1_{C}}\times\sigma))\\
& =(\sigma^{-1}\circ(\sigma^{-1}\times1_{1_{C}})\circ(\sigma^{-1}%
\times1_{1_{C}}\times1_{1_{C}}))\\
& \cdot(\alpha^{-1}\circ1_{1_{C}\times1_{C}\times\otimes}\circ1_{g}%
)\cdot(\alpha^{-1}\circ1_{\otimes\times1_{C}\times1_{C}}\circ1_{g})\\
& \cdot(\sigma\circ(1_{1_{C}}\times\sigma)\circ(1_{1_{C}}\times1_{1_{C}}%
\times\sigma))\\
& =(\sigma^{-1}\circ(\sigma^{-1}\times1_{1_{C}})\circ(\sigma^{-1}%
\times1_{1_{C}}\times1_{1_{C}}))\\
& \cdot([(\alpha\circ1_{\otimes\times1_{C}\times1_{C}})\cdot(\alpha
\circ1_{1_{C}\times1_{C}\times\otimes})]^{-1}\circ1_{g})\\
& \cdot(\sigma\circ(1_{1_{C}}\times\sigma)\circ(1_{1_{C}}\times1_{1_{C}}%
\times\sigma))\\
& =(\sigma^{-1}\circ(\sigma^{-1}\times1_{1_{C}})\circ(\sigma^{-1}%
\times1_{1_{C}}\times1_{1_{C}}))\\
& \cdot([(1t\circ(\alpha\times1_{1_{C}}))\cdot(\alpha\circ1_{1_{C}%
\times\otimes\times1_{C}})\cdot(1_{\otimes}\circ(1_{1_{C}}\times\alpha
))]^{-1}\circ1_{g})\\
& \cdot(\sigma\circ(1_{1_{C}}\times\sigma)\circ(1_{1_{C}}\times1_{1_{C}}%
\times\sigma))\\
& =(\sigma^{-1}\circ(\sigma^{-1}\times1_{1_{C}})\circ(\sigma^{-1}%
\times1_{1_{C}}\times1_{1_{C}}))\cdot(1_{\otimes}\circ(1_{1_{C}}\times
\alpha^{-1})\circ1_{g})\\
& \cdot(\alpha^{-1}\circ1_{1_{C}\times\otimes\times1_{C}}\circ1_{g}%
)\cdot(1_{\otimes}\circ(\alpha^{-1}\times1_{1_{C}})\circ1_{g})\\
& \cdot(\sigma\circ(1_{1_{C}}\times\sigma)\circ(1_{1_{C}}\times1_{1_{C}}%
\times\sigma))
\end{align*}

For evaluating the righthand side of the MacLane condition we need the
following two identities
\begin{align*}
& (1_{\otimes}\circ((\sigma\circ(1_{1_{C}}\times\sigma))\times1_{1_{C}}%
))\cdot(\sigma^{-1}\circ(\sigma^{-1}\times1_{1_{C}})\circ1_{1_{C}\times
\otimes\times1_{C}})\\
& =(1_{\otimes}\circ(\sigma\times1_{1_{C}})\circ(1_{1_{C}}\times1_{\otimes
}\times1_{1_{C}}))\\
& \cdot(\sigma^{-1}\circ(\sigma^{-1}\times1_{1_{C}})\circ(1_{1_{C}}%
\times1_{\otimes}\times1_{1_{C}}))\\
& =(\sigma^{-1}\circ(1_{t\otimes}\times1_{1_{C}})\circ(1_{1_{C}}\times
\sigma\times1_{1_{C}}))\\
& =(\sigma^{-1}\circ(1_{t\otimes}\times1_{1_{C}})\circ(1_{1_{C}}%
\times1_{t\otimes}\times1_{1_{C}}))\\
& \cdot(1_{t\otimes}\circ(1_{t\otimes}\times1_{1_{C}})\circ(1_{1_{C}}%
\times1_{\otimes}\times1_{1_{C}}))
\end{align*}

and
\begin{align*}
& (\sigma\circ(1_{1_{C}}\times\sigma)\circ1_{1_{C}\times\otimes\times1_{C}%
})\cdot(1_{\otimes}\circ(1_{1_{C}}\times(\sigma^{-1}\circ(\sigma^{-1}%
\times1_{1_{C}}))))\\
& =(\sigma\circ(1_{1_{C}}\times\sigma)\circ(1_{1_{C}}\times1_{\otimes}%
\times1_{1_{C}}))\\
& \cdot(1_{\otimes}\circ(1_{1_{C}}\times\sigma^{-1})\circ(1_{1_{C}}%
\times\sigma^{-1}\times1_{1_{C}}))\\
& =(\sigma\circ(1_{1_{C}}\times1_{t\otimes})\circ(1_{1_{C}}\times\sigma
^{-1}\times1_{1_{C}}))\\
& =(1_{t\otimes}\circ(1_{1_{C}}\times1_{t\otimes})\circ(1_{1_{C}}\times
\sigma^{-1}\times1_{1_{C}}))\\
& \cdot(\sigma\circ(1_{1_{C}}\times1_{t\otimes})\circ(1_{1_{C}}\times
1_{t\otimes}\times1_{1_{C}}))
\end{align*}

Using these identities we have for the righthand side of the MacLane
condition
\begin{align*}
& (1_{\otimes}\circ(\widehat{\alpha}\times1_{1_{C}}))\cdot(\widehat{\alpha
}\circ1_{1_{C}\times\otimes\times1_{C}})\cdot(1_{\otimes}\times(1_{1_{C}%
}\times\widehat{\alpha}))\\
& =(1_{\otimes}\circ([(\sigma^{-1}\circ(\sigma^{-1}\times1_{1_{C}}))\cdot
t\alpha^{-1}\cdot(\sigma\circ(1_{1_{C}}\times\sigma))]\times1_{1_{C}}))\\
& \cdot([(\sigma^{-1}\circ(\sigma^{-1}\times1_{1_{C}}))\cdot t\alpha^{-1}%
\cdot(\sigma\circ(1_{1_{C}}\times\sigma))]\circ1_{1_{C}\times\otimes
\times1_{C}})\\
& \cdot(1_{\otimes}\circ(1_{1_{C}}\times\lbrack(\sigma^{-1}\circ(\sigma
^{-1}\times1_{1_{C}}))\cdot t\alpha^{-1}\cdot(\sigma\circ(1_{1_{C}}%
\times\sigma))]))\\
& =(1_{\otimes}\circ((\sigma^{-1}\circ(\sigma^{-1}\times1_{1_{C}}%
))\times1_{1_{C}})\cdot(1_{\otimes}\circ(t\alpha^{-1}\times1_{1_{C}}))\\
& \cdot(1_{\otimes}\circ((\sigma\circ(1_{1_{C}}\times\sigma))\times1_{1_{C}%
}))\cdot(\sigma^{-1}\circ(\sigma^{-1}\times1_{1_{C}})\circ1_{1_{C}%
\times\otimes\times1_{C}})\\
& \cdot(t\alpha^{-1}\circ1_{1_{C}\times\otimes\times1_{C}})\cdot(\sigma
\circ(1_{1_{C}}\times\sigma)\circ1_{1_{C}\times\otimes\times1_{C}})\\
& \cdot(1_{\otimes}\circ(1_{1_{C}}\times(\sigma^{-1}\circ(\sigma^{-1}%
\times1_{1_{C}}))))\cdot(1_{\otimes}\circ(1_{1_{C}}\times t\alpha^{-1}))\\
& \cdot(1_{\otimes}\circ(1_{1_{C}}\times(\sigma\circ(1_{1_{C}}\times
\sigma))))\\
& =(1_{\otimes}\circ(\sigma^{-1}\times1_{1_{C}})\circ(\sigma^{-1}%
\times1_{1_{C}}\times1_{1_{C}}))\cdot(1_{\otimes}\circ(t\alpha^{-1}%
\times1_{1_{C}}))\\
& \cdot(\sigma^{-1}\circ(1_{t\otimes}\times1_{1_{C}})\circ(1_{1_{C}}%
\times1_{t\otimes}\times1_{1_{C}}))\\
& \cdot(1_{t\otimes}\circ(1_{t\otimes}\times1_{1_{C}})\circ(1_{1_{C}}%
\times\sigma\times1_{1_{C}}))\cdot(t\alpha^{-1}\circ(1_{1_{C}}\times
1_{\otimes}\times1_{1_{C}}))\\
& \cdot(1_{t\otimes}\circ(1_{1_{C}}\times1_{t\otimes})\circ(1_{1_{C}}%
\times\sigma^{-1}\times1_{1_{C}}))\\
& \cdot(\sigma\circ(1_{1_{C}}\times1_{t\otimes})\circ(1_{1_{C}}\times
1_{t\otimes}\times1_{1_{C}}))\cdot(1_{\otimes}\circ(1_{1_{C}}\times
t\alpha^{-1}))\\
& \cdot(1_{\otimes}\circ(1_{1_{C}}\times\sigma)\circ(1_{1_{C}}\times1_{1_{C}%
}\times\sigma))\\
& =(1_{\otimes}\circ(\sigma^{-1}\times1_{1_{C}})\circ(\sigma^{-1}%
\times1_{1_{C}}\times1_{1_{C}}))\cdot(1_{\otimes}\circ(t\alpha^{-1}%
\times1_{1_{C}}))\\
& \cdot(\sigma^{-1}\circ((1_{t\otimes}\circ(1_{1_{C}}\times1_{t\otimes
}))\times1_{1_{C}}))\cdot(1_{t\otimes}\circ(1_{t\otimes}\times1_{1_{C}}%
)\circ(1_{1_{C}}\times\sigma\times1_{1_{C}}))\\
& \cdot(t\alpha^{-1}\circ(1_{1_{C}}\times1_{\otimes}\times1_{1_{C}}%
))\cdot(1_{t\otimes}\circ(1_{1_{C}}\times1_{t\otimes})\circ(1_{1_{C}}%
\times\sigma^{-1}\times1_{1_{C}}))\\
& \cdot(\sigma\circ(1_{1_{C}}\times(1_{t\otimes}\circ(1_{t\otimes}%
\times1_{1_{C}}))))\cdot(1_{\otimes}\circ(1_{1_{C}}\times t\alpha^{-1}))\\
& \cdot(1_{\otimes}\circ(1_{1_{C}}\times\sigma)\circ(1_{1_{C}}\times1_{1_{C}%
}\times\sigma))\\
& =(1_{\otimes}\circ((\sigma^{-1}\circ(\sigma^{-1}\times1_{1_{C}}%
))\times1_{1_{C}}))\cdot(\sigma^{-1}\circ(t\alpha^{-1}\times1_{1_{C}}))\\
& \cdot(t\alpha^{-1}\circ(1_{1_{C}}\times1_{t\otimes}\times1_{1_{C}}%
))\cdot(\sigma\circ(1_{1_{C}}\times t\alpha^{-1}))\\
& \cdot(1_{\otimes}\circ(1_{1_{C}}\times(\sigma\circ(1_{1_{C}}\times
\sigma))))\\
& =(\sigma^{-1}\circ((\sigma^{-1}\circ(\sigma^{-1}\times1_{1_{C}}%
))\times1_{1_{C}}))\cdot(1_{t\otimes}\circ(t\alpha^{-1}\times1_{1_{C}}))\\
& \cdot(t\alpha^{-1}\circ(1_{1_{C}}\times1_{t\otimes}\times1_{1_{C}}%
))\cdot(1_{t\otimes}\circ(1_{1_{C}}\times t\alpha^{-1}))\\
& \cdot(\sigma\circ(1_{1_{C}}\times(\sigma\circ(1_{1_{C}}\times\sigma))))\\
& =(\sigma^{-1}\circ((\sigma^{-1}\circ(\sigma^{-1}\times1_{1_{C}}%
))\times1_{1_{C}}))\cdot(1_{\otimes}\circ(1_{1_{C}}\times\alpha^{-1}%
)\circ1_{g})\\
& \cdot(\alpha^{-1}\circ1_{1_{C}\times\otimes\times1_{C}}\circ1_{g}%
)\cdot(1_{\otimes}\circ(\alpha^{-1}\times1_{1_{C}})\circ1_{g})\\
& \cdot(\sigma\circ(1_{1_{C}}\times(\sigma\circ(1_{1_{C}}\times\sigma))))
\end{align*}

The lefthand side and the righthand side are thus equal and this proves the proposition.
\end{proof}

Let us define $S_{C,\otimes}=\{\alpha\;|\;\langle C,\otimes,\alpha\rangle\;$is
a semimonoidal category $\}$. Then the previous proposition show that for each
natural isomorphism $\sigma:\otimes\longrightarrow t\otimes$ we have a mapping
of $S_{C,\otimes}$ to itself.

Let us next consider the case of a monoidal category $\langle C,\otimes
,K_{e},\alpha,\beta,\gamma\rangle$. Using the natural isomorphism $\sigma$ we
can define new natural isomorphisms
\begin{align*}
\widehat{\beta}  & =(t\gamma)\cdot(\sigma\circ1_{K_{e}\times1_{C}}%
):\otimes\circ(K_{e}\times1_{C})\longrightarrow Q\\
\widehat{\gamma}  & =(t\beta)\cdot(\sigma\circ1_{1_{C}\times K_{e}}%
):\otimes\circ(1_{C}\times K_{e})\longrightarrow P
\end{align*}

For $\widehat{\alpha}$ and the two natural isomorphisms $\widehat{\beta}$ and
$\widehat{\gamma}$ we have

\begin{proposition}
$\langle C,\otimes,K_{e},\widehat{\alpha},\widehat{\beta},\widehat{\gamma
}\rangle$ is a monoidal category
\end{proposition}

\begin{proof}
The First MacLane coherence condition has already been verified. For the
second MacLane condition we need the identities
\begin{align*}
& (1_{\otimes}\circ((\sigma\circ1_{1_{C}\times K_{e}})\times1_{1_{C}}%
))\cdot((\sigma^{-1}\circ(\sigma^{-1}\times1_{1_{C}}))\circ1_{1_{C}\times
K_{e}\times1_{C}})\\
& =(1_{\otimes}\circ(\sigma\times1_{1_{C}})\circ(1_{1_{C}}\times1_{K_{e}%
}\times1_{1_{C}}))\\
& \cdot(\sigma^{-1}\circ(\sigma^{-1}\times1_{1_{C}})\circ(1_{1_{C}}%
\times1_{K_{e}}\times1_{1_{C}}))\\
& =(\sigma^{-1}\circ(1_{t\otimes}\times1_{1_{C}})\circ(1_{1_{C}}\times
1_{K_{e}}\times1_{1_{C}}))\\
& =(\sigma^{-1}\circ(1_{t\otimes}\circ(1_{1_{C}}\times1_{K_{e}})\times
1_{1_{C}}))
\end{align*}

and
\begin{align*}
& (1_{t\otimes}\circ(t\beta\times1_{1_{C}}))\cdot(t\alpha^{-1}\circ
1_{1_{C}\times K_{e}\times1_{C}})\\
& =(1_{\otimes}\circ(1_{1_{C}}\times\beta)\circ1_{T_{3}})\cdot(\alpha
^{-1}\circ1_{1_{C}\times K_{e}\times1_{C}}\circ1_{T_{3}})\\
& =(((1_{\otimes}\circ(1_{1_{C}}\times\beta))\cdot(\alpha^{-1}\circ
1_{1_{C}\times K_{e}\times1_{C}}))\circ1_{T_{3}})\\
& =(1_{\otimes}\circ(\gamma\times1_{1_{C}})\circ1_{T_{3}})
\end{align*}

Using these two identities we have
\begin{align*}
& (1_{\otimes}\circ(\widehat{\gamma}\times1_{1_{C}}))\cdot(\widehat{\alpha
}\circ1_{1_{C}\times K_{e}\times1_{C}})\\
& =(1_{\otimes}\circ((t\beta\cdot(\sigma\circ1_{1_{C}\times K_{e}}%
))\times1_{1_{C}}))\\
& \cdot(((\sigma^{-1}\circ(\sigma^{-1}\times1_{1_{C}}))\cdot t\alpha^{-1}%
\cdot(\sigma\circ(1_{1_{C}}\times\sigma)))\circ1_{1_{C}\times K_{e}\times
1_{C}})\\
& =(1_{\otimes}\circ(t\beta\times1_{1_{C}}))\cdot(1_{\otimes}\circ
((\sigma\circ1_{1_{C}\times K_{e}})\times1_{1_{C}}))\\
& \cdot((\sigma^{-1}\circ(\sigma^{-1}\times1_{1_{C}}))\circ1_{1_{C}\times
K_{e}\times1_{C}})\cdot(t\alpha^{-1}\circ1_{1_{C}\times K_{e}\times1_{C}})\\
& \cdot((\sigma\circ(1_{1_{C}}\times\sigma))\circ1_{1_{C}\times K_{e}%
\times1_{C}})\\
& =(1_{\otimes}\circ(t\beta\times1_{1_{C}}))\cdot(\sigma^{-1}\circ
(1_{t\otimes}\circ(1_{1_{C}}\times1_{K_{e}})\times1_{1_{C}}))\cdot
(t\alpha^{-1}\circ1_{1_{C}\times K_{e}\times1_{C}})\\
& \cdot((\sigma\circ(1_{1_{C}}\times\sigma))\circ1_{1_{C}\times K_{e}%
\times1_{C}})\\
& =(\sigma^{-1}\circ(1_{P}\times1_{1_{C}}))\cdot(1_{t\otimes}\circ
(t\beta\times1_{1_{C}}))\cdot(t\alpha^{-1}\circ1_{1_{C}\times K_{e}\times
1_{C}})\\
& \cdot((\sigma\circ(1_{1_{C}}\times\sigma))\circ1_{1_{C}\times K_{e}%
\times1_{C}})\\
& =(\sigma^{-1}\circ(1_{P}\times1_{1_{C}}))\cdot(1_{\otimes}\circ(\gamma
\times1_{1_{C}})\circ1_{T_{3}})\cdot((\sigma\circ(1_{1_{C}}\times\sigma
))\circ1_{1_{C}\times K_{e}\times1_{C}})\\
& =(\sigma^{-1}\circ1_{P\times1_{C}})\cdot(1_{t\otimes}\circ(1_{1_{C}}\times
t\gamma))\cdot(\sigma\circ(1_{1_{C}}\times\sigma)\circ1_{1_{C}\times
K_{e}\times1_{C}})\\
& =(\sigma^{-1}\circ1_{1_{C}\times Q})\cdot(1_{t\otimes}\circ(1_{1_{C}}\times
t\gamma))\cdot(\sigma\circ(1_{1_{C}}\times(\sigma\circ(1_{K_{e}}\times
1_{1_{C}}))))\\
& =(1_{\otimes}\circ(1_{1_{C}}\times\lbrack t\gamma\cdot(\sigma\circ
1_{K_{e}\times1_{C}})]))\\
& =(1_{\otimes}\circ(1_{1_{C}}\times\widehat{\beta}))
\end{align*}

For the last MacLane condition we have
\begin{align*}
& \widehat{\beta}\circ1_{1_{C}\times K_{e}}\\
& =(t\gamma\circ1_{1_{C}\times K_{e}})\cdot(\sigma\circ1_{K_{e}\times1_{C}%
}\circ1_{1_{C}\times K_{e}})\\
& =(\gamma\circ1_{\tau}\circ1_{1_{C}\times K_{e}})\cdot(\sigma\circ
1_{K_{e}\times K_{e}})\\
& =(\gamma\circ1_{K_{e}\times1_{C}}\circ1_{\tau})\cdot(\sigma\circ
1_{K_{e}\times K_{e}})\\
& =(\beta\circ1_{1_{C}\times K_{e}}\circ1_{\tau})\cdot(\sigma\circ
1_{K_{e}\times K_{e}})\\
& =(t\beta\circ1_{K_{e}\times1_{C}})\cdot(\sigma\circ1_{1_{C}\times K_{e}%
}\circ1_{K_{e}\times1_{C}})\\
& =(t\beta\cdot(\sigma\circ1_{1_{C}\times K_{e}}))\circ1_{K_{e}\times1_{C}}\\
& =\widehat{\gamma}\circ1_{K_{e}\times1_{C}}%
\end{align*}
\end{proof}

Let $M_{C,\otimes,e}=\{(\alpha,\beta,\gamma)\;|\;\langle C,\otimes
,K_{e},\alpha,\beta,\gamma\rangle$ is a monoidal category$\rangle$. Then the
previous proposition show that for each natural isomorphism $\sigma
:\otimes\longrightarrow t\otimes$ we have a map
\[
T_{t}(\sigma):M_{C,\otimes,e}\longrightarrow M_{C,\otimes,e}%
\]
defined by $T_{t}(\sigma)(\alpha,\beta,\gamma)=(\widehat{\alpha}%
,\widehat{\beta},\widehat{\gamma})$. Let us next for each $\rho:\otimes
\longrightarrow\otimes$ define a map on elements in $M_{C,\otimes,e}$
\[
T_{1}(\rho)(\alpha,\beta,\gamma)=(\widetilde{\alpha},\widetilde{\beta
},\widetilde{\gamma})
\]

where we have
\begin{align*}
\widetilde{\alpha}  & =(\rho^{-1}\circ(\rho^{-1}\times1_{1_{C}}))\cdot
\alpha\cdot(\rho\circ(1_{1_{C}}\times\rho))\\
\widetilde{\beta}  & =\beta\cdot(\rho\circ1_{K_{e}\times1_{C}})\\
\widetilde{\gamma}  & =\gamma\cdot(\rho\circ1_{1_{C}\times K_{e}})
\end{align*}

For this map we have

\begin{proposition}
$T_{1}(\rho):M_{C,\otimes,e}\longrightarrow M_{C,\otimes,e}$
\end{proposition}

The proof of this proposition is similar to the one for the map $T_{1}%
(\sigma)$ and is not reproduced here.

Let
\[
G_{C,\otimes,e}=\{T_{t}(\sigma),T_{1}(\rho)\;|\;\sigma:\otimes\longrightarrow
t\otimes,\rho:\otimes\longrightarrow\otimes\text{ , }\sigma,\rho\text{ natural
isomorphisms}\}
\]
>From the construction it is evident that all maps in $G_{C,\otimes,e}$ are
bijections. The next proposition show that $G_{C,\otimes,e}$ is closed under
composition of maps.

\begin{proposition}
Let $\sigma_{1},\sigma_{2}:\otimes\longrightarrow t\otimes$ and $\rho_{1}%
,\rho_{2}:\otimes\longrightarrow\otimes$ be natural isomorphisms. Then we
have
\begin{align*}
T_{t}(\sigma_{2})\circ T_{t}(\sigma_{1})  & =T_{1}(t\sigma_{1}\cdot\sigma
_{2})\\
T_{t}(\sigma_{1})\circ T_{1}(\rho_{1})  & =T_{t}(\rho_{1}\cdot t\sigma_{1})\\
T_{1}(\rho_{1})\circ T_{t}(\sigma_{1})  & =T_{t}(\sigma_{1}\cdot\rho_{1})\\
T_{1}(\rho_{2})\circ T_{1}(\rho_{1})  & =T_{1}(\rho_{1}\cdot\rho_{2})
\end{align*}
\end{proposition}

The proof of this proposition is routine and is left out. The set
$G_{C,\otimes,e}$ is thus closed under composition and contains the identity
map $T_{1}(1_{\otimes})=1_{M_{C,\otimes,e}}$.All maps in the set
$G_{C,\otimes,e}$ are invertible by construction and $G_{C,\otimes,e}$ is
closed under the operation of taking the inverse of a map. We have
\begin{align*}
T_{1}(\rho)\circ T_{1}(\rho^{-1})  & =1_{M_{C,\otimes,e}}\\
T_{t}(\sigma)\circ T_{t}((t\sigma)^{-1})  & =1_{M_{C,\otimes,e}}%
\end{align*}

The previous propositions can now be restated in the following way.

\begin{corollary}
\label{symcor}The set $M_{C,\otimes,e}$ of monoidal structures on $C$
corresponding to a fixed product $\otimes$ and unit $e$ is invariant under the
action of the $S_{2}$-graded group $G_{C,\otimes,e}$.\label{symdefprop}
\end{corollary}

We can use the $S_{2}$-graded group $G_{C,\otimes,e}$ to give an
interpretation of the notion of a symmetric monoidal category.

\begin{proposition}
Let $\langle C,\otimes,K_{e},\alpha,\beta,\gamma,\sigma\rangle$ be a symmetric
monoidal category. Then $H=\{T_{t}(\sigma),1_{\otimes}\}$ is a $S_{2}$ graded
subgroup of $G_{C,\otimes,e}$ and $(\alpha,\beta,\gamma)\in M_{C,\otimes,e}$
is a fixpoint for the action of $H$.
\end{proposition}

This gives an interpretation of the Yang-Baxter equation and the
two unit conditions in terms of invariance with respect to the
action by the group $H$. No such interpretation appears to be
possible for the first two conditions from the definition
\ref{symcat}, of a symmetry. These two conditions appear to be of
a technical nature.

\subsection{$\sigma$-commutative comonoids in symmetric monoidal categories}

Recall that a comonoid in a monoidal category is a triple $\langle
A,\delta_{A},\epsilon_{A}\rangle$ where $A$ is a object in the category and
$\delta_{A}:A\longrightarrow A\otimes A$ \ and $\epsilon_{A}:A\longrightarrow
e$ are morphisms in the category such that the following diagrams commute%

{\tiny
\begin{diagram}
A\otimes(A\otimes A)       &      \lTo^{1_A\otimes\delta_A}    &      A\otimes
A  & \lTo^{\delta_A} &   A\\
\dTo^{\alpha_{A,A,A}}     &        \ldTo(3,3)_{\delta_A\otimes1_A}%
&                        &                        &       \\                                                                                  &                 &       \\
(A\otimes A)\otimes
A      &                                                    &                        &                         &       \\
\end{diagram}
}%

{\tiny
\begin{diagram}
e\otimes A   &   \lTo^{\epsilon_A\otimes1_A} & A\otimes A & \lTo
^{1_A\otimes\epsilon_A} & A\otimes e \\
&   \rdTo_{\beta_A}           &  \uTo_{\delta_A}  &    \ldTo_{\gamma_A}%
&\\
&                                        &     A                    &                                                  &\\
\end{diagram}
}%

The simpler structure $\langle A,\delta_{A}\rangle$ is called a cosemigroup.
The morphism $\epsilon_{A}$ is the counit for the comonoid and $\delta_{A}$ is
called the coproduct.

Before we proceed with formal developments we will first consider
some examples of these constructions. Let us first consider the
case of sets. The category $Sets$ is a monoidal category with
cartesian product, $\times$ as bifunctor. The neutral object is
the one point set $e=\{\ast\}$. The associativity constraints
$\alpha_{A,B,C}:A\times(B\times C)\longrightarrow (A\times
B)\times C$ and unit constraints $\beta_{A}:e\otimes
A\longrightarrow
A$ and $\gamma_{A}:A\otimes e\longrightarrow A$ given by%

\begin{align*}
\alpha_{A,B,C}(x,(y,z))  & =((x,y),z)\\
\beta_{A}(\ast,x)  & =x\\
\gamma_{A}(x,\ast)  & =x
\end{align*}

Finite sets offer many examples of cosemigroups. Let $A=\{a,b,c\}$ and define
a map $\delta_{A}:A\longrightarrow A\times A$ by%

\begin{align*}
\delta_{A}(a)  & =(a,a)\\
\delta_{A}(b)  & =(b,a)\\
\delta_{A}(c)  & =(a,c)
\end{align*}

A direct calculation show that $\langle A,\delta_{A}\rangle$ is a cosemigroup.
There is only one possible map $\epsilon_{A}:A\longrightarrow e $ since
$e=\{\ast\}$ is terminal is $Sets$ and this is the map $\epsilon_{A}(x)=\ast$
for all $x\in A$. But for this map we find
\[
\lbrack\beta_{A}\circ(\epsilon_{A}\otimes1_{A})\circ\delta_{A}](b)=[\beta
_{A}\circ(\epsilon_{A}\otimes1_{A})](b,a)=\beta_{A}(\ast,a)=a
\]

so $\langle A,\delta_{A},\epsilon_{A}\rangle$ is not a comonoid.

Let $A$ be any set. Define the map $\delta_{A}:A\longrightarrow A\times A$ by%

\[
\delta_{A}(x)=(x,x)
\]

This is the diagonal map in $Sets$. We then have%

\begin{align*}
\lbrack\alpha_{A,A,A}\circ(1_{A}\times\delta_{A})\circ\delta_{A}](x)  &
=[\alpha_{A,A,A}\circ(1_{A}\times\delta_{A})](x,x)=((x,x),x)\\
\lbrack(\delta_{A}\times1_{A})\circ\delta_{A}](x)  & =(\delta_{A}\times
1_{A})(x,x)=((x,x),x)
\end{align*}

so $\langle A,\delta_{A}\rangle$ is a cosemigroup. The only
possible counit satisfy
\begin{align*}
\lbrack\beta_{A}\circ(\epsilon_{A}\otimes1_{A})\circ\delta_{A}](x)  &
=[\beta_{A}\circ(\epsilon_{A}\otimes1_{A})](x,x)=\beta_{A}(\ast,x)=x\\
\lbrack\gamma_{A}\circ(1_{A}\otimes\epsilon_{A})\circ\delta_{A}](x)  &
=[\gamma_{A}\circ(1_{A}\otimes\epsilon_{A})](x,x)=\gamma_{A}(x,\ast)=x
\end{align*}

so $\langle A,\delta_{A},\epsilon_{A}\rangle$ is a comonoid. Let $\delta
_{A}:A\longrightarrow A\times A$, $\epsilon_{A}:A\longrightarrow\{\ast\}$ be
any comonoid structure on $A$. We have $\delta_{A}(a)=(f(a),g(a))$ and
$\epsilon_{A}(a)=\ast$. The first counit condition $\beta_{A}\circ
(\epsilon_{A}\times1_{A})\circ\delta_{A}=1_{A}$ gives $g(a)=a$ for all $a $.
Similarly the second counit condition gives $f(a)=a$ for all $a$. So the
previous example in fact gives the only possible comonoid structure in this
category. We will always assume that the objects in $Sets$ are comonoid with
this structure.

As our next example let us consider a pointed set. This is a set
$A$ with a chosen point $x_{0}\in A$. Define a map
$\delta_{A}:A\longrightarrow A\times
A$ by $\delta_{A}(x)=(x_{0},x)$. Then we have%

\begin{align*}
\lbrack(1_{A}\times\delta_{A})\circ\delta_{A}](x)  & =(1_{A}\times\delta
_{A})(x_{0},x)=(x_{0},x_{0},x)\\
\lbrack(\delta_{A}\times1_{A})\circ\delta_{A}](x)  & =(\delta_{A}\times
1_{A})(x_{0},x)=(x_{0},x_{0},x)
\end{align*}

so $\langle A,\delta_{A}\rangle$ is a cosemigroup. It is not a comonoid
because the only possible map $\epsilon_{A}:A\longrightarrow e$ gives
\[
\lbrack\gamma_{A}\circ(1_{A}\otimes\epsilon_{A})\circ\delta_{A}](x)=[\gamma
_{A}\circ(1_{A}\otimes\epsilon_{A})](x_{0},x)=\gamma_{A}(x_{0},\ast)=x_{0}%
\]

so if there are any elements in $A$ different from $x_{0}$ then $A$ is not a
comonoid. This construction only gives a comonoid when $A=e$. This fact is
true for any monoidal category.

Let us next consider the category $Vect_{k}$. This is the category of
vectorspaces over a field $k$ with morphisms given by linear maps. This
category is monoidal with product bifunctor given by the tensorproduct of
vectorspaces $\otimes=\otimes_{k}$. The neutral object is $k$. The
associativity constraint $\alpha$ and unit constraints $\beta$ and $\gamma$
for this case are the linear maps $\alpha_{A,B,C}:A\otimes(B\otimes
C)\longrightarrow(A\otimes B)\otimes C$, $\beta_{A}:k\otimes A\longrightarrow
A$ and $\gamma_{A}:A\otimes k\longrightarrow A$ given on generators by
\begin{align*}
\alpha_{A,B,C}(x\otimes(y\otimes z))  & =(x\otimes y)\otimes z\\
\beta_{A}(r\otimes x)  & =rx\\
\gamma_{A}(x\otimes r)  & =rx
\end{align*}

Let $A$ be any finite dimensional vectorspace in $Vect_{k}$. Let $\Omega$ be a
finite index set and let $\{a_{i}\}_{i\in\Omega}$ be a basis for $A$ indexed
by $\Omega$. Then $\{a_{i}\otimes a_{i^{\prime}}\}_{i,i^{\prime}\in\Omega}$ is
a basis for $A\otimes A$. Define a linear map $\delta_{A}:A\longrightarrow
A\otimes A$ by
\[
\delta_{A}(a_{i})=a_{i}\otimes a_{i}%
\]

Then evidently $\langle A,\delta_{A}\rangle$ is a cosemigroup. Define a linear
map $\epsilon_{A}:A\longrightarrow k$ on generators by $\epsilon_{A}%
(a_{i})=1\in k$. Then we have
\begin{align*}
\lbrack\beta_{A}\circ(\epsilon_{A}\otimes1_{A})\circ\delta_{A}](a_{i})  &
=[\beta_{A}\circ(\epsilon_{A}\otimes1_{A})](a_{i},a_{i})=\beta_{A}(1\otimes
a_{i})=a_{i}\\
\lbrack\gamma_{A}\circ(1_{A}\otimes\epsilon_{A})\circ\delta_{A}](a_{i})  &
=[\gamma_{A}\circ(1_{A}\otimes\epsilon_{A})](a_{i},a_{i})=\gamma_{A}%
(a_{i}\otimes1)=a_{i}%
\end{align*}

so $\langle A,\delta_{A},\epsilon_{A}\rangle$ is a comonoid. In contrast to
the case of $Sets$ we can have many nonisomorphic comonoid structures on a
given object in $Vect_{k}$. Let $\delta_{A}:A\longrightarrow A\otimes A$ and
$\epsilon_{A}:A\longrightarrow k$ be linear maps. We have thus
\begin{align*}
\delta_{A}(a_{i})  & =\sum_{j,k}r_{j,k}^{i}a_{j}\otimes a_{k}\\
\epsilon_{A}(a_{i})  & =q_{i}%
\end{align*}

where all indices run from $1$ to $m$, the dimension of $A$.

Then $\langle A,\delta_{A},\epsilon_{A}\rangle$ is a comonoid if
$\{r_{j,k}^{i}\}$ and $\{q_{i}\}$ are solutions of the following
system of quadratic equations.
\begin{align*}
\sum_{j}(r_{j,k}^{i}r_{l,n}^{j}-r_{l,j}^{i}r_{n,k}^{j})  & =0\text{ \ \ for
all }i,k,l,n\text{.}\\
\sum_{j}r_{j,k}^{i}q_{j}  & =\delta_{i,k}\text{ for all }i,k\text{.}\\
\sum_{j}r_{k,j}^{i}q_{j}  & =\delta_{i,k}\text{ for all }i,k\text{.}%
\end{align*}

For $m=2$ this system have four different families of solutions. One of these
families is the following
\begin{align*}
\delta_{A}(a_{1})  & =a_{1}\otimes a_{1}\\
\delta_{A}(a_{2})  & =-xa_{1}\otimes a_{1}+a_{1}\otimes a_{2}+a_{2}\otimes
a_{1}\\
q_{A}(a_{1})  & =1\\
q_{A}(a_{2})  & =x
\end{align*}
\label{cm2}

where $x$ is a arbitrary element of $k$.

Let now $G$ be a finite group and let $A=\mathcal{F}(G)$ be the vectorspace of
$k$ valued functions on $G$.

Note that since $G$ is finite we have $\mathcal{F}(G\times G)\approx
\mathcal{F}(G)\otimes_{k}\mathcal{F}(G)$. Define a linear map $\delta
_{\mathcal{F}(G)}:\mathcal{F}(G)\longrightarrow\mathcal{F}(G)\otimes
_{k}\mathcal{F}(G)$ by
\[
\delta_{\mathcal{F}(G)}(f)(x,y)=f(xy)
\]

This clearly makes $\mathcal{F}(G)$ into a cosemigroup. The linear map
$\epsilon_{\mathcal{F}(G)}:\mathcal{F}(G)\longrightarrow k$%
\[
\epsilon_{\mathcal{F}(G)}(f)=f(e)
\]

where $e\in G$ is the unit of the group $G$, makes $\langle$ $\mathcal{F}%
(G),\delta_{\mathcal{F}(G)},\epsilon_{\mathcal{F}(G)}\rangle$ into
a comonoid. Note that this conclusion depends strongly on the
identification $\mathcal{F}(G\times
G)\approx\mathcal{F}(G)\otimes_{k}\mathcal{F}(G)$. For infinite
groups this relation does not hold in general but for some
infinite groups it does. For these cases we also get comonoids.

The tensorproduct is not the only monoidal structure on $Vect_{k}$. Let
$\oplus$ be the direct sum of vectorspaces. This is a monoidal structure with
the neutral object given by the zero dimensional vectorspace $e=\{0\}$. The
maps $\alpha,\beta$ and $\gamma$ are the standard identifications used for the
direct sum. The symmetry is the linear map $\sigma(u,v)=(v,u)$. These
structures defines the structure of a symmetric monoidal category on
$Vect_{k}$. A cosemigroup is a pair $\langle A,\delta_{A}\rangle$ with
$\delta_{A}:A\rightarrow A\oplus A$ a coassociative linear map. Any such map
is determined by a pair of linear maps $f,g:A\rightarrow A$ through
$\delta_{A}(a)=(f(a),g(a))$. The coassociativity gives the following
conditions on the maps $f$ and $g$.
\begin{align*}
f\circ f  & =f\\
g\circ g  & =g\\
f\circ g  & =g\circ f
\end{align*}

So any pair of commuting projectors on $A$ define the structure of a
cosemigroup on $A$. There are thus in general many nontrivial cosemigroup
structures on a linear space. The comonoidstructure is however much more
restrictive. This is because the neutral object for $\oplus$ is also the
terminal object for the category. This means that there is only one possible
counit for any comonoid. It is straight forward to see that the counit
property for the only possible counit gives $f=g=1_{A}$. So there is only one
comonoid structure on $A$ and this is the diagonal map
\[
\delta_{A}(a)=(a,a)
\]

In all the examples we have seen that coproduct for the comonoids have been
monomorphisms. This is true in general

\begin{proposition}
Let $\langle B,\delta_{B},\epsilon_{B}\rangle$ be a comonoid. Then the
coproduct is a monomorphism.
\end{proposition}

\begin{proof}
Let $D$ be any object in $\mathcal{C}$ and let $\varphi,\psi:D\longrightarrow
B$ be two morphisms in $\mathcal{C}$ such that $\delta_{B}\circ\varphi
=\delta_{B}\circ\psi$. Then we have
\begin{align*}
\psi & =1_{B}\circ\psi\\
& =\beta_{B}\circ(\epsilon_{B}\otimes1_{B})\circ\delta_{B}\circ\psi\\
& =\beta_{B}\circ(\epsilon_{B}\otimes1_{B})\circ\delta_{B}\circ\varphi\\
& =1_{B}\circ\varphi\\
& =\varphi
\end{align*}

so $\delta_{B}$ is by definition mono.
\end{proof}

We will in general only be interested in comonoids where the
coproduct has the additional property of being commutative. Only
such comonoids carry enough structure to support a full theory of
relations. We express this property by using the symmetry
$\sigma$.

\begin{definition}
A comonoid $\langle A,\delta_{A},\epsilon_{A}\rangle$ in a symmetric monoidal
category is $\sigma$-commutative if $\sigma_{A,A}\circ\delta_{A}=\delta_{A}$.
\end{definition}

\subsection{\bigskip C-categories and M-categories}

In $Sets$ each object is a $\sigma$-commutative comonoid in one
and only one way . For the case of a general symmetric monoidal
category we have seen that objects may have several
$\sigma$-commutative comonoid structures defined on them. We need
to preserve the unique relation between objects and structures
when we generalize from $Sets$. This relation is expressed in
terms of functors and natural transformations. To any category $C$
we have associated a set of functors. These are the projection
functors $P:$ $C\times C\longrightarrow C$ and $Q:C\times
C\longrightarrow C$ ,the diagonal functor $\Delta:$
$C\longrightarrow C\times C$ defined by $\Delta(X)=(X,X)$ and the
transposition functor $\tau:C\times C\longrightarrow C\times C$ .
Let $e$ be a fixed object in the category $C$. To this object we
associate the constant functor $K_{e}:C\longrightarrow C$ .
Finally let us assume that $\otimes:$ $C\times C\longrightarrow C$
is a bifunctor and let $H=(1_{C}\times\tau
\times1_{C})\circ(\Delta\times\Delta)$. We are now ready to define
the notion of a C-category.

\begin{definition}
A C-category is a collection $\langle$ $C,\otimes,K_{e},\alpha,\beta
,\gamma,\sigma,\delta,\epsilon\rangle$ where $\langle$ $C,\otimes,K_{e}%
,\alpha,\beta,\gamma,\sigma\rangle$ is a symmetric monoidal category and where
$\delta,\epsilon$ are natural transformations
\begin{align*}
\delta & :1_{C}\longrightarrow\otimes\circ\Delta\\
\epsilon & :1_{C}\longrightarrow K_{e}%
\end{align*}

such that the following relations holds
\begin{align*}
(1_{\otimes}\circ(\delta\times1_{1_{C}})\circ1_{\Delta})\cdot\delta &
=(\alpha\circ1_{(1_{C}\times\Delta)\circ\Delta})\cdot(1_{\otimes}%
\circ(1_{1_{C}}\times\sigma)\circ1_{\Delta})\cdot\delta\\
1_{1_{C}}  & =(\beta\circ1_{\Delta})\cdot(1_{\otimes}\circ(\epsilon
\times1_{1_{C}})\circ1_{\Delta})\cdot\delta\\
1_{1_{C}}  & =(\gamma\circ1_{\Delta})\cdot(1_{\otimes}\circ(1_{1_{C}}%
\times\epsilon)\circ1_{\Delta})\cdot\delta\\
\delta & =(\sigma^{-1}\circ1_{\Delta})\cdot\delta\\
\delta\circ1_{\otimes}  & =(\alpha^{-1}\circ1_{\otimes\times1_{C}\times1_{C}%
}\circ1_{H})\cdot(1_{\otimes}\circ(\alpha\times1_{1_{C}})\circ1_{H})\\
& \cdot(1_{\otimes\circ(\otimes\times1_{C})}\circ(1_{1_{C}}\times\sigma
\times1_{1_{C}})\circ1_{\Delta\times\Delta})\\
& \cdot(1_{\otimes}\circ(\alpha^{-1}\times1_{1_{C}})\circ1_{\Delta\times
\Delta})\cdot(\alpha\circ1_{\otimes\times1_{C}\times1_{C}}\circ1_{\Delta
\times\Delta})\\
& \cdot(1_{\otimes}\circ(\delta\times\delta))\\
\epsilon\circ1_{\otimes}  & =(\beta\circ1_{1_{C}\times K_{e}})\cdot
(1_{\otimes}\circ(\epsilon\times\epsilon))
\end{align*}
\end{definition}

The four first relations ensure that for each object in $C$ there is fixed a
unique commutative comonoid structure. The last two relations say that if a
object $X$ can be decomposed as $X=A\otimes B$, then we can express the unique
comonoid structure on $X$ in terms of the comonoid structures on $A$ and $B$.
For the strict case they take the following form in terms of objects
\begin{align*}
\delta_{A\otimes B}  & =(1_{A}\otimes\sigma_{A,B}\otimes1_{B})\circ(\delta
_{A}\otimes\delta_{B})\\
\epsilon_{A\otimes B}  & =\epsilon_{A}\otimes\epsilon_{B}%
\end{align*}

A M-category is the dual of a C-category. We get its defining equations by
reversing all arrows. It is a category where for each object there is fixed a
unique monoid structure and where the monoid structure on a object of the form
$X=A\otimes B$ can be expressed in terms of the structures on $A$ and $B $.

\section{Categorical theory of relations}

In this part of the paper we use the categorical framework
described in the previous section to define a category of
relations and develop its properties. We first define the notion
of a relation and a corelation in a C-category. In a similar way
relations and corelations can be developed in a M-category. The
notions of C-categories and M-categories are dual concepts so that
any definitions made or propositions proved in one of them hold in
a dualized version in the other. Since the notion of relation and
corelation also are dual of each other it is clear that it is
enough to develop the theory of relations in C-categories. The
other cases follow by duality. We start this section by defining
relations on a object $A$ in a C-category $C$ in terms of arrows
and collect such arrows into a category of relations
$\mathcal{R}^{A}(C)$. This category of relations is then shown to
be isomorphic to the category $\mathcal{S}^{A}(C)$ of $A-A$
bicomodules in $C$. A semimonoidal structure $\boxtimes^{A}$is
introduced in this category and by isomorphism into the category
of relations. This semimonoidal structure is then used to
introduce a bifunctor $\otimes^{A}$on $\mathcal{S}^{A}(C)$ and by
isomorphism on $\mathcal{R}^{A}(C)$. This bifunctor is used to
introduce a monoidal structure on the category of relations.
Certain properties of relations like transitivity and reflexivity
are formulated in algebraic terms using the monoidal structure. In
the final sections a generalized notion of symmetry is introduced,
this notion of symmetry use in a essential way the formulation of
the Yang-Baxter equation in terms of action of a $S_{2}$ graded
group. The new notion of symmetry is then used to further
categorize properties of relations. Equivalence relations appears
as commutative and associative algebras with units.

\subsection{Relations}

Let $\langle C,\otimes,k_{e},\alpha,\beta,\gamma,\sigma,\delta,\epsilon
\rangle$ be a C-category and let $A$ be a object in $C$.

\begin{definition}
A relation on $A$ is a arrow in $C$ \ with codomain $A\otimes A$.
\end{definition}

Note that we will use the same symbol for an arrow in $C$ and the
corresponding morphism of relations. Also note that a given arrow
$f:B\longrightarrow B^{\prime}$ in $C$ can give rise to more than
one morphism of relations. This can happen because we might have
$r_{1}=r_{1}^{\prime
}\circ f$ and $r_{2}=r_{2}^{\prime}\circ f$ where $r_{1},r_{2}%
:B\longrightarrow A\otimes A$ and $r_{1}^{\prime},r_{2}^{\prime}:B^{\prime
}\longrightarrow A\otimes A$ are two pairs of relations on $A$. \ In this
sense we can write $1_{r}=1_{B}$ where $B$ is the domain of the arrow $r$. Let
us now consider a few examples of this construction.

Let us first consider the case of $Sets$. This is a C-category with
$\delta_{X}(x)=(x,x)$ and $\epsilon_{X}(x)=\ast$ for all objects $X$ $\in C$.
Let $A$ and $B$ be sets and let $r:B\longrightarrow A\times A$ be a map of
sets. We can write $r(b)=(f(b),g(b))$. We have
\begin{align*}
& [(r\times r)\circ\delta_{B}](b)\\
& =(f(b),f(b),g(b),g(b))\\
& =(\delta_{A}\times\delta_{A})(f(b),g(b))\\
& =[\delta_{A\times A}\circ r](b)
\end{align*}

so $r$ is a arrow in the C-category $Sets$ and is therefore a
relation in $Sets$ in our sense. A relation on $A$ in the usual
sense is a subset of $A\times A$. This is equivalent to assuming
that the map $r$ is a monomorphism. In general the map $r$ assign
to each element in $B$ a source and a target. Several elements in
$B$ can be assigned the same source and target. In fact we observe
that in $Sets$ a relation in our sense is the same as a directed
labelled graph.

Let us next consider the C-category $Vect_{k}$ with direct sum as monoidal
structure and $\delta$ and $\epsilon$ defined as for $Sets$. A relation on a
linear space $A$ is any linear map $r:B\longrightarrow A\oplus A$. Let
$L:A\longrightarrow A$ be a endomorphism on $A$. Let $B=A$ and define
$r:B\longrightarrow A\oplus A$ by
\[
r(a)=(a,L(a))
\]

Then $r$ is a linear map and therefore defines a relation on $A$ in our sense.
Note that the image of $A$ under $r$ is by definition the graph of the linear
map $L$. More generally, let $L$ be a linear subspace of $A\oplus A$. Let
$B=L$ and $r:B\longrightarrow A\oplus A$ the inclusion map. Then $r$ is
evidently a relation on $A$. In general a relation on $A$ is like a graph,
where the set of vertices and the set of labels have a vectorspace structure
and the source and target maps respect these structures.

As with any categorical concept the notion of a relation has a dual.

\begin{definition}
A corelation on a $A$ is a arrow in $C$ with domain $A\otimes A$.
\end{definition}

Let $r:S\longrightarrow\Omega\times\Omega$ be a relation on $\Omega$ in
$Sets$. We assume now that the sets $S$ and $\Omega$ are finite. The algebraic
description of the sets $S$ and $\Omega$ are given by the space of $k$ valued
functions $B=\mathcal{F}(S)$ on $S$ and $A=\mathcal{F}(\Omega) $ on $\Omega$.
Let $c$ defined by $c(f\otimes g)(x)=(f\otimes g)(r(x))$. Then $c:A\otimes
A\longrightarrow B$ is a linear map and by duality a morphism of the induced
algebra structures on $B$ and $A\otimes A$.

Therefore the algebraic image of the relation $r$ in $Sets$ is a corelation
$c$ in $Vect_{k}$. This example show that corelations arise naturally by
algebraization of relations in $Sets$. \ Note that in general a corelation
$c:A\otimes A\longrightarrow B$ in $Vect_{k}$ with the tensor product as
monoidal structure is in algebra usually called a $A\otimes A$ algebra.

\subsection{Categories of relations}

Let $r:B\longrightarrow A\otimes A$ and $r^{\prime}:B^{\prime}\longrightarrow
A\otimes A$ be two relations on $A$. A morphism $f:r\longrightarrow r^{\prime
}$ is a arrow $f:B\longrightarrow B^{\prime}$ in $C$ such that the following
diagram commute%

{\tiny
\begin{diagram}
B             &                       &          \rTo^{f}%
&                      &          B'            \\
&    \rdTo_{r}  &                           & \ldTo_{r'}%
&                         \\
&                       &         A\otimes
A &                      &                         \\
\end{diagram}
}%

Let $\mathcal{R}^{A}(C)$ be the category of relations on $A$. \
This is a category whose objects are relations and morphisms are
morphisms of relations as just defined. It is evident that to each
diagram in $\mathcal{R}^{A}(C)$ there is a corresponding diagram
of arrows in $C$ and commutativity of diagrams in
$\mathcal{R}^{A}(C)$ follows from commutativity of the
corresponding diagrams in $C$. For now there is no restriction on
the object $A$ or the arrows that are relations on $A$. We will
introduce further restrictions as we develop the properties of the
category of relations.

Morphisms of corelations are defined by dualizing the corresponding diagrams
for morphisms of relations. Corelations and morphisms of corelations form the
category of corelations on $A,$ $\mathcal{R}_{A}(C)$.

We will now proceed to develop some formal properties of the
category $\mathcal{R}^{A}(C)$. The corresponding dualized
properties holds for the category $\mathcal{R}_{A}(C)$.

Let $r$ be a object in $\mathcal{R}^{A}(C)$ with domain $B$. Define two arrows
$\delta^{l}:B\longrightarrow A\otimes B$ and $\delta^{r}:B\longrightarrow
B\otimes A$ in $C$ by
\begin{align*}
\delta^{l}  & =(\gamma_{A}\otimes1_{B})\circ((1_{A}\otimes\epsilon_{A}%
)\otimes1_{B})\circ(r\otimes1_{B})\circ\delta_{B}\\
\delta^{r}  & =(1_{B}\otimes\beta_{A})\circ(1_{B}\otimes(\epsilon_{A}%
\otimes1_{A}))\circ(1_{B}\otimes r)\circ\delta_{B}%
\end{align*}

Define $\theta^{l}:B\longrightarrow(A\otimes A)\otimes B$ and $\theta
^{r}:B\longrightarrow B\otimes(A\otimes A)$ by
\begin{align*}
\theta^{l}  & =(r\otimes1_{B})\circ\delta_{B}\\
\theta^{r}  & =(1_{B}\otimes r)\circ\delta_{B}%
\end{align*}

We first prove the identities

\begin{lemma}%
\begin{align*}
(\delta_{A\otimes A}\otimes1_{B})\circ\theta^{l}  & =\alpha_{A\otimes
A,A\otimes A,B}\circ(1_{A\otimes A}\otimes\theta^{l})\circ\theta^{l}\\
(1_{B}\otimes\delta_{A\otimes A})\circ\theta^{r}  & =\alpha_{B,A\otimes
A,A\otimes A}\circ(\theta^{r}\otimes1_{A\otimes A})\circ\theta^{r}\\
(\theta^{l}\otimes1_{A\otimes A})\circ\theta^{r}  & =\alpha_{A\otimes
A,B,A\otimes A}\circ(1_{A\otimes A}\otimes\theta^{r})\circ\theta^{l}%
\end{align*}
\end{lemma}

\begin{proof}
Since $r$ is a morphism in $C$ we have the diagram%

{\tiny
\begin{diagram}
B \times B               &   \rTo^{r\otimes r}            &   A\otimes
A\otimes A \otimes A    \\
\uTo^{\delta_B}    &                                                   & \uTo
_{\delta_{A\otimes A}}             \\
B                            &      \rTo_{r}%
&  A\otimes A                                    \\
\end{diagram}
}%

But then we have
\begin{align*}
& \alpha_{A\otimes A,A\otimes A,B}\circ(1_{A\otimes A}\otimes\theta^{l}%
)\circ\theta^{l}\\
& =\alpha_{A\otimes A,A\otimes A,B}\circ(1_{A\otimes A}\otimes(r\otimes
1_{B}))\circ(1_{A\otimes A}\otimes\delta_{B})\circ(r\otimes1_{B})\circ
\delta_{B}\\
& =\alpha_{A\otimes A,A\otimes A,B}\circ(r\otimes(r\otimes1_{B}))\circ
(1_{B}\otimes\delta_{B})\circ\delta_{B}\\
& =\alpha_{A\otimes A,A\otimes A,B}\circ(r\otimes(r\otimes1_{B}))\circ
\alpha_{B,B,B}^{-1}\circ(\delta_{B}\otimes1_{B})\circ\delta_{B}\\
& =((r\otimes r)\circ\delta_{B}\otimes1_{B})\circ\delta_{B}\\
& =(\delta_{A\otimes A}\otimes1_{B})\circ(r\otimes1_{B})\circ\delta_{B}%
\end{align*}

The proof of the second relation proceeds in a similar way. For the third we
have
\begin{align*}
& (\theta^{l}\otimes1_{A\otimes A})\circ\theta^{r}\\
& =((r\otimes1_{B})\otimes1_{A\otimes A})\circ(\delta_{B}\otimes1_{A\otimes
A})\circ(1_{B}\otimes r)\circ\delta_{B}\\
& =((r\otimes1_{B})\otimes r)\circ(\delta_{B}\otimes1_{B})\circ\delta_{B}\\
& =((r\otimes1_{B})\otimes r)\circ\alpha_{B,B,B}\circ(1_{B}\otimes\delta
_{B})\circ\delta_{B}\\
& =\alpha_{A\otimes A,B,A\otimes A}\circ(r\otimes(1_{B}\otimes r))\circ
(1_{B}\otimes\delta_{B})\circ\delta_{B}\\
& =\alpha_{A\otimes A,B,A\otimes A}\circ(r\otimes\theta^{r})\circ\delta_{B}%
\end{align*}
\end{proof}

we can now prove the following

\begin{proposition}
Let $r$ be a relation. Then the following diagrams commute

\label{bi1}%
{\tiny
\begin{equation*}
\begin{aligned}
\begin{diagram}
&                                  & A\otimes B   & \rTo^{\delta_A\otimes1_B}%
&  (A\otimes A)\otimes B   \\
& \ruTo^{\delta^l}%
&                     &                                                  &                                       \\
B                &                                  &                     &                                                  &   \uTo
_{\alpha_{A,A,B}}  \\
&  \rdTo^{\delta^l   }%
&                     &                                                  &                                        \\
&                                &  A\otimes B  &  \rTo_{1_A\otimes\delta^l}%
&   A\otimes(A\otimes B)     \\
\end{diagram}
\end{aligned}
\;
\;
\;
\begin{aligned}
\begin{diagram}
e\otimes B              &\lTo^{\epsilon_A\otimes1_B}    &   A\otimes
B            \\
&  \rdTo_{\beta_B}                      &    \uTo_{\delta^l }  \\
&                                                   &        B                          \\
\end{diagram}
\end{aligned}
\end{equation*}
}

\label{bi2}%
{\tiny
\begin{equation*}
\begin{aligned}
\begin{diagram}
&                                  & B\otimes A   & \rTo^{1_B\otimes\delta_A}%
&  B\otimes(A\otimes A)   \\
& \ruTo^{\delta^r}%
&                     &                                                  &                                       \\
B                &                                  &                     &                                                  &   \dTo
_{\alpha_{B,A,A}}  \\
&  \rdTo^{\delta^r   }%
&                     &                                                  &                                        \\
&                                &  B\otimes A  &  \rTo_{\delta^r\otimes1_A}%
&   (B\otimes A)\otimes A    \\
\end{diagram}
\end{aligned}
\;
\;
\;
\begin{aligned}
\begin{diagram}
B\otimes A              &\rTo^{1_B\otimes\epsilon_A}    &   B\otimes
e           \\
\uTo^{\delta^r}&  \ldTo_{\gamma_B}%
&                               \\
B                            &                                                   &                             \\
\end{diagram}
\end{aligned}
\end{equation*}
}

\label{bi3}%
{\tiny
\begin{diagram}
&                                  & A\otimes B   & \rTo^{1_A\otimes\delta^r}%
& A\otimes(B\otimes A)  \\
& \ruTo^{\delta^l}%
&                     &                                                  &                                       \\
B                &                                  &                     &                                                  &   \dTo
_{\alpha_{A,B,A}}  \\
&  \rdTo^{\delta^r   }%
&                     &                                                  &                                        \\
&                                &  B\otimes A  &  \rTo_{\delta^l\otimes1_A}%
&   (A\otimes B)\otimes A   \\
\end{diagram}
}
\end{proposition}

\begin{proof}
Using the lemma and the naturality of $\alpha$ and $\delta$ we have
\begin{align*}
& \alpha_{A,A,B}\circ(1_{A}\otimes\delta^{l})\circ\delta^{l}\\
& =\alpha_{A,A,B}\circ(1_{A}\otimes(\gamma_{A}\otimes1_{B}))\circ(1_{A}%
\otimes((1_{A}\otimes\epsilon_{A})\otimes1_{B}))\\
& \circ(1_{A}\otimes\theta^{l})\circ(\gamma_{A}\otimes1_{B})\circ
((1_{A}\otimes\epsilon_{A})\otimes1_{B})\circ\theta^{l}\\
& =(((\gamma_{A}\circ(1_{A}\otimes\epsilon_{A}))\otimes(\gamma_{A}\circ
(1_{A}\otimes\epsilon_{A})))\otimes1_{B})\\
& \circ\alpha_{A\otimes A,A\otimes A,B}\circ(1_{A\otimes A}\otimes\theta
^{l})\circ\theta^{l}\\
& =(((\gamma_{A}\circ(1_{A}\otimes\epsilon_{A}))\otimes(\gamma_{A}\circ
(1_{A}\otimes\epsilon_{A})))\otimes1_{B})\\
& \circ(\delta_{A\otimes A}\otimes1_{B})\circ\theta^{l}\\
& =(\delta_{A}\otimes1_{B})\circ((\gamma_{A}\circ(1_{A}\otimes\epsilon
_{A}))\otimes1_{B})\circ\theta^{l}\\
& =(\delta_{A}\otimes1_{B})\circ\delta^{l}\\
& =(\delta_{A}\otimes1_{B})\circ(1_{A}\otimes\epsilon_{A}\otimes1_{B}%
)\circ\theta^{l}\\
& =(\delta_{A}\otimes1_{B})\circ\delta^{l}%
\end{align*}

Since $\epsilon$ is natural and $r$ a morphism in $C$ we have the identities
\begin{align*}
\epsilon_{A\otimes e}  & =\epsilon_{A}\circ\gamma_{A}\\
\epsilon_{A\otimes A}  & =\epsilon_{A\otimes e}\circ(1_{A}\otimes\epsilon
_{A})\\
\epsilon_{B}  & =\epsilon_{A\otimes A}\circ r
\end{align*}

But then we have
\begin{align*}
& \beta_{B}\circ(\epsilon_{A}\otimes1_{B})\circ\delta^{l}\\
& =\beta_{B}\circ(\epsilon_{A}\otimes1_{B})\circ(\gamma_{A}\otimes1_{B}%
)\circ((1_{A}\otimes\epsilon_{A})\otimes1_{B})\circ(r\otimes1_{B})\circ
\delta_{B}\\
& =\beta_{B}\circ(\epsilon_{A\otimes e}\otimes1_{B})\circ((1_{A}%
\otimes\epsilon_{A})\otimes1_{B})\circ(r\otimes1_{B})\circ\delta_{B}\\
& =\beta_{B}\circ(\epsilon_{A\otimes A}\otimes1_{B})\circ(r\otimes1_{B}%
)\circ\delta_{B}\\
& =\beta_{B}\circ(\epsilon_{B}\otimes1_{B})\circ\delta_{B}\\
& =1_{B}%
\end{align*}

so the first pair of diagrams are commutative. The proof of the commutativity
of the second pair of diagrams is similar. For the last diagram we have
\begin{align*}
& \alpha_{A,B,A}\circ(1_{A}\otimes\delta^{r})\circ\delta^{l}\\
& \alpha_{A,B,A}\circ(1_{A}\otimes(1_{B}\otimes\beta_{A}))\circ(1_{A}%
\otimes(1_{B}\otimes(\epsilon_{A}\otimes1_{A})))\\
& \circ(1_{A}\otimes\theta^{r})\circ(\gamma_{A}\otimes1_{B})\circ
((1_{A}\otimes\epsilon_{A})\otimes1_{B})\circ\theta^{l}\\
& =((1_{A}\otimes1_{B})\otimes\beta_{A})\circ((1_{A}\otimes1_{B}%
)\otimes(\epsilon_{A}\otimes1_{A}))\\
& \circ((\gamma_{A}\otimes1_{B})\otimes1_{A\otimes A})\circ(((1_{A}%
\otimes\epsilon_{A})\otimes1_{B})\otimes1_{A\otimes A})\\
& \circ\alpha_{A\otimes A,B,A\otimes A}\circ(1_{A\otimes A}\otimes\theta
^{r})\circ\theta^{l}\\
& =(((\gamma_{A}\circ(1_{A}\otimes\epsilon_{A}))\otimes1_{B})\otimes(\beta
_{A}\circ(\epsilon_{A}\otimes1_{A})))\\
& \circ(\theta^{l}\otimes1_{A\otimes A})\circ\theta^{r}\\
& =(\delta^{l}\otimes1_{A})\circ(1_{B}\otimes(\beta_{A}\circ(\epsilon
_{A}\otimes1_{A})))\circ\theta^{r}\\
& =(\delta^{l}\otimes1_{A})\circ\delta^{r}%
\end{align*}

so this diagram is also commutative.
\end{proof}

The previous proposition show that the pair $\{\delta^{l},\delta^{r}\}$ define
the structure of a $A-A$ bicomodule on $B$.

\begin{definition}
Let $\delta^{l}:B\longrightarrow A\otimes B,\delta^{r}:B\longrightarrow
B\otimes A$ be arrows in $C$. Then $\{\delta^{l},\delta^{r}\}$ define the
structure of a $A-A$ bicomodule on $B$ if the diagrams in proposition
\ref{bi1} commute
\end{definition}

We call the object $B$ in the previous definition the underlying object for
the $A-A$ bicomodule $\delta=\{\delta^{l},\delta^{r}\}$.

Let now $\delta$ and $\gamma$ be $A-A$ bicomodules with underlying objects $B$
and $E$. A morphism $f:\delta\longrightarrow\gamma$ of $A-A$ bicomodules is a
arrow $f:B\longrightarrow E$ in $C$ such that the following diagrams commute%

{\tiny
\begin{equation*}
\begin{aligned}
\begin{diagram}
A\otimes B                        &  \rTo^{1_A\otimes f}%
&  A\otimes E                    \\
\uTo^{\delta^l}%
&                                              &   \uTo_{\gamma^l}        \\
B                                       &   \rTo_{f}%
&     E                                 \\
\end{diagram}
\end{aligned}
\;
\;
\;
\begin{aligned}
\begin{diagram}
B\otimes A                      &  \rTo^{f\otimes1_A}            &  E\otimes
A                  \\
\uTo^{\delta^r}%
&                                              &   \uTo_{\gamma^r}        \\
B                                       &   \rTo_{f}%
&     E                                 \\
\end{diagram}
\end{aligned}
\end{equation*}
}

We now form a new category where objects are $A-A$ bicomodules and where
morphisms are morphisms of $A-A$ bicomodules. Let this category be named
$\mathcal{S}^{A}(C)$.

\subsection{Relations in terms of $A-A$ bicomodules.}

\ To each object $r$ in $\mathcal{R}^{A}(C)$ there corresponds a object
$\delta$ in $\mathcal{S}^{A}(C)$. For morphisms of relations we have the following.

\begin{proposition}
Let $r$ and $s$ be two relations with domains $B$ and $E$ and let
$f:r\longrightarrow s$ be a morphism of relations. Let $\delta$ and $\gamma$
be the objects in $\mathcal{S}^{A}(C)$ corresponding to $r$ and $s$. Then the
corresponding arrow $f:B\longrightarrow E$ in $C$ defines a morphism
$f:\delta\longrightarrow\gamma$ in $\mathcal{S}^{A}(C)$.
\end{proposition}

\begin{proof}
\bigskip Since $f$ is a morphism in $C$ and also a morphism from $r$ to $s$ we
have the following commutative diagrams%

{\tiny
\begin{equation*}
\begin{aligned}
\begin{diagram}
B             &                       &          \rTo^{f}%
&                      &          E            \\
&    \rdTo_{r}  &                           & \ldTo_{s}%
&                         \\
&                       &         A\otimes
A &                      &                         \\
\end{diagram}
\end{aligned}
\;
\;
\;
\begin{aligned}
\begin{diagram}
B\otimes B            &   \rTo^{f\otimes f}                         & E\otimes
E           \\
\uTo^{\delta_B}%
&                                                       & \uTo_{\delta_E}%
\\
B                           &   \rTo_{f}%
&     E                      \\
\end{diagram}
\end{aligned}
\end{equation*}
}

Using these identities we have e then have%

\begin{align*}
& (1_{A}\otimes f)\circ\delta^{l}\\
& =(1_{A}\otimes f)\circ(\gamma_{A}\otimes1_{B})\circ((1_{A}\otimes
\epsilon_{A})\otimes1_{B})\circ(r\otimes1_{B})\circ\delta_{B}\\
& =(\gamma_{A}\otimes1_{E})\circ((1_{A}\otimes1_{e})\otimes f)\circ
(s\otimes1_{B})\circ(f\otimes1_{B})\circ\delta_{B}\\
& =(\gamma_{A}\otimes1_{E})\circ((1_{A}\otimes\epsilon_{A})\otimes1_{E}%
)\circ(s\otimes1_{B})\circ(f\otimes f)\circ\delta_{B}\\
& =(\gamma_{A}\otimes1_{E})\circ((1_{A}\otimes\epsilon_{A})\otimes1_{E}%
)\circ(s\otimes1_{B})\circ\delta_{E}\circ f\\
& =\gamma^{l}\circ f
\end{align*}

In a similar way we prove the identity $(f\otimes1_{A})\circ\delta^{r}%
=\gamma^{r}\circ f$.
\end{proof}

The previous definition show that we have a well defined functor
$\Phi:\mathcal{R}^{A}(C)\longrightarrow\mathcal{S}^{A}(C)$, where $\Phi(r)$ is
the $A-A$ bicomodule corresponding to $r$ and where $\Phi(f)=f$. We will next
construct a functor from $\mathcal{S}^{A}(C)$ to $\mathcal{R}^{A}(C)$.

Let $\delta$ be a object in $\mathcal{S}^{A}(C)$. define a morphism
$r:B\longrightarrow A\otimes A$ by
\[
r=(1_{A}\otimes\beta_{A})\circ(1_{A}\otimes(\epsilon_{B}\otimes1_{A}%
))\circ(1_{A}\otimes\delta^{r})\circ\delta^{l}%
\]

We have proved that $\beta_{A}$ and $\epsilon_{B}$ are arrows in $C$ and have
therefore the following result.

\begin{proposition}
$r$ is a object in $\mathcal{R}^{A}(C)$.
\end{proposition}

Using this result we can define a map of objects $\Psi:\mathcal{S}%
^{A}(C)\longrightarrow\mathcal{R}^{A}(C)$ by $\Psi(\delta)=r$. For morphisms
in $\mathcal{S}^{A}(C)$ we have

\begin{proposition}
Let $\delta$, and $\gamma$ be two objects in $\mathcal{S}^{A}(C)$ and let
$f:\delta\longrightarrow\gamma$ be a morphism. Then the corresponding arrow in
$C$ defines a morphism of the objects $r=\Psi(\delta)$ and $s=\Psi(\gamma)$ in
$\mathcal{R}^{A}(C)$.
\end{proposition}

\begin{proof}
Let the domains of $r$ and $s$ be $B$ and $E$. We have the following
commutative diagrams%

{\tiny
\begin{equation*}
\begin{aligned}
\begin{diagram}
A\otimes B                        &  \rTo^{1_A\otimes f}%
&  A\otimes E                   \\
\uTo^{\delta^l}%
&                                              &   \uTo_{\gamma^l}        \\
B                                       &   \rTo_{f}%
&     E                                \\
\end{diagram}
\end{aligned}
\;
\;
\;
\begin{aligned}
\begin{diagram}
B\otimes A                      &  \rTo^{f\otimes1_A}            &  E\otimes
A                  \\
\uTo^{\delta^r}%
&                                              &   \uTo_{\gamma^r}        \\
B                                       &   \rTo_{f}%
&     E                                \\
\end{diagram}
\end{aligned}
\end{equation*}
}

{\tiny
\begin{diagram}
B                                   &  \rTo^{f}                     &   E\\
\dTo^{\epsilon_B}        &  \ldTo_{\epsilon_E}   &       \\
e                                    &                                   &       \\
\end{diagram}
}%

But then we have
\begin{align*}
& s\circ f\\
& =(1_{A}\otimes\beta_{A})\circ(1_{A}\otimes(\epsilon_{E}\otimes1_{A}%
))\circ(1_{A}\otimes\gamma^{r})\circ\gamma^{l}\circ f\\
& =(1_{A}\otimes\beta_{A})\circ(1_{A}\otimes(\epsilon_{E}\otimes1_{A}%
))\circ(1_{A}\otimes\gamma^{r})\circ(1_{A}\otimes f)\circ\delta^{l}\\
& =(1_{A}\otimes\beta_{A})\circ(1_{A}\otimes(\epsilon_{E}\otimes1_{A}%
))\circ(1_{A}\otimes(f\otimes1_{A}))\circ(1_{A}\otimes\delta^{r})\circ
\delta^{l}\\
& =(1_{A}\otimes\beta_{A})\circ(1_{A}\otimes(\epsilon_{E}\circ f\otimes
1_{A}))\circ(1_{A}\otimes\delta^{r})\circ\delta^{l}\\
& =(1_{A}\otimes\beta_{A})\circ(1_{A}\otimes(\epsilon_{B}\otimes1_{A}%
))\circ(1_{A}\otimes\delta^{r})\circ\delta^{l}\\
& =r
\end{align*}

so $f$ is a morphism of relations.
\end{proof}

We use this result to extend $\Psi$ to a functor from $\mathcal{S}^{A}(C)$ to
$\mathcal{R}^{A}(C)$ by defining $\Psi(f)=f$. We will now show that
$\mathcal{R}^{A}(C)$ and $\mathcal{S}^{A}(C)$ are isomorphic categories. We
need the following lemma

\begin{lemma}
\label{delta_b morphism}%
\begin{align*}
(\gamma_{A}\otimes\beta_{A})\circ((1_{A}\otimes\epsilon_{A})\otimes
(\epsilon_{A}\otimes1_{A}))\circ\delta_{A\otimes A}  & =1_{A\otimes A}\\
(\delta^{l}\otimes1_{B})\circ\delta_{B}  & =\alpha_{A,B,B}\circ(1_{A}%
\otimes\delta_{B})\circ\delta^{l}\\
(1_{B}\otimes\delta^{r})\circ\delta_{B}  & =\alpha_{B,B,A}^{-1}\circ
(\delta_{B}\otimes1_{A})\circ\delta^{r}%
\end{align*}
\end{lemma}

\begin{proof}
for the first part of the lemma we have
\begin{align*}
& (\gamma_{A}\otimes\beta_{A})\circ((1_{A}\otimes\epsilon_{A})\otimes
(\epsilon_{A}\otimes1_{A}))\circ\delta_{A\otimes A}\\
& =(\gamma_{A}\otimes\beta_{A})\circ((1_{A}\otimes\epsilon_{A})\otimes
(\epsilon_{A}\otimes1_{A}))\circ\alpha_{A\otimes A,A,A}^{-1}\\
& \circ(\alpha_{A,A,A}\otimes1_{A})\circ((1_{A}\otimes\sigma_{A,A}%
)\otimes1_{A})\circ(\alpha_{A,A,A}^{-1}\otimes1_{A})\\
& \circ\alpha_{A\otimes A,A,A}\circ(\delta_{A}\otimes\delta_{A})\\
& =(1_{A}\otimes\beta_{A})\circ(\gamma_{A}\otimes(1_{e}\otimes1_{A}%
))\circ\alpha_{A\otimes e,e,A}^{-1}\circ(\alpha_{A,e,e}\otimes1_{A})\\
& \circ((1_{A}\otimes\sigma_{e,e})\otimes1_{A})\circ(\alpha_{A,e,e}%
^{-1}\otimes1_{A})\circ\alpha_{A\otimes e,A,A}\\
& \circ((1_{A}\otimes\epsilon_{A})\otimes(\epsilon_{A}\otimes1_{A}%
))\circ(\delta_{A}\otimes\delta_{A})\\
& =(1_{A}\otimes\beta_{A})\circ(\gamma_{A}\otimes(1_{e}\otimes1_{A}%
))\circ\alpha_{A\otimes e,e,A}^{-1}\circ(\alpha_{A,e,e}\otimes1_{A})\\
& \circ((1_{A}\otimes\sigma_{e,e})\otimes1_{A})\circ(\alpha_{A,e,e}%
^{-1}\otimes1_{A})\circ\alpha_{A\otimes e,A,A}\\
& \circ(\gamma_{A}^{-1}\otimes\beta_{A}^{-1})\\
& =(1_{A}\otimes\beta_{A})\circ\alpha_{A,e,A}^{-1}\circ((1_{A}\otimes\beta
_{e})\otimes1_{A})\circ((1_{A}\otimes\sigma_{e,e})\otimes1_{A})\\
& \circ(\alpha_{A,e,e}^{-1}\otimes1_{A})\circ\alpha_{A\otimes e,e,A}%
\circ(\gamma_{A}^{-1}\otimes\beta_{A}^{-1})\\
& =(1_{A}\otimes\beta_{A})\circ\alpha_{A,e,A}^{-1}\circ((1_{A}\otimes
\gamma_{e})\otimes1_{A})\circ\\
& \circ(\alpha_{A,e,e}^{-1}\otimes1_{A})\circ\alpha_{A\otimes e,e,A}%
\circ(\gamma_{A}^{-1}\otimes\beta_{A}^{-1})\\
& =(1_{A}\otimes\beta_{A})\circ(1_{A}\otimes(\gamma_{e}\otimes1_{A}))\circ\\
& \circ\alpha_{A,e\otimes e,A}^{-1}\circ(\alpha_{A,e,e}^{-1}\otimes1_{A}%
)\circ\alpha_{A\otimes e,e,A}\circ(\gamma_{A}^{-1}\otimes\beta_{A}^{-1})\\
& =(1_{A}\otimes\beta_{A})\circ(1_{A}\otimes(\gamma_{e}\otimes1_{A}%
))\circ(1_{A}\otimes\alpha_{e,e,A})\\
& \circ\alpha_{A,e,e\otimes A}^{-1}\circ((1_{A}\otimes1_{e})\otimes\beta
_{A}^{-1})\circ(\gamma_{A}^{-1}\otimes1_{A})\\
& =(1_{A}\otimes\beta_{A})\circ(1_{A}\otimes(1_{e}\otimes\beta_{A}%
))\circ(1_{A}\otimes(1_{e}\otimes\beta_{A}^{-1}))\\
& \circ\alpha_{A,e,A}^{-1}\circ(\gamma_{A}^{-1}\otimes1_{A})\\
& =(1_{A}\otimes\beta_{A})\circ\alpha_{A,e,A}^{-1}\circ(\gamma_{A}^{-1}%
\otimes1_{A})\\
& =(1_{A}\otimes\beta_{A})\circ(1_{A}\otimes\beta_{A}^{-1})\\
& =1_{A\otimes A}%
\end{align*}

for the second part of the lemma we have
\begin{align*}
& (\delta^{l}\otimes1_{B})\circ\delta_{B}\\
& =((\gamma_{A}\otimes1_{B})\otimes1_{B})\circ(((1_{A}\otimes\epsilon
_{A})\otimes1_{B})\otimes1_{B})\circ((r\otimes1_{B})\otimes1_{B})\\
& \circ(\delta_{B}\otimes1_{B})\circ\delta_{B}\\
& =((\gamma_{A}\otimes1_{B})\otimes1_{B})\circ(((1_{A}\otimes\epsilon
_{A})\otimes1_{B})\otimes1_{B})\circ((r\otimes1_{B})\otimes1_{B})\\
& \circ\alpha_{B,B,B}\circ(1_{B}\otimes\delta_{B})\circ\delta_{B}\\
& =\alpha_{A,B,B}\circ(\gamma_{A}\otimes1_{B\otimes B})\circ((1_{A}%
\otimes\epsilon_{A})\otimes1_{B\otimes B})\circ(r\otimes1_{B\otimes B}%
)\circ(1_{B}\otimes\delta_{B})\circ\delta_{B}\\
& =(1_{B}\otimes\delta_{B})\circ(\gamma_{A}\otimes1_{B})\circ((1_{A}%
\otimes\epsilon_{A})\otimes1_{B})\circ(r\otimes1_{B})\circ\delta_{B}\\
& =(1_{B}\otimes\delta_{B})\circ\delta^{l}%
\end{align*}

The proof of the third part of the lemma is similar to the second part.
\end{proof}

We now use the lemma to prove the following theorem

\begin{theorem}
The functor $\Phi:\mathcal{R}^{A}(C)\longrightarrow\mathcal{S}^{A}(C)$ is
invertible with inverse $\Psi:\mathcal{S}^{A}(C)\longrightarrow\mathcal{R}%
^{A}(C)$.
\end{theorem}

\begin{proof}
We only need to prove that $\Phi$ is bijective with inverse $\Psi$ on objects
since $\Psi$ is obviously the inverse of $\Phi$ on morphisms.

Let $r$ be a object in $\mathcal{R}^{A}(C)$ with domain $B$. Using lemma
\ref{delta_b morphism} we have
\begin{align*}
& (\Psi\circ\Phi)(r)\\
& =(1_{A}\otimes\beta_{A})\circ(1_{A}\otimes(\epsilon_{B}\otimes1_{A}%
))\circ(1_{A}\otimes(\Phi(r))^{r})\circ(\Phi(r))^{l}\\
& =(1_{A}\otimes\beta_{A})\circ(1_{A}\otimes(\epsilon_{B}\otimes1_{A}%
))\circ(1_{A}\otimes(1_{B}\otimes\beta_{A}))\\
& \circ(1_{A}\otimes(1_{B}\otimes(\epsilon_{A}\otimes1_{A})))\circ
(1_{A}\otimes(1_{B}\otimes r))\circ(1_{A}\otimes\delta_{B})\\
& \circ(\gamma_{A}\otimes1_{B})\circ((1_{A}\otimes\epsilon_{A})\otimes
1_{B})\circ(r\otimes1_{B})\circ\delta_{B}\\
& =(1_{A}\otimes\beta_{A})\circ(1_{A}\otimes(1_{e}\otimes\beta_{A}%
))\circ(1_{A}\otimes(1_{e}\otimes(\epsilon_{A}\otimes1_{A})))\\
& \circ(1_{A}\otimes(1_{e}\otimes r))\circ(1_{A}\otimes(\epsilon_{B}%
\otimes1_{B})\circ\delta_{B})\circ(\gamma_{A}\otimes1_{B})\\
& \circ((1_{A}\otimes\epsilon_{A})\otimes1_{B})\circ(r\otimes1_{B})\circ
\delta_{B}\\
& =(1_{A}\otimes\beta_{A})\circ(1_{A}\otimes(1_{e}\otimes\beta_{A}%
))\circ(1_{A}\otimes(1_{e}\otimes(\epsilon_{A}\otimes1_{A})))\\
& \circ(1_{A}\otimes(1_{e}\otimes r))\circ(1_{A}\otimes\beta_{B}^{-1}%
)\circ(\gamma_{A}\otimes1_{B})\\
& \circ((1_{A}\otimes\epsilon_{A})\otimes1_{B})\circ(r\otimes1_{B})\circ
\delta_{B}\\
& =(1_{A}\otimes\beta_{A})\circ(1_{A}\otimes(1_{e}\otimes\beta_{A}%
))\circ(1_{A}\otimes(1_{e}\otimes(\epsilon_{A}\otimes1_{A})))\\
& \circ(1_{A}\otimes\beta_{A\otimes A}^{-1})\circ(\gamma_{A}\otimes1_{A\otimes
A})\circ((1_{A}\otimes\epsilon_{A})\otimes1_{A\otimes A})\circ(r\otimes
r)\circ\delta_{B}\\
& =(1_{A}\otimes\beta_{A}\circ(1_{e}\otimes\beta_{A})\circ\beta_{e\otimes
A}^{-1})\circ(1_{A}\otimes\beta_{e\otimes A}^{-1})\circ(\gamma_{A}%
\otimes(1_{e}\otimes1_{A}))\\
& ((1_{A}\otimes\epsilon_{A})\otimes(\epsilon_{A}\otimes1_{A}))\circ
\delta_{A\otimes A}\circ r\\
& =(\gamma_{A}\otimes\beta_{A})\circ((1_{A}\otimes\epsilon_{A})\otimes
(\epsilon_{A}\otimes1_{A}))\circ\delta_{A\otimes A}\circ r\\
& =r
\end{align*}

so $\Psi\circ\Phi=1_{\mathcal{R}^{A}(C)}$. Next let $\delta$ be any object in
$\mathcal{S}^{A}(C)$ with underlying object $B$. We have
\begin{align*}
& (\Phi(\Psi(\delta)))^{l}\\
& =(\gamma_{A}\otimes1_{B})\circ((1_{A}\otimes\epsilon_{A})\otimes1_{B}%
)\circ(\Psi(\delta)\otimes1_{B})\circ\delta_{B}\\
& =(\gamma_{A}\otimes1_{B})\circ((1_{A}\otimes\epsilon_{A})\otimes1_{B}%
)\circ((1_{A}\otimes\beta_{A})\otimes1_{B})\\
& \circ((1_{A}\otimes(\epsilon_{B}\otimes1_{A}))\otimes1_{B})\circ
((1_{A}\otimes\delta^{r})\otimes1_{B})\circ(\delta^{l}\otimes1_{B})\circ
\delta_{B}\\
& =(\gamma_{A}\otimes1_{B})\circ((1_{A}\otimes\epsilon_{A})\otimes1_{B}%
)\circ((1_{A}\otimes\beta_{A})\otimes1_{B})\\
& \circ((1_{A}\otimes(\epsilon_{B}\otimes1_{A}))\otimes1_{B})\circ
((1_{A}\otimes\delta^{r})\otimes1_{B})\\
& \circ\alpha_{A,B,B}\circ(1_{A}\otimes\delta_{B})\circ\delta^{l}\\
& =(\gamma_{A}\otimes1_{B})\circ((1_{A}\otimes\epsilon_{e\otimes A}%
)\otimes1_{B})\circ((1_{A}\otimes(\epsilon_{B}\otimes1_{A}))\otimes1_{B})\\
& \circ((1_{A}\otimes\delta\rangle)\otimes1_{B})\circ\alpha_{A,B,B}\circ
(1_{A}\otimes\delta_{B})\circ\delta^{l}\\
& =(\gamma_{A}\otimes1_{B})\circ((1_{A}\otimes\epsilon_{B\otimes A}%
)\otimes1_{B})\circ((1_{A}\otimes\delta^{r})\otimes1_{B})\\
& \circ\alpha_{A,B,B}\circ(1_{A}\otimes\delta_{B})\circ\delta^{l}\\
& =(\gamma_{A}\otimes1_{B})\circ((1_{A}\otimes\epsilon_{B\otimes e}%
)\otimes1_{B})\circ((1_{A}\otimes(1_{B}\otimes\epsilon_{A}))\otimes1_{B})\\
& \circ((1_{A}\otimes\delta^{r})\otimes1_{B})\circ\alpha_{A,B,B}\circ
(1_{A}\otimes\delta_{B})\circ\delta^{l}\\
& =(\gamma_{A}\otimes1_{B})\circ((1_{A}\otimes\epsilon_{B\otimes e}%
)\otimes1_{B})\circ((1_{A}\otimes\gamma_{B}^{-1})\otimes1_{B})\\
& \circ\alpha_{A,B,B}\circ(1_{A}\otimes\delta_{B})\circ\delta^{l}\\
& =(\gamma_{A}\otimes1_{B})\circ((1_{A}\otimes\epsilon_{B})\otimes1_{B}%
)\circ\alpha_{A,B,B}\circ(1_{A}\otimes\delta_{B})\circ\delta^{l}\\
& =(\gamma_{A}\otimes1_{B})\circ\alpha_{A,e,B}\circ(1_{A}\otimes(\epsilon
_{B}\otimes1_{B}))\circ(1_{A}\otimes\delta_{B})\circ\delta^{l}\\
& =(1_{A}\otimes\beta_{B})\circ(1_{A}\otimes(\epsilon_{B}\otimes1_{B}%
))\circ(1_{A}\otimes\delta_{B})\circ\delta^{l}\\
& =\delta^{l}%
\end{align*}

and similarly we find that $(\Phi(\Psi(\delta)))^{r}=\delta^{r}$. This proves
that $\Phi\circ\Psi=1_{\mathcal{S}^{A}(C)}$.
\end{proof}

By definition $\delta_{A}:A\longrightarrow A\otimes A$ is a object in
$\mathcal{R}^{A}(C)$. Its image by $\Phi$ is therefore an object in
$\mathcal{S}^{A}(C)$.

\begin{proposition}
$\Phi(\delta_{A})=\{\delta_{A},\delta_{A}\}$
\end{proposition}

\begin{proof}
We have
\begin{align*}
& \Phi(\delta_{A})^{l}\\
& =(\gamma_{A}\otimes1_{A})\circ((1_{A}\otimes\epsilon_{A})\otimes1_{A}%
)\circ(\delta_{A}\otimes1_{A})\circ\delta_{A}\\
& =(\gamma_{A}\circ(1_{A}\otimes\epsilon_{A})\circ\delta_{A}\otimes1_{A}%
)\circ\delta_{A}\\
& =(1_{A}\otimes1_{A})\circ\delta_{A}\\
& =\delta_{A}%
\end{align*}

In a similar way we prove that $\Phi(\delta_{A})^{r}=\delta_{A}$.
\end{proof}

The object $\{\delta_{A},\delta_{A}\}$ in $\mathcal{S}^{A}(C)$ will play an
important role for a product we will define later and is given a special
name.
\[
a=\{\delta_{A},\delta_{A}\}
\]

In the following we will mostly work in the category $\mathcal{S}^{A}(C)$ and
use the isomorphism to induce the corresponding structures on the category of
relations $\mathcal{R}^{A}(C)$. Note that the category of corelations on
$A$,$\mathcal{R}_{A}(C)$, is by duality isomorphic to the category of $A-A$
bimodules. Denote this category by $\mathcal{S}_{A}(C)$.

\subsection{The $\boxtimes^{A}$ product of relations}

Let $\delta$ and $\gamma$ be two objects in $\mathcal{S}^{A}(C)$ with
underlying objects $B$ and $E$. Define two arrows in $C$ $(\delta\boxtimes
^{A}\gamma)^{l}:B\otimes E\longrightarrow A\otimes(B\otimes E)$ and
$(\delta\boxtimes^{A}\gamma)^{r}:B\otimes E\longrightarrow(B\otimes E)\otimes
A$ by
\begin{align*}
(\delta\boxtimes^{A}\gamma)^{l}  & =\alpha_{A,B,E}^{-1}\circ(\delta^{l}%
\otimes1_{E})\\
(\delta\boxtimes^{A}\gamma)^{r}  & =\alpha_{B,E,A}\circ(1_{B}\otimes\gamma
^{r})
\end{align*}

Then we have

\begin{proposition}
$\delta\boxtimes^{A}\gamma=\{(\delta\boxtimes^{A}\gamma)^{l},(\delta
\boxtimes^{A}\gamma)^{r}\}$ is a object in $\mathcal{S}^{A}(C)$.
\end{proposition}

\begin{proof}
Using the naturality and the MacLane coherence condition for $\alpha$ we have
\begin{align*}
& \alpha_{A,A,B\otimes E}\circ(1_{A}\otimes(\delta\boxtimes^{A}\gamma
)^{l})\circ(\delta\boxtimes^{A}\gamma)^{l}\\
& =\alpha_{A,A,B\otimes E}\circ(1_{A}\otimes\alpha_{A,B,E}^{-1})\circ
(1_{A}(\delta^{l}\otimes1_{E}))\circ\alpha_{A,B,E}^{-1}\circ(\delta^{l}%
\otimes1_{E})\\
& =\alpha_{A,A,B\otimes E}\circ(1_{A}\otimes\alpha_{A,B,E}^{-1})\circ
\alpha_{A,A\otimes B,E}^{-1}\circ((1_{A}\otimes\delta^{l})\circ\delta
^{l}\otimes1_{E})\\
& =\alpha_{A,A,B\otimes E}\circ(1_{A}\otimes\alpha_{A,B,E}^{-1})\circ
\alpha_{A,A\otimes B,E}^{-1}\circ(\alpha_{A,A,B}^{-1}\otimes1_{E})\\
& \circ((\delta_{A}\otimes1_{B})\otimes1_{E})\circ(\delta^{l}\otimes1_{E})\\
& =\alpha_{A\otimes A,B,E}^{-1}\circ((\delta_{A}\otimes1_{B})\otimes
1_{E})\circ(\delta^{l}\otimes1_{E})\\
& =(\delta_{A}\otimes1_{B\otimes E})\circ(\alpha_{A,B,E}^{-1}\circ(\delta
^{l}\otimes1_{E})\\
& =(\delta_{A}\otimes1_{B\otimes E})\circ(\delta\boxtimes^{A}\gamma)^{l}%
\end{align*}

and
\begin{align*}
& \beta_{B\otimes E}\circ(\epsilon_{A}\otimes1_{B\otimes E})\circ
(\delta\boxtimes^{A}\gamma)^{l}\\
& =\beta_{B\otimes E}\circ(\epsilon_{A}\otimes(1_{B}\otimes1_{E}))\circ
\alpha_{A,B,E}^{-1}\circ(\delta^{l}\otimes1_{E})\\
& =\beta_{B\otimes E}\circ\alpha_{e,B,E}^{-1}\circ((\epsilon_{A}\otimes
1_{B})\circ\delta^{l}\otimes1_{E})\\
& =\beta_{B\otimes E}\circ\alpha_{e,B,E}^{-1}\circ(\beta_{B}^{-1}\otimes
1_{E})\\
& =\beta_{B\otimes E}\circ\beta_{B\otimes E}^{-1}\\
& =1_{B\otimes E}%
\end{align*}

so $B\otimes E$ is a left $A$ comodule. In a similar way we show that
$B\otimes E$ is a right $A$ comodule. For the compatibility between the two
structures we have
\begin{align*}
& \alpha_{A,B\otimes E,A}\circ(1_{A}\otimes(\delta\boxtimes^{A}\gamma
)^{r})\circ(\delta\boxtimes^{A}\gamma)^{l}\\
& =\alpha_{A,B\otimes E,A}\circ(1_{A}\otimes\alpha_{B,E,A})\circ(1_{A}%
\otimes(1_{B}\otimes\gamma^{r}))\\
& \circ\alpha_{A,B,E}^{-1}\circ(\delta^{l}\otimes1_{E})\\
& =\alpha_{A,B\otimes E,A}\circ(1_{A}\otimes\alpha_{B,E,A})\circ
\alpha_{A,B,E\otimes A}^{-1}\\
& \circ(\delta^{l}\otimes(1_{E}\otimes1_{A}))\circ(1_{B}\otimes\gamma^{r})\\
& =(\alpha_{A,B,E}^{-1}\otimes1_{A})\circ\alpha_{A\otimes B,E,A}\\
& \circ(\delta^{l}\otimes(1_{E}\otimes1_{A}))\circ(1_{B}\otimes\gamma^{r})\\
& =(\alpha_{A,B,E}^{-1}\otimes1_{A})\circ((\delta^{l}\otimes1_{E})\otimes
1_{A})\\
& \circ\alpha_{B,E,A}\circ(1_{B}\otimes\gamma^{r})\\
& =((\delta\boxtimes^{A}\gamma)^{l}\otimes1_{A})\circ(\delta\boxtimes
^{A}\gamma)^{r}%
\end{align*}
\end{proof}

Using the previous proposition we can define a object map $\boxtimes
^{A}:\mathcal{S}^{A}(C)\times\mathcal{S}^{A}(C)\longrightarrow\mathcal{S}%
^{A}(C) $ by
\[
\boxtimes^{A}(\delta,\gamma)=\delta\boxtimes^{A}\gamma
\]

Let $\delta,\gamma,\rho$ and $\theta$ be objects in $\mathcal{S}^{A}(C)$ with
underlying objects $B,E,D$ and $T$ in $C$ and let $f:\delta\longrightarrow
\rho$ and $g:\gamma\longrightarrow\theta$ be two morphisms in $\mathcal{S}%
^{A}(C)$. Let $f:B\longrightarrow E$ and $g:D\longrightarrow T$ be the
corresponding arrows in $C$ and let $f\otimes g:B\otimes E\longrightarrow
D\otimes T$ \ be their product in $C$. Define $f\boxtimes^{A}g=f\otimes g$.

\begin{proposition}
$f\boxtimes^{A}g:^{A}\delta\boxtimes^{A}\gamma\longrightarrow\rho\boxtimes
^{A}\theta$ is a morphism in $\mathcal{S}^{A}(C)$
\end{proposition}

\begin{proof}
We have previously proved that $f\otimes g$ is a arrow in $C$. It is also a
morphism in $\mathcal{S}^{A}(C)$
\begin{align*}
& (\rho\boxtimes^{A}\theta)^{l}\circ(f\boxtimes^{A}g)\\
& =\alpha_{A,D,T}^{-1}\circ(\rho^{l}\otimes1_{T})\circ(f\otimes g)\\
& =\alpha_{B,D,T}^{-1}\circ(\rho^{l}\circ f\otimes g)\\
& =\alpha_{B,D,T}^{-1}\circ((1_{A}\otimes f)\otimes g)\circ(\delta^{l}%
\otimes1_{E})\\
& =(1_{A}\otimes(f\otimes g))\circ\alpha_{A,B,E}^{-1}\circ(\delta^{l}%
\otimes1_{E})\\
& =(1_{A}\otimes(f\boxtimes^{A}g))\circ(\delta\boxtimes^{A}\gamma)^{l}%
\end{align*}

The identity $(\rho\boxtimes^{A}\theta)^{r}\circ(f\otimes g)=((f\boxtimes
^{A}g)\otimes1_{A})\circ(\delta\boxtimes^{A}\gamma)^{r}$ is proved in a
similar way.
\end{proof}

Using this proposition we can extend the object map $\boxtimes^{A}$ to a
bifunctor \ by defining
\[
\boxtimes^{A}(f,g)=f\boxtimes^{A}g
\]

In terms of this bifunctor we have the following immediate consequence of
lemma \ref{delta_b morphism}

\begin{corollary}
Let $\delta$ be a object in $\mathcal{S}^{A}(C)\mathcal{\ }$with underlying
object $B$. Then $\delta_{B}:\delta\longrightarrow\delta\boxtimes^{A}\delta$
is a morphism in $\mathcal{S}^{A}(C)$.
\end{corollary}

In general there exists no neutral object for $\boxtimes^{A}$. This is clearly
seen in the case of $Sets$. Let $B$ be a set and let $f:B\longrightarrow A$ be
a injective map of sets. Define a $A-A$ bicomodule structure on $B$ by
$\delta^{l}(x)=(f(x),x)$ and $\delta^{r}(x)=(x,f(x))$. Assume that $\omega$ is
a neutral object for $\boxtimes^{A}$and let the underlying object for $\omega$
be $S$. Then there must exist a isomorphism $h:\delta\boxtimes^{A}%
\omega\longrightarrow\delta$ and therefore bijective map $h:B\times
S\longrightarrow B$ that is a morphism of $A-A$ bicomodules. But this implies
that $f(h(x,s))=f(x)$ for all $x$ and $s$. But since $f$ is injective we must
have $h(x,s)=x$ and this is not possible if there is more than one element in
$S$. A neutral element $\omega$ for $\boxtimes^{A}$ therefore must have
$e=\{\ast\}$ as underlying object. Any $A-A$ bicomodule structure $\omega$ on
$e$ must be of the form $\omega^{l}(\ast)=(x_{0},\ast)$ and $\omega^{r}%
(\ast)=(\ast,y_{0})$ for some elements $x_{0},y_{0}\in A$. Let $B$ be a set
with more than one point and let $f:B\longrightarrow A$ be a map of sets that
is not constant. Define a $A-A$ bicomodule structure on $B$ by $\delta
^{l}(x)=(f(x),x)$ and $\delta^{r}(x)=(x,f(x))$. If $\omega$ is a neutral
object for $\boxtimes^{A}$ there must exist a isomorphism $h:e\times
B\longrightarrow B$ that is a morphism of $A-A$ bicomodules. But this implies
that for all $b\in B$ we have $f(h(\ast,b))=x_{0}$ and this implies that $f$
is constant since $h$ is bijective. This is a contradiction and this proves
that $\boxtimes^{A}$ does not have a neutral object in the case of $Sets$.

\subsection{Semimonoidal structures on the category of relations}

Recall that $\langle\mathcal{S}^{A}(C),\boxtimes^{A},M^{A}\rangle$
is a semimonoidal category if
$M^{A}:\boxtimes^{A}\circ(1\times\boxtimes
^{A})\longrightarrow\boxtimes^{A}\circ(\boxtimes^{A}\times1)$ is a
natural isomorphism such that the first of the MacLane Coherence
conditions is satisfied. We will in general assume that the
category $\mathcal{S}^{A}(C) $ is a semimonoidal category with
respect to some choice of natural isomorphism $M^{A}$.

At least one semimonoidal structure for $\boxtimes^{A}$ will always exists.

\begin{proposition}
Let $\delta$ and $\gamma$ be objects in $\mathcal{S}^{A}(C)$ with underlying
objects $B$ and $E$ in $C$. Define
\[
M_{\delta,\gamma,\rho}^{A}=\alpha_{B,E,D}%
\]
where $\alpha$ is the associativity constraint for the category $C$. Then
$\langle\mathcal{S}^{A}(C),\boxtimes^{A},M^{A},S^{A}\rangle$ is a symmetric
semimonoidal category.
\end{proposition}

\begin{proof}
We have proved previously that $\alpha_{B,E,D}$ is a $C$arrow in $C$. Next we
need to show that the associativity constraint for $\otimes$ on $C$ is in fact
a morphism in $\mathcal{S}^{A}(C)$. If we use the naturality and the MacLane
coherence condition for $\alpha$ we have
\begin{align*}
& ((\delta\boxtimes^{A}\gamma)\boxtimes^{A}\rho)\circ\alpha_{B,E,D}\\
& =\alpha_{A,B\otimes E,D}^{-1}\circ(\alpha_{A,B,E}^{-1}\otimes1_{D}%
)\circ((\delta^{l}\otimes1_{B})\otimes1_{D})\circ\alpha_{B,E,D}\\
& =\alpha_{A,B\otimes E,D}^{-1}\circ(\alpha_{A,B,E}^{-1}\otimes1_{D}%
)\circ\alpha_{A\otimes B,E,D}\circ(\delta^{l}\otimes(1_{E}\otimes1_{D}))\\
& =(1_{A}\otimes\alpha_{B,E,D})\circ\alpha_{A,B,E\otimes D}^{-1}\circ
(\delta_{B}\otimes1_{E\otimes D})\\
& =(1_{A}\otimes\alpha_{B,E,D})\circ(\delta\boxtimes^{A}(\gamma\boxtimes
^{A}\rho))
\end{align*}

$M^{A}$ is clearly an isomorphism and is a natural transformation if the
following identity holds \ $((f\boxtimes^{A}g)\boxtimes^{A}h)\circ
M_{\delta,\gamma,\rho}^{A}=M_{\delta^{\prime},\gamma^{\prime},\rho^{\prime}%
}^{A}\circ(f\boxtimes^{A}(g\boxtimes^{A}h))$ for all morphisms $f:\delta
\longrightarrow\delta^{\prime},g:\gamma\longrightarrow\gamma^{\prime}$ and
$h:\rho\longrightarrow\rho^{\prime}$in $\mathcal{S}^{A}(C)$. But the
corresponding identity in $C$ is $((f\otimes g)\otimes h)\circ\alpha
_{B,E,D}=\alpha_{B^{\prime},E^{\prime},D^{\prime}}\circ(f\otimes(g\otimes h))$
and this identity holds because $\alpha$ is a natural transformation.
\end{proof}

The previous proposition leads us to the following definition

\begin{definition}
The semimonoidal $\langle\mathcal{S}^{A}(C),\boxtimes^{A},M^{A}\rangle$
category is external if for all objects $\delta,\gamma$ and $\rho$ we have
\[
M_{\delta,\gamma,\rho}^{A}=\alpha_{B,E,D}%
\]

where $\alpha$ is the associativity constraint for the product $\otimes$ on
the category $C$ and where $B$,$E$ and $D$ are the underlying objects for
$\delta,\gamma$ and $\rho$.
\end{definition}

Since $\mathcal{S}^{A}(C)$ is isomorphic to the category of relations, a
semimonoidal structure on $\mathcal{S}^{A}(C)$ will induce one on the category
of relations. Let the product in $\mathcal{R}^{A}(C)$ corresponding to
$\boxtimes^{A}$ be $\boxdot^{A}:\mathcal{R}^{A}(C)\times\mathcal{R}%
^{A}(C)\longrightarrow\mathcal{R}^{A}(C)$. We thus have
\[
\boxdot^{A}=\Psi\circ\boxtimes^{A}\circ(\Phi\times\Phi)
\]

We have the following explicit expression for the product

\begin{proposition}
Let $r$ and $s$ be two objects in $\mathcal{R}^{A}(C)$. Then we have
\[
r\boxdot^{A}s=(\gamma_{A}\otimes\beta_{A})\circ((1_{A}\otimes\epsilon
_{A})\otimes(\epsilon_{A}\otimes1_{A}))\circ(r\otimes s)
\]
\end{proposition}

\begin{proof}
Since $\varphi$ and $\Psi$ are isomorphisms with $\Phi=\Psi^{-1}$we only need
to verify that
\[
\Phi(r\boxdot^{A}s)=\Phi(r)\boxtimes^{A}\Phi(s)
\]

for all objects $r$ and $s$ in $\mathcal{R}^{A}(C)$. Note that the naturality
of $\epsilon$ gives the following relations
\begin{align*}
& \epsilon_{A}\circ\beta_{A}\circ(\epsilon_{A}\otimes1_{A})\\
& =\epsilon_{e\otimes A}\circ(\epsilon_{A}\otimes1_{A})\\
& =\epsilon_{A\otimes A}\\
& =\epsilon_{A\otimes e}\circ(1_{A}\otimes\epsilon_{A})\\
& =\epsilon_{A}\circ\gamma_{A}\circ(\epsilon_{A}\otimes1_{A})
\end{align*}

W then have
\begin{align*}
& (\Phi(r\boxdot^{A}s))^{l}\\
& =(\gamma_{A}\otimes1_{B\otimes E})\circ((1_{A}\otimes\epsilon_{A}%
)\otimes1_{B\otimes E})\circ((t\boxdot^{A}s)\otimes1_{B\otimes E})\circ
\delta_{B\otimes E}\\
& =(\gamma_{A}\otimes1_{B\otimes E})\circ((1_{A}\otimes\epsilon_{A}%
)\otimes1_{B\otimes E})\circ((\gamma_{A}\otimes\beta_{A})\otimes1_{B\otimes
E})\\
& \circ(((1_{A}\otimes\epsilon_{A})\otimes(\epsilon_{A}\otimes1_{A}%
))\otimes1_{B\otimes E})\circ((r\otimes s)\otimes1_{B\otimes E})\circ
\delta_{B\otimes E}\\
& =(\gamma_{A}\otimes1_{B\otimes E})\circ((1_{A}\otimes\epsilon_{A}%
)\otimes1_{B\otimes E})\circ((\gamma_{A}\otimes\beta_{A})\otimes1_{B\otimes
E})\\
& \circ(((1_{A}\otimes\epsilon_{A})\otimes(\epsilon_{A}\otimes1_{A}%
))\otimes1_{B\otimes E})\circ((r\otimes s)\otimes1_{B\otimes E})\\
& \circ\alpha_{B\otimes E,B,E}^{-1}\circ(\alpha_{B,E,B}\otimes1_{E}%
)\circ((1_{B}\otimes\sigma_{B,E})\otimes1_{E})\\
& \circ(\alpha_{B,B,E}^{-1}\otimes1_{E})\circ\alpha_{B\otimes B,E,E}%
\circ(\delta_{B}\otimes\delta_{E})\\
& =\alpha_{A,B,E}^{-1}\circ((\gamma_{A}\otimes1_{B})\otimes1_{E})\circ
(\alpha_{A,e,B}\otimes1_{E})\circ((1_{A}\otimes\sigma_{B,e})\otimes1_{E})\\
& \circ((1_{A}\otimes(1_{B}\otimes\epsilon_{A}))\otimes1_{E})\circ
(\alpha_{A,B,A}^{-1}\otimes1_{E})\circ\alpha_{A\otimes B,A,E}\\
& \circ((\gamma_{A}\otimes1_{B})\otimes(\beta_{A}\otimes1_{E}))\circ
(((1_{A}\otimes\epsilon_{A})\otimes1_{B})\otimes((\epsilon_{A}\otimes
1_{A})\otimes1_{E}))\\
& \circ((r\otimes1_{B})\otimes(s\otimes1_{E}))\circ(\delta_{B}\otimes
\delta_{E})\\
& =\alpha_{A,B,E}^{-1}\circ((1_{A}\otimes\gamma_{B})\otimes1_{E})\circ
(\alpha_{A,B,e}^{-1}\otimes1_{E})\circ\alpha_{A\otimes B,e,E}\circ(1_{A\otimes
B}\otimes(\epsilon_{A}\otimes1_{E}))\\
& \circ((\gamma_{A}\otimes1_{B})\otimes(\beta_{A}\otimes1_{E}))\circ
(((1_{A}\otimes\epsilon_{A})\otimes1_{B})\otimes((\epsilon_{A}\otimes
1_{A})\otimes1_{E}))\\
& \circ((r\otimes1_{B})\otimes(s\otimes1_{E}))\circ(\delta_{B}\otimes
\delta_{E})\\
& =(1_{A}\otimes(\gamma_{B}\otimes1_{E}))\circ\alpha_{A,B\otimes e,E}%
^{-1}\circ(\alpha_{A,B,e}^{-1}\otimes1_{E})\circ\alpha_{A\otimes B,e,E}%
\circ(1_{A\otimes B}\otimes(\epsilon_{A}\otimes1_{E}))\\
& \circ((\gamma_{A}\otimes1_{B})\otimes(\beta_{A}\otimes1_{E}))\circ
(((1_{A}\otimes\epsilon_{A})\otimes1_{B})\otimes((\epsilon_{A}\otimes
1_{A})\otimes1_{E}))\\
& \circ((r\otimes1_{B})\otimes(s\otimes1_{E}))\circ(\delta_{B}\otimes
\delta_{E})\\
& =(1_{A}\otimes(\gamma_{B}\otimes1_{E}))\circ(1_{A}\otimes\alpha
_{B,e,E})\circ\alpha_{A,B,e\otimes E}^{-1}\circ(1_{A\otimes B}\otimes
(\epsilon_{A}\otimes1_{E}))\\
& \circ((\gamma_{A}\otimes1_{B})\otimes(\beta_{A}\otimes1_{E}))\circ
(((1_{A}\otimes\epsilon_{A})\otimes1_{B})\otimes((\epsilon_{A}\otimes
1_{A})\otimes1_{E}))\\
& \circ((r\otimes1_{B})\otimes(s\otimes1_{E}))\circ(\delta_{B}\otimes
\delta_{E})\\
& =\alpha_{A,B,E}^{-1}\circ(\Phi(r)^{l}\otimes\beta_{E}\circ((\epsilon
_{A}\circ\beta_{A}\circ(\epsilon_{A}\otimes1_{A}))\circ s)\otimes1_{E}%
)\circ\delta_{E})\\
& =\alpha_{A,B,E}^{-1}\circ(\Phi(r)^{l}\otimes\beta_{E}\circ((\epsilon
_{A}\circ\gamma_{A}\circ(1_{A}\otimes\epsilon_{A}))\circ s)\otimes1_{E}%
)\circ\delta_{E})\\
& =\alpha_{A,B,E}^{-1}\circ(\Phi(r)^{l}\otimes\beta_{E}\circ(\epsilon
_{A}\otimes1_{E})\circ\Phi(s)^{l})\\
& =\alpha_{A,B,E}^{-1}\circ(\Phi(r)^{l}\otimes1_{E})\\
& =\Phi(r)^{l}\boxtimes^{A}\Phi(s)^{l}%
\end{align*}

The proof of the identity $(\Phi(r\boxdot^{A}s))^{r}=(\Phi(r)\boxtimes^{A}%
\Phi(s))^{r}$is similar.
\end{proof}

By duality the category of $C$-corelations $\mathcal{R}_{A}(C)$ is isomorphic
to the category of $A-A$ bimodules $\mathcal{S}_{A}(C)$.

\subsection{The tensor product of relations}

For categories of bimodules over rings we have a standard
construction of a tensor product. This construction is categorical
in nature and can in a
natural way be generalized to the category of $A-A$ bimodules $\mathcal{S}%
_{A}(C)$. By dualizing this construction we arrive at our definition of a
tensor product of $A-A$ bicomodules. The isomorphism $\Psi:\mathcal{S}%
^{A}(C)\longrightarrow\mathcal{R}^{A}(C)$ is used to define the tensor product
of relations.

The following lemma is fundamental for the construction of the tensor product.\

\begin{lemma}
Let $\delta$ be a object in $\mathcal{S}^{A}(C)$. Then $\delta^{l}%
:\delta\longrightarrow a\boxtimes^{A}\delta$ and $\delta^{r}:\delta
\longrightarrow\delta\boxtimes^{A}a$ are morphisms in $\mathcal{S}^{A}(C)$ and
if $f:\delta\longrightarrow\gamma$ is a morphism in $\mathcal{S}^{A}(C)$ the
following diagrams in $\mathcal{S}^{A}(C)$ are commutative.%

{\tiny
\begin{equation*}
\begin{aligned}
\begin{diagram}
a \boxtimes\delta&  \rTo^{1_a\boxtimes^A f}            &  a\boxtimes^A\gamma\\
\uTo^{\delta^l}%
&                                              &   \uTo_{\gamma^l}        \\
\delta&   \rTo_{f}                              &     \gamma\\
\end{diagram}
\end{aligned}
\;
\;
\;
\begin{aligned}
\begin{diagram}
\delta\boxtimes^A a                      &  \rTo^{f\boxtimes^A 1_a}%
&  \gamma\boxtimes^A a                 \\
\uTo^{\delta^r}%
&                                              &   \uTo_{\gamma^r}        \\
\delta&   \rTo_{f}                              &     \gamma\\
\end{diagram}
\end{aligned}
\end{equation*}
}%
\end{lemma}

\begin{proof}
There are four diagrams that need to be commutative for the first part of the
lemma to be true. It is seen by inspection that this set of diagrams is
included in the set of diagrams defining $\delta$ to be a $A-A$ bicomodule.
The second part of the lemma is clearly true since the diagrams in $C$
corresponding to the given diagrams are the conditions for $f$ to be a
morphism of the $A-A$ bicomodules $\delta$ and $\gamma$.
\end{proof}

Let now $\delta$ and $\gamma$ be any pair of objects in $\mathcal{S}^{A}(C) $.
>From the previous lemma we can conclude that $\mathcal{P}_{\delta,\gamma}^{A}$
given by%

{\tiny
\begin{diagram}
\delta\boxtimes^A(a\boxtimes^A\gamma
)  &                                                      &   \rTo
^{M^A_{\delta,a,\gamma}}%
&                                                    &    (\delta
\boxtimes^A a)\boxtimes\gamma\\
& \luTo_{1_{\delta}\boxtimes^A\gamma^l}%
&                                         & \ruTo_{\delta^r\boxtimes
^A 1_{\gamma}}  &                                                     \\
&                                                      &       \delta
\boxtimes\gamma
&                                                     &                                                     \\
\end{diagram}
}%

is a diagram in $\mathcal{S}^{A}(C)$. The limit of this diagram, when it
exists, is determined by a object in $\mathcal{S}^{A}(C)$ denoted by
$\delta\otimes^{A}\gamma$ and a morphism $\pi_{\delta,\gamma}^{A}%
:\delta\otimes^{A}\gamma\longrightarrow\delta\boxtimes^{A}\gamma$.

\begin{definition}
Let $\delta$ and $\gamma$ be two objects in $\mathcal{S}^{A}(C)$. The tensor
product of $\delta$ and $\gamma$ is given by
\[
\otimes^{A}(\delta,\gamma)=\delta\otimes^{A}\gamma
\]
\end{definition}

The following property of $\pi_{\delta,\gamma}^{A}$ is important.

\begin{proposition}
\label{mono}$\pi_{\delta,\gamma}^{A}$ is \ a monomorphism.
\end{proposition}

\begin{proof}
Let $\rho$ be a object in $\mathcal{S}^{A}(C)$ and let $f,g:\rho
\longrightarrow\delta\otimes^{A}\gamma$ be a pair of morphisms such that
$\pi_{\delta,\gamma}^{A}\circ f=\pi_{\delta,\gamma}^{A}\circ g$. Define
$h=\pi_{\delta,\gamma}^{A}\circ g$. Then $\langle\rho,h\rangle$ is a cone on
$\mathcal{P}_{\delta,\gamma}^{A}$and therefore the equation
\[
\pi_{\delta,\gamma}^{A}\circ f=h
\]

has a unique solution. But both $f$ and $g$ are solutions and therefore by
uniqueness we can conclude that $f=g$ and this proves that $\pi_{\delta
,\gamma}^{A}$ is a monomorphism.
\end{proof}

We now want to extend the tensorproduct to morphisms. Let $\delta$,$\gamma
$,$\theta$ and $\rho$ be objects in $\mathcal{S}^{A}(C)$ and let
$f:\delta\longrightarrow\theta$ and $g:\gamma\longrightarrow\rho$ be
morphisms. Then we have

\begin{lemma}
$(f\boxtimes^{A}g)\circ\pi_{\delta,\gamma}^{A}$ is a cone on the diagram
$\mathcal{P}_{\theta,\rho}^{A}$.
\end{lemma}

\begin{proof}
We have
\begin{align*}
& M_{\theta,a,\rho}^{A}\circ(1_{\theta}\boxtimes^{A}\rho^{l})\circ
(f\boxtimes^{A}g)\circ\pi_{\delta,\gamma}^{A}\\
& =M_{\theta,a,\rho}^{A}\circ(f\boxtimes^{A}(\rho^{l}\circ g))\circ\pi
_{\delta,\gamma}^{A}\\
& =M_{\theta,a,\rho}^{A}\circ(f\boxtimes^{A}(1_{a}\boxtimes^{A}g)\circ
\gamma^{l})\circ\pi_{\delta,\gamma}^{A}\\
& =M_{\theta,a,\rho}^{A}\circ(f\boxtimes^{A}(1_{a}\boxtimes^{A}g))\circ
(1_{\delta}\boxtimes^{A}\gamma^{l})\circ\pi_{\delta,\gamma}^{A}\\
& =((f\boxtimes^{A}1_{a})\boxtimes^{A}g)\circ M_{\delta,a,\gamma}^{A}%
\circ(1_{\delta}\boxtimes^{A}\gamma^{l})\circ\pi_{\delta,\gamma}^{A}\\
& =((f\boxtimes^{A}1_{a})\boxtimes^{A}g)\circ(\delta^{r}\boxtimes^{A}%
1_{\gamma})\circ\pi_{\delta,\gamma}^{A}\\
& =((f\boxtimes^{A}1_{a})\circ\delta^{r}\boxtimes^{A}g)\circ\pi_{\delta
,\gamma}^{A}\\
& =((\theta^{r}\circ f)\boxtimes^{A}g)\circ\pi_{\delta,\gamma}^{A}\\
& =(\theta^{r}\boxtimes^{A}1_{\rho})\circ(f\boxtimes^{A}g)\circ\pi
_{\delta,\gamma}^{A}%
\end{align*}
\end{proof}

Let $f\otimes^{A}g:\delta\otimes^{A}\gamma\longrightarrow\theta\otimes^{A}%
\rho$ be the unique morphism that exists by the universality of the cone
$\theta\otimes^{A}\rho$. For this morphism we have the commutative diagram%

{\tiny
\begin{diagram}
\delta\boxtimes^A\gamma&   \rTo^{f\boxtimes^A g} &         \theta
\boxtimes^A\rho\\
\uTo^{\pi_{\delta,\gamma}^A} &                                &         \uTo
_{\pi_{\theta,\rho}^A}  \\
\delta\otimes^A\gamma&  \rTo_{f\otimes^A g}        &        \theta
\otimes^A \rho\\
\end{diagram}
}%
\label{pinat}

In general $\delta\otimes^{A}\gamma$ will not exist for all pairs of objects
in $\mathcal{S}^{A}(C)$. In order for it to exists and have reasonable
properties we need to restrict the notion relation as we have defined it. Our
first restriction is to assume that $\otimes^{A}$ is defined for all pairs of
objects in $\mathcal{S}^{A}(C)$. Our second restriction involves the arrow
$\pi_{\delta,\gamma}^{A}$. \ Let $\delta,\gamma$ and $\rho$ be relations. We
require that the morphism $\pi_{\delta,\gamma}^{A}\boxtimes^{A}1_{\rho}$ is
mono. We have proved that $\pi_{\delta,\gamma}^{A}$ is always mono, but
requiring that $\pi_{\delta,\gamma}^{A}\boxtimes^{A}1_{\rho}$ is mono is a
nontrivial restriction in general. It can be thought of a some kind of
''flatness'' condition on $A$.

Given the above restrictions we can define a map $\otimes^{A}:\mathcal{S}%
^{A}(C)\times\mathcal{S}^{A}(C)\longrightarrow\mathcal{S}^{A}(C)$ by
\begin{align*}
\otimes^{A}(\delta,\gamma)  & =\delta\otimes^{A}\gamma\\
\otimes^{A}(f,g)  & =f\otimes^{A}g
\end{align*}

\begin{proposition}
The map $\otimes^{A}$ is a bifunctor
\end{proposition}

\begin{proof}
Let $\delta,\delta^{\prime},\delta^{\prime\prime},\gamma,\gamma^{\prime}$ and
$\gamma^{\prime\prime}$ be four objects in $\mathcal{S}^{A}(C)$ and let
$f:\delta\longrightarrow\delta^{\prime},f^{\prime}:\delta^{\prime
}\longrightarrow\delta^{\prime\prime},g:\gamma\longrightarrow\gamma^{\prime}$
and $g^{\prime}:\gamma^{\prime}\longrightarrow\gamma^{\prime\prime}$ be
morphisms. By universality we know that the equation
\[
\pi_{\delta^{\prime\prime},\gamma^{\prime\prime}}\circ\varphi=((f^{\prime
}\circ f)\boxtimes^{A}(g^{\prime}\circ g))\circ\pi_{\delta,\gamma}^{A}%
\]

has a unique solution. One solution is by definition $(f^{\prime}\circ
f)\otimes^{A}(g^{\prime}\circ g)$. But we also have
\begin{align*}
& \pi_{\delta^{\prime\prime},\gamma^{\prime\prime}}^{A}\circ((f^{\prime
}\otimes^{A}g^{\prime})\circ(f\otimes^{A}g))\\
& =(f^{\prime}\boxtimes^{A}g^{\prime})\circ\pi_{\delta^{\prime},\gamma
^{\prime}}^{A}\circ(f\otimes^{A}g)\\
& =(f^{\prime}\boxtimes^{A}g^{\prime})\circ(f\boxtimes^{A}g)\circ\pi
_{\delta,\gamma}^{A}\\
& =((f^{\prime}\circ f)\boxtimes^{A}(g^{\prime}\circ g))\circ\pi
_{\delta,\gamma}^{A}%
\end{align*}

By uniqueness we must have
\[
(f^{\prime}\circ f)\otimes^{A}(g^{\prime}\circ g)=(f^{\prime}\otimes
^{A}g^{\prime})\circ(f\otimes^{A}g)
\]

Also by universality the following equation has a unique solution.
\[
\pi_{\delta,\gamma}^{A}\circ\varphi=1_{\delta\boxtimes^{A}\gamma}\circ
\pi_{\delta,\gamma}^{A}%
\]

One solution is clearly $1_{\delta\otimes^{A}\gamma}$. But we have
\begin{align*}
& \pi_{\delta,\gamma}^{A}\circ(1_{\delta}\otimes^{A}1_{\gamma})\\
& =(1_{\delta}\boxtimes^{A}1_{\gamma})\circ\pi_{\delta,\gamma}^{A}\\
& =1_{\delta\boxtimes^{A}\gamma}\circ\pi_{\delta,\gamma}^{A}%
\end{align*}

so by uniqueness we have $1_{\delta\otimes^{A}\gamma}=1_{\delta}\otimes
^{A}1_{\gamma}$.
\end{proof}

We will call $\delta\otimes^{A}\gamma$ for the tensor product of the $A-A$
bicomodules $\delta$ and $\gamma$.

We defined the map $\otimes^{A}$ using universal cones in the category
$\mathcal{S}^{A}(C)$ but we will now prove that it can be constructed from
universal cones in the category $C$.

Let $\delta$ and $\gamma$ be two objects in $\mathcal{S}^{A}(C)$ with
underlying objects $B$ and $E$ in $C$. Let the diagram $\mathcal{P}_{B,E}^{A}
$ in $C$ be given by%

{\tiny
\begin{diagram}
B\otimes(A\otimes
E) &                                                      &    \rTo
^{\alpha_{B,A,E}}%
&                                                    &   (B\otimes A)\otimes
E                \\
& \luTo_{1_B\otimes\gamma^l}%
&                                         & \ruTo_{\delta^r\otimes1_E}%
&                                                     \\
&                                                      &       B\otimes
E                &                                                     &                                                     \\
\end{diagram}
}%

Assume that there exists a universal cone $\langle X,h\rangle$ on the diagram
$\mathcal{P}_{B,E}^{A}$ in $C$. Define
\begin{align*}
\Theta^{l}  & =(\gamma_{A}\otimes1_{X})\circ((1_{A}\otimes\epsilon_{B\otimes
E})\otimes1_{X})\circ(\alpha_{A,B,E}^{-1}\otimes1_{X})\\
& \circ((\delta^{l}\otimes1_{E})\otimes1_{X})\circ(h\otimes1_{X})\circ
\delta_{X}\\
\Theta^{r}  & =(1_{X}\otimes\beta_{A})\circ(1_{X}\otimes(\epsilon_{B\otimes
E}\otimes1_{A}))\circ(1_{X}\otimes\alpha_{B,E,A})\\
& \circ(1_{X}\otimes(1_{B}\otimes\gamma^{r}))\circ(1_{X}\otimes h)\circ
\delta_{X}%
\end{align*}

\begin{proposition}
$\Theta=\{\Theta^{l},\Theta^{r}\}$ is a $A-A$ bicomodule with underlying
object $X$.\label{Ldef}
\end{proposition}

\begin{proof}
Let us define morphisms $L,M:X\longrightarrow A$ by
\begin{align*}
L  & =\gamma_{A}\circ(1_{A}\otimes\epsilon_{B\otimes E})\circ\alpha
_{A,B,E}^{-1}\circ(\delta^{l}\otimes1_{E})\circ h\\
M  & =\beta_{A}\circ(\epsilon_{B\otimes E}\otimes1_{A})\circ\alpha
_{B,E,A}\circ(1_{B}\otimes\gamma^{r})\circ h
\end{align*}

\bigskip Then $\Theta^{l}=(L\otimes1_{X})\circ\delta_{X}$ and $\Theta
^{r}=(1_{X}\otimes M)\circ\delta_{X}$ and we have for the left structure
\begin{align*}
& \alpha_{A,A,X}\circ(1_{A}\otimes\Theta^{l})\circ\Theta^{l}\\
& =\alpha_{A,A,X}\circ(1_{A}\otimes(L\otimes1_{X})\circ\delta_{X}%
)\circ(L\otimes1_{X})\circ\delta_{X}\\
& =\alpha_{A,A,X}\circ(1_{A}\otimes(L\otimes1_{X}))\circ(1_{A}\otimes
\delta_{X})\circ(L\otimes1_{X})\circ\delta_{X}\\
& =\alpha_{A,A,X}\circ(1_{A}\otimes(L\otimes1_{X}))\circ(L\otimes(1_{X}%
\otimes1_{X}))\circ(1_{X}\otimes\delta_{X})\circ\delta_{X}\\
& =\alpha_{A,A,X}\circ(1_{A}\otimes(L\otimes1_{X}))\circ(L\otimes(1_{X}%
\otimes1_{X}))\circ\alpha_{X;X;X}^{-1}\circ(\delta_{X}\otimes1_{X})\circ
\delta_{X}\\
& =((1_{A}\otimes L)\otimes1_{X})\circ((L\otimes1_{X})\otimes1_{X}%
)\circ(\delta_{X}\otimes1_{X})\circ\delta_{X}\\
& =((L\otimes L)\circ\delta_{X}\otimes1_{X})\circ\delta_{X}\\
& =(\delta_{A}\circ L\otimes1_{X})\circ\delta_{X}\\
& =(\delta_{A}\otimes1_{X})\circ\Theta^{l}%
\end{align*}

The proof for the right structure is similar. For the compatibility of the
left and right structure we have
\begin{align*}
& (\Theta^{l}\otimes1_{A})\circ\Theta^{r}\\
& =((L\otimes1_{X})\circ\delta_{X}\otimes1_{A})\circ(1_{X}\otimes
M)\circ\delta_{X}\\
& =((L\otimes1_{X})\otimes1_{A})\circ(\delta_{X}\otimes1_{A})\circ
(1_{X}\otimes M)\circ\delta_{X}\\
& =((L\otimes1_{X})\otimes M)\circ(\delta_{X}\otimes1_{X})\circ\delta_{X}\\
& =((L\otimes1_{X})\otimes M)\circ\alpha_{X;X;X}^{-1}\circ(1_{X}\otimes
\delta_{X})\circ\delta_{X}\\
& =\alpha_{A,X,A}^{-1}\circ(L\otimes(1_{X}\otimes M)\circ\delta_{X}%
)\circ\delta_{X}\\
& =\alpha_{A,X,A}^{-1}\circ(1_{A}\otimes\Theta^{r})\circ\Theta^{l}%
\end{align*}
\end{proof}

We will next show that $h$ is a morphism in $\mathcal{S}^{A}(C)$. For this we
need the following lemma

\begin{lemma}%
\[
(\gamma_{A}\circ(1_{A}\otimes\epsilon_{B\otimes E})\otimes(\beta_{B}%
\circ(\epsilon_{A}\otimes1_{B})\otimes1_{E})\circ\alpha_{A,B,E})\circ
\delta_{A\otimes(B\otimes E)}=1_{A\otimes(B\otimes E)}%
\]
\end{lemma}

\begin{proof}
Let $T$ be defined as
\begin{align*}
T  & =((1_{A}\otimes\sigma_{e,e})\otimes1_{B\otimes E})\circ(\alpha
_{A,e,e}^{-1}\otimes1_{B\otimes E})\circ\alpha_{A\otimes e,e,B\otimes E}\\
& \circ((1_{A}\otimes\epsilon_{A})\otimes(\epsilon_{B\otimes E}\otimes
1_{B\otimes E}))\circ(\delta_{A}\otimes\delta_{B\otimes E})
\end{align*}

Then we have
\begin{align*}
& (\gamma_{A}\circ(1_{A}\otimes\epsilon_{B\otimes E})\otimes(\beta_{B}%
\circ(\epsilon_{A}\otimes1_{B})\otimes1_{E})\circ\alpha_{A,B,E})\\
& =(\gamma_{A}\otimes(\beta_{B}\otimes1_{E}))\circ((1_{A}\otimes
\epsilon_{B\otimes E})\otimes((\epsilon_{A}\otimes1_{B})\otimes1_{E}))\\
& \circ(1_{A\otimes(B\otimes E)}\otimes\alpha_{A,B,E})\circ\delta
_{A\otimes(B\otimes E)}\\
& =(\gamma_{A}\otimes(\beta_{B}\otimes1_{E}))\circ(1_{A\otimes e}\otimes
\alpha_{e,B,E})\circ((1_{A}\otimes\epsilon_{B\otimes E})\otimes(\epsilon
_{A}\otimes1_{B\otimes E}))\\
& \circ\alpha_{A\otimes(B\otimes E),A,B\otimes E}^{-1}\circ(\alpha_{A,B\otimes
E,A}\otimes1_{B\otimes E})\circ((1_{A}\otimes\sigma_{A,B\otimes E}%
)\otimes1_{B\otimes E})\\
& \circ(\alpha_{A,A,B\otimes E}^{-1}\otimes1_{B\otimes E})\alpha_{A\otimes
A,B\otimes E,B\otimes E}\circ(\delta_{A}\otimes\delta_{B\otimes E})\\
& =(\gamma_{A}\otimes(\beta_{B}\otimes1_{E}))\circ(1_{A\otimes e}\otimes
\alpha_{e,B,E})\circ\alpha_{A\otimes e,e,B\otimes E}\circ(\alpha
_{A,e,e}\otimes1_{B\otimes E})\circ T\\
& =(\gamma_{A}\otimes(\beta_{B}\otimes1_{E}))\circ\alpha_{A\otimes e,\eta
pB,E}^{-1}\circ(\alpha_{A\otimes e,e,B}\otimes1_{E})\circ\alpha_{(A\otimes
e)\otimes e,B,E}\\
& \circ(\alpha_{A,e,e}\otimes1_{B\otimes E})\circ T\\
& =\alpha_{A,B,E}^{-1}\circ((\gamma_{A}\otimes1_{B})\otimes1_{E}%
)\circ((1_{A\otimes e}\otimes\beta_{B})\otimes1_{E})\circ(\alpha_{A\otimes
e,e,B}\otimes1_{E})\\
& \circ\alpha_{(A\otimes e)\otimes e,B,E}\circ(\alpha_{A,e,e}\otimes
1_{B\otimes E})\circ T\\
& =\alpha_{A,B,E}^{-1}\circ((\gamma_{A}\otimes1_{B})\otimes1_{E})\circ
((\gamma_{A\otimes e}\otimes1_{B})\otimes1_{E})\circ\alpha_{(A\otimes
e)\otimes e,B,E}\\
& \circ(\alpha_{A,e,e}\otimes1_{B\otimes E})\circ T\\
& =(\gamma_{A}\otimes1_{B\otimes E})\circ(\gamma_{A\otimes e}\otimes
1_{B\otimes E})\circ(\alpha_{A,e,e}\otimes1_{B\otimes E})\circ T\\
& =(\gamma_{A}\otimes1_{B\otimes E})\circ((\gamma_{A}\otimes1_{e}%
)\otimes1_{B\otimes E})\circ(\alpha_{A,e,e}\otimes1_{B\otimes E})\circ T\\
& =(\gamma_{A}\otimes1_{B\otimes E})\circ((1_{A}\otimes\beta_{e}%
)\otimes1_{B\otimes E})\circ((1_{A}\otimes\sigma_{e,e})\otimes1_{B\otimes
E})\\
& \circ(\alpha_{A,e,e}^{-1}\otimes1_{B\otimes E})\circ\alpha_{A\otimes
e,e,B\otimes E}\circ((1_{A}\otimes\epsilon_{A})\circ\delta_{A}\otimes
(\epsilon_{B\otimes E}\otimes1_{B\otimes E})\circ\delta_{B\otimes E})\\
& =(\gamma_{A}\otimes1_{B\otimes E})\circ((1_{A}\otimes\gamma_{e}%
)\otimes1_{B\otimes E})\circ(\alpha_{A,e,e}^{-1}\otimes1_{B\otimes E}%
)\circ\alpha_{A\otimes e,e,B\otimes E}\\
& \circ(\gamma_{A}^{-1}\otimes\beta_{B\otimes E}^{-1})\\
& =(\gamma_{A}\otimes1_{B\otimes E})\circ((\gamma_{A}\otimes1_{e}%
)\otimes1_{B\otimes E})\circ\alpha_{A\otimes e,e,B\otimes E}\circ(\gamma
_{A}^{-1}\otimes\beta_{B\otimes E}^{-1})\\
& =(\gamma_{A}\otimes1_{B\otimes E})\circ\alpha_{A,e,B\otimes E}\circ
(\gamma_{A}\otimes(1_{e}\otimes1_{B\otimes E}))\circ(\gamma_{A}^{-1}%
\otimes\beta_{B\otimes E}^{-1})\\
& =(\gamma_{A}\otimes\beta_{B\otimes E})\circ(\gamma_{A}^{-1}\otimes
\beta_{B\otimes E}^{-1})\\
& =1_{A\otimes(B\otimes E)}%
\end{align*}
\end{proof}

We can now prove that $h$ is a morphism in $\mathcal{S}^{A}(C)$.

\begin{proposition}
$h$ defines a monomorphism in $\mathcal{S}^{A}(C)$ with domain $\Theta$ and
codomain $\delta\boxtimes^{A}\gamma$
\end{proposition}

\begin{proof}
The fact that $h$ is a monomorphism in $C$ follows from the universality as i
did for $\pi_{\delta,\gamma}^{A}$ in proposition \ref{mono}.

For the left structure we have
\begin{align*}
& (1_{A}\otimes h)\circ\Theta^{l}\\
& (1_{A}\otimes h)\circ(L\otimes1_{X})\circ\delta_{X}\\
& =(L\otimes1_{B\otimes E})\circ(1_{X}\otimes h)\circ\delta_{X}\\
& =(\gamma_{A}\circ(1_{A}\otimes\epsilon_{B\otimes E})\otimes1_{B\otimes
E})\circ(\alpha_{A,B,E}^{-1}\otimes1_{B\otimes E})\circ((\delta^{l}%
\otimes1_{E})\otimes1_{B\otimes E})\\
& \circ(h\otimes h)\circ\delta_{X}\\
& =(\gamma_{A}\circ(1_{A}\otimes\epsilon_{B\otimes E})\otimes1_{B\otimes
E})\circ(\alpha_{A,B,E}^{-1}\otimes1_{B\otimes E})\circ((\delta^{l}%
\otimes1_{E})\otimes1_{B\otimes E})\\
& \circ\delta_{B\otimes E}\circ h\\
& =(\gamma_{A}\circ(1_{A}\otimes\epsilon_{B\otimes E})\otimes1_{B\otimes
E})\circ(\alpha_{A,B,E}^{-1}\otimes1_{B\otimes E})\\
& \circ((1_{A\otimes B}\otimes1_{E})\otimes(\beta_{B}\otimes1_{E}%
)\circ((\epsilon_{A}\otimes1_{B})\otimes1_{E}))\circ((\delta^{l}\otimes
1_{E})\otimes(\delta^{l}\otimes1_{E}))\\
& \circ\delta_{B\otimes E}\circ h\\
& =(\gamma_{A}\circ(1_{A}\otimes\epsilon_{B\otimes E})\otimes(\beta_{B}%
\circ(\epsilon_{A}\otimes1_{B})\otimes1_{E})\circ\alpha_{A,B,E})\circ
(\alpha_{A,B,E}^{-1}\otimes\alpha_{A,B,E}^{-1})\\
& \circ\delta_{(A\otimes B)\otimes E}\circ(\delta^{l}\otimes1_{E})\circ h\\
& =(\gamma_{A}\circ(1_{A}\otimes\epsilon_{B\otimes E})\otimes(\beta_{B}%
\circ(\epsilon_{A}\otimes1_{B})\otimes1_{E})\circ\alpha_{A,B,E})\circ
\delta_{A\otimes(B\otimes E)}\\
& \circ\alpha_{A,B,E}^{-1}\circ(\delta^{l}\otimes1_{E})\circ h\\
& =\alpha_{A,B,E}^{-1}\circ(\delta^{l}\otimes1_{E})\circ h\\
& =(\delta\boxtimes^{A}\gamma)^{l}\circ h
\end{align*}

In a similar way we show the identity $(h\otimes1_{A})\circ\Theta^{r}%
=(\delta\boxtimes^{A}\gamma)^{r}\circ h$.
\end{proof}

\begin{proposition}
Let the semimonoidal category $\langle\mathcal{S}^{A}(C),\boxtimes^{A}%
,\alpha\rangle$ be external. Then $\langle\Theta,h\rangle$ is a universal cone
on $\mathcal{P}_{\delta,\gamma}^{A}$.
\end{proposition}

\begin{proof}
It is evident that $\langle\Theta,h\rangle$ is a cone on the diagram
$\mathcal{P}_{\delta,\gamma}^{A}$. Let $\langle\theta,u\rangle$ be any cone on
$\mathcal{P}_{\delta,\gamma}^{A}$. Let $\varphi,\psi:\theta\longrightarrow
\Theta$ be two morphisms in $\mathcal{S}^{A}(C)$ such that
\begin{align*}
h\circ\varphi & =\theta\\
h\circ\psi & =\theta
\end{align*}

Then $h\circ\varphi=h\circ\psi$ and since $h$ is mono we have $\varphi=\psi$.
Therefore the equation $h\circ\varphi=\theta$ has at most one solution.

The fact that $\langle\theta,u\rangle$ is a cone gives us the relation
\[
M_{\delta,a,\gamma}^{A}\circ(1_{\delta}\boxtimes^{A}\gamma^{l})\circ
u=(\delta^{l}\boxtimes^{A}1_{\gamma})\circ u
\]

If the underlying objects for $\delta,\gamma$ and $\theta$ are $B,E$ and $D$,
then $u:D\longrightarrow B\otimes E$ and the previous identity corresponds to
the following

identity in $C$
\[
(1_{B}\otimes\gamma^{l})\circ u=(\delta^{l}\otimes1_{E})\circ u
\]

and therefore $\langle D,u\rangle$ is a cone on the diagram $\mathcal{P}%
_{B,E}^{A}$ in $C$. By universality there exists a unique morphism
$\varphi:D\longrightarrow X$ in $C$such that
\[
h\circ\varphi=u
\]

The fact that $\varphi$ is a morphism in $C$ and $u$ and $\delta_{D}$ are
morphisms in $\mathcal{S}^{A}(C)$ gives us the following four commutative diagrams%

{\tiny
\begin{equation*}
\begin{aligned}
\begin{diagram}
D                               &     \rTo^{u}          &   B\otimes
E       \\
\dTo^{\theta^l}      &                                 &   \dTo_{\delta
^l\otimes1_E}  \\
A\otimes D               &    \rTo_{1_A \otimes u} &  A\otimes B\otimes E \\
\end{diagram}
\end{aligned}
\;
\;
\;
\begin{aligned}
\begin{diagram}
D                               &     \rTo^{\delta_D}          &   D\otimes
D       \\
\dTo^{\theta^l}      &                                 &   \dTo_{\delta
^l\otimes1_D}  \\
A\otimes D               &    \rTo_{1_A \otimes\delta_D} &  A\otimes D\otimes
D \\
\end{diagram}
\end{aligned}
\end{equation*}
}

{\tiny
\begin{equation*}
\begin{aligned}
\begin{diagram}
D                               &     \rTo^{\phi}          &   X       \\
\dTo^{\delta_D}      &                                 &   \dTo_{\delta_X}  \\
D\otimes D               &    \rTo_{\phi\otimes\phi} &  X \otimes X \\
\end{diagram}
\end{aligned}
\;
\;
\;
\begin{aligned}
\begin{diagram}
D &                                                      &    \rTo^{u}%
&                                                    &   B\otimes
E                \\
& \rdTo_{\epsilon_D}    &                                         & \ldTo
_{\epsilon_{B\otimes E}}%
&                                                     \\
&                                                      &       e                &
\end{diagram}
\end{aligned}
\end{equation*}
}

If we define $L$ as in proposition \ref{Ldef} we have the following
identities
\begin{align*}
& L\circ\varphi\\
& =\gamma_{A}\circ(1_{A}\otimes\epsilon_{B\otimes E})\circ\alpha_{A,B,E}%
^{-1}\circ(\delta^{l}\otimes1_{E})\circ h\circ\varphi\\
& =\gamma_{A}\circ(1_{A}\otimes\epsilon_{B\otimes E})\circ\alpha_{A,B,E}%
^{-1}\circ(\delta^{l}\otimes1_{E})\circ u\\
& =\gamma_{A}\circ(1_{A}\otimes\epsilon_{B\otimes E})\circ(1_{A}\otimes
u)\circ\theta^{l}\\
& =\gamma_{A}\circ(1_{A}\otimes\epsilon_{D})\circ\theta^{l}%
\end{align*}

But then we have
\begin{align*}
& \Theta^{l}\circ\varphi\\
& =(L\otimes1_{X})\circ\delta_{X}\circ\varphi\\
& =(L\otimes1_{X})\circ(\varphi\otimes\varphi)\circ\delta_{D}\\
& =(L\circ\varphi\otimes\varphi)\circ\delta_{D}\\
& =(1_{A}\otimes\varphi)\circ(\gamma_{A}\otimes1_{D})\circ((1_{A}%
\otimes\epsilon_{D})\otimes1_{D})\circ(\theta^{l}\otimes1_{D})\circ\delta
_{D}\\
& =(1_{A}\otimes\varphi)\circ(\gamma_{A}\otimes1_{D})\circ\alpha_{A,e,D}%
\circ(1_{A}\otimes(\epsilon_{D}\otimes1_{D}))\circ(1_{A}\otimes\delta
_{D})\circ\theta^{l}\\
& =(1_{A}\otimes\varphi)\circ(1_{A}\otimes\beta_{D}\circ(\epsilon_{D}%
\otimes1_{D})\circ\delta_{D})\circ\theta^{l}\\
& =(1_{A}\otimes\varphi)\circ\theta^{l}%
\end{align*}

In a similar way we prove the identity $\Theta^{r}\circ\varphi=(\varphi
\otimes1_{A})\circ\theta^{r}$. This proves that $\varphi$ is a morphism in
$\mathcal{S}^{A}(C)$ and therefore that the equation
\[
h\circ\varphi=u
\]

has a unique solution in $\mathcal{S}^{A}(C)$.
\end{proof}

This proposition show that $\delta\otimes^{A}\gamma\approx\Theta$ since
universal cones are determined up to isomorphism. This is the way the tensor
product is usually computed.

We now have two bifunctors $\boxtimes^{A}$and $\otimes^{A}$ defined on
$\mathcal{S}^{A}(C)$. These two structures are related at the functorial level
as the next proposition show.

\begin{proposition}
$\pi_{\delta,\gamma}^{A}$ are the components of a natural monomorphism
\[
\pi^{A}:\otimes^{A}\longrightarrow\boxtimes^{A}%
\]
\end{proposition}

\begin{proof}
The proposition follows directly from the commutative diagram \ref{pinat}.
\end{proof}

\subsection{Monoidal structures on the category of relations}

We have seen that $\boxtimes^{A}$ defines a semimonoidal structure on
$\mathcal{S}^{A}(C)$ with associativity constraint $M^{A}$. We will in the
following only consider the case when the product $\otimes^{A}$ defines a
monoidal structure on $\mathcal{S}^{A}(C)$ with neutral object $a$. This is a
further restriction on the category $\mathcal{S}^{A}(C)$ and thus on the
category of relations. Recall that the pair $\otimes^{A},a$ defines a monoidal
structure on $\mathcal{S}^{A}(C)$ if for all objects $\delta,\gamma$ and
$\rho$ there exists isomorphisms
\begin{align*}
m_{\delta,\gamma,\rho}^{A}  & :\delta\otimes^{A}(\gamma\otimes^{A}%
\rho)\longrightarrow(\delta\otimes^{A}\gamma)\otimes^{A}\rho\\
l_{\delta}^{A}  & :a\otimes^{A}\delta\longrightarrow\delta\\
r_{\delta}^{A}  & :\delta\otimes^{A}a\longrightarrow\delta
\end{align*}

that are natural in $\delta,\gamma$ and $\rho$ and such that the MacLane
coherence conditions are satisfied. The coherence conditions are a set of
equations for the morphisms $m^{A},l^{A}$ and $r^{A}$ and these equations may
have no solutions, a unique solution or many solutions depending on the
category$\mathcal{\ }C$ and the coalgebra $A$.

\begin{definition}
A monoidal structure $\langle\mathcal{S}^{A}(C),\otimes^{A},a,m^{A}%
,l^{A},r^{A}\rangle$ on the category $\mathcal{S}^{A}(C)$ is induced if
\label{external}\label{adm}for all objects $\delta,\gamma$ and $\rho$ in
$\mathcal{S}^{A}(C)$ the following diagrams commute%

{\tiny
\begin{equation*}
\begin{aligned}
\begin{diagram}
a\otimes^A\delta
&                                                      &   \rTo^{l^A_{\delta}}%
&                                                    &    \delta\\
& \rdTo_{\pi^A_{a,\delta}}%
&                                         & \ldTo_{\delta^l}%
&                                                     \\
&                                                      &       a\boxtimes
^A\delta
&                                                     &                                                     \\
\end{diagram}
\end{aligned}
\;
\;
\;
\begin{aligned}
\begin{diagram}
\delta\otimes
^A a  &                                                      &   \rTo
^{r^A_{\delta}}%
&                                                    &    \delta\\
& \rdTo_{\pi^A_{\delta,a}}%
&                                         & \ldTo_{\delta^r}%
&                                                     \\
&                                                      &       \delta
\boxtimes
^A a              &                                                     &                                                     \\
\end{diagram}
\end{aligned}
\end{equation*}
}

{\tiny
\begin{diagram}
\delta\boxtimes^A(\gamma\boxtimes^A\rho)               &\rTo^{M^A_{\delta
,\gamma,\rho}}    &             (\delta\boxtimes^A\gamma)\boxtimes\rho\\
\uTo^{1_{\delta}\boxtimes^A\pi_{\gamma,\rho}^A}%
&                                        &              \uTo_{\pi
_{\delta,\gamma}^A\otimes1_{\rho}}     \\
\delta\boxtimes^A(\gamma\otimes^A\rho
)          &                                        &              (\delta
\otimes^A\gamma)\boxtimes^A\rho\\
\uTo^{\pi_{\delta,\gamma\otimes^A \rho}^A}%
&                                        &              \uTo_{\pi
_{\delta\otimes^A \gamma,\rho}^A}   \\
\delta\otimes^A(\gamma\otimes^A\rho)        & \rTo_{m^A_{\delta,\gamma,\rho}%
}&            (\delta\otimes^A\gamma)\otimes^A\rho\\
\end{diagram}
}%
\end{definition}

We will in the following derive a necessary and sufficient
condition for induced constraints to exist in the external case. \
Let us assume that the semimonoidal category
$\langle\mathcal{S}^{A}(C),\boxtimes^{A},M^{A}\rangle$
is external. Let $\mathcal{P}_{a,\delta}^{A}$ be the diagram%

{\tiny
\begin{diagram}
a\boxtimes^A(a\boxtimes^A\delta
)  &                                                      &   \rTo
^{M^A_{a,a,\delta}}%
&                                                    &    (\delta
\boxtimes^A a)\boxtimes\gamma\\
& \luTo_{1_{a}\boxtimes^A\delta^l}%
&                                         & \ruTo_{\delta_A\boxtimes
^A 1_{\delta}}  &                                                     \\
&                                                      &       a \boxtimes
\delta
&                                                     &                                                     \\
\end{diagram}
}%

Then $\langle\delta,\delta^{l}\rangle$ is clearly a cone on this diagram since
this is equivalent to the condition that $\delta^{l}$ is a left comodule
structure on the underlying object of $\delta$. But we have also a stronger condition.

\begin{proposition}
Let the semimonoidal category $\langle\mathcal{S}^{A}(C),\boxtimes^{A}%
,M^{A}\rangle$ be external. Then $\langle\delta,\delta^{l}\rangle$ is a
universal cone on $\mathcal{P}_{a,\delta}^{A}$.
\end{proposition}

\begin{proof}
In order to prove that $\langle\delta,\delta^{l}\rangle$ is a universal cone
we must show that the equation
\[
\delta^{l}\circ\varphi=f
\]

has a unique solution $\varphi:\gamma\longrightarrow\delta$ for any
$f:\gamma\longrightarrow a\boxtimes^{A}\delta$ such that
\[
(\delta_{A}\otimes1_{\delta})\circ f=M_{a,a,\delta}^{A}\circ(1_{a}%
\otimes\delta^{l})\circ f
\]

Since $\beta_{B}\circ(\epsilon_{A}\otimes1_{B})\circ\delta^{l}=1_{B}$ the
equation can have only one solution and this solution must be
\[
\varphi=\beta_{B}\circ(\epsilon_{A}\otimes1_{B})\circ f
\]

The universality is proved if we can show that this is in fact a solution and
also a morphism in $\mathcal{S}^{A}(C).$
\begin{align*}
& \delta^{l}\circ\varphi\\
& =\delta^{l}\circ\beta_{B}\circ(\epsilon_{A}\otimes1_{B})\circ f\\
& =\beta_{A\otimes B}\circ(1_{e}\otimes\delta^{l})\circ(\epsilon_{A}%
\otimes1_{B})\circ f\\
& =\beta_{A\otimes B}\circ(\epsilon_{A}\otimes(1_{A}\otimes1_{B}))\circ
(1_{A}\otimes\delta^{l})\circ f\\
& =\beta_{A\otimes B}\circ(\epsilon_{A}\otimes(1_{A}\otimes1_{B}))\circ
\alpha_{A,A,B}^{-1}\circ(\delta_{A}\otimes1_{B})\circ f\\
& =\beta_{A\otimes B}\circ\alpha_{e,A,B}^{-1}\circ(\beta_{A}^{-1}\otimes
1_{B})\circ f\\
& =f
\end{align*}

so $\varphi$ is a solution. Here we have used the identity $\beta_{A\otimes
B}=(\beta_{A}\otimes1_{B})\circ\alpha_{e,A,B}$. By construction $\varphi$ is a
arrow in $C$, but we also have
\begin{align*}
& (1_{A}\otimes\varphi)\circ\gamma^{l}\\
& =(1_{A}\otimes\beta_{B})\circ(1_{A}\otimes(\epsilon_{A}\otimes1_{B}%
))\circ(1_{A}\otimes f)\circ\gamma^{l}\\
& =(1_{A}\otimes\beta_{B})\circ(1_{A}\otimes(\epsilon_{A}\otimes1_{B}%
))\circ(1_{A}\otimes\delta^{l})\circ f\\
& =(1_{A}\otimes\beta_{B})\circ(1_{A}\otimes(\epsilon_{A}\otimes1_{B}%
))\circ\alpha_{A,A,B}^{-1}\circ(\delta_{A}\otimes1_{B})\circ f\\
& =(1_{A}\otimes\beta_{B})\circ\alpha_{A,e,B}^{-1}\circ((1_{A}\otimes
\epsilon_{A})\otimes1_{B})\circ(\delta_{A}\otimes1_{B})\circ f\\
& =(\gamma_{A}\otimes1_{B})\circ(\gamma_{A}^{-1}\otimes1_{B})\circ f\\
& =f\\
& =\delta^{l}\circ\varphi
\end{align*}

so $\varphi$ is a morphism in $\mathcal{S}^{A}(C)$.
\end{proof}

We have the following two corollaries to the previous proposition

\begin{corollary}
Let the semimonoidal category $\langle\mathcal{S}^{A}(C),\boxtimes^{A}%
,M^{A}\rangle$ be external and let the underlying object for $\delta$ be $B$.
If induced unit constraints exists they must be of the form
\begin{align*}
l_{\delta}^{A}  & =\beta_{B}\circ(\epsilon_{A}\otimes1_{B})\circ\pi_{a,\delta
}^{A}\\
r_{\delta}^{A}  & =\gamma_{B}\circ(1_{B}\otimes\epsilon_{A})\circ\pi
_{\delta,a}^{A}%
\end{align*}
\end{corollary}

\begin{corollary}
Let the semimonoidal category $\langle\mathcal{S}^{A}(C),\boxtimes^{A}%
,M^{A}\rangle$ be external.Then the morphism $\delta^{l}:\delta\longrightarrow
a\otimes^{A}\delta$ is a monomorphism.
\end{corollary}

\begin{proof}
Let $\gamma$ be any object in $\mathcal{S}^{A}(C)$ and let
$f,g:\gamma \longrightarrow\delta$ be any pair of morphisms.
Assume that $\delta^{l}\circ f=\delta^{l}\circ g$ and define
$h=\delta^{l}\circ g$. Then both $f$ and $h$ satisfy the equation
\[
\delta^{l}\circ\varphi=h
\]

By universality we can conclude that $f=g$.
\end{proof}

\begin{proposition}
Let the semimonoidal category $\langle\mathcal{S}^{A}(C),\boxtimes^{A}%
,M^{A}\rangle$ be external. Then induced unit and associativity constraints
are unique if they exist.
\end{proposition}

\begin{proof}
By definition a external unit constraint $l_{\delta}^{A}$ is a solution of the
equation
\[
\delta^{l}\circ l_{\delta}^{A}=\pi_{a,\delta}^{A}%
\]

But by universality this equation has a unique solution. The uniqueness of
$r_{\delta}^{A}$ is proved in a similar way. For $m_{r,s,t}^{A}$ we note that
the morphism $t=(\pi_{\delta,\gamma}^{A}\otimes1_{\rho})\circ\pi
_{\delta\otimes^{A}\gamma,\rho}^{A}$ is mono. Let $f$ and $g$ be two morphisms
such that the third diagram in definition \ref{external} commutes. Then we
have
\begin{align*}
t\circ f  & =(\pi_{\delta,\gamma}^{A}\otimes1_{\rho})\circ\pi_{\delta
\otimes^{A}\gamma,\rho}^{A}\circ f=\alpha_{B,E,D}\circ(1_{\delta}\otimes
\pi_{\gamma,\rho}^{A})\circ\pi_{\delta,\gamma\otimes^{A}\rho}^{A}\\
t\circ g  & =(\pi_{\gamma,\delta}^{A}\otimes1_{\rho})\circ\pi_{\delta
\otimes^{A}\gamma,\rho}^{A}\circ g=\alpha_{B,E,D}\circ(1_{\delta}\otimes
\pi_{\gamma,\rho}^{A})\circ\pi_{\delta,\gamma\otimes^{A}\rho}^{A}%
\end{align*}

so $t\circ f=t\circ g$. But $t$ is mono and therefore $f=g$
\end{proof}

We can now give sufficient conditions for the existence of induced unit and
associativity constraints.

\begin{theorem}
Let the semimonoidal category $\langle\mathcal{S}^{A}(C),\boxtimes^{A}%
,M^{A}\rangle$ be external and assume that for all $\delta,\gamma$ and $\rho$
there exists a isomorphism $m_{\delta,\gamma,\rho}^{A}:\delta\otimes
^{A}(\gamma\otimes^{A}\rho)\longrightarrow(\delta\otimes^{A}\gamma)\otimes
^{A}\rho$ such that the following diagram commute.%

{\tiny
\begin{diagram}
\delta\boxtimes^A(\gamma\boxtimes^A\rho)               &\rTo^{M^A_{\delta
,\gamma,\rho}}    &             (\delta\boxtimes^A\gamma)\boxtimes\rho\\
\uTo^{1_{\delta}\boxtimes^A\pi_{\gamma,\rho}^A}%
&                                        &              \uTo_{\pi
_{\delta,\gamma}^A\otimes1_{\rho}}     \\
\delta\boxtimes^A(\gamma\otimes^A\rho
)          &                                        &              (\delta
\otimes^A\gamma)\boxtimes^A\rho\\
\uTo^{\pi_{\delta,\gamma\otimes^A \rho}^A}%
&                                        &              \uTo_{\pi
_{\delta\otimes^A \gamma,\rho}^A}   \\
\delta\otimes^A(\gamma\otimes^A\rho)        & \rTo_{m^A_{\delta,\gamma,\rho}%
}&            (\delta\otimes^A\gamma)\otimes^A\rho\\
\end{diagram}
}%
\label{external 2}

Then a induced monoidal structure $\langle\mathcal{S}^{A}(C),\otimes
^{A},a,m^{A},l^{A},r^{A}\rangle$ on the category $\mathcal{S}^{A}(C)$ exists.
\end{theorem}

\begin{proof}
Let us first prove that $m_{\delta,\gamma,\rho}^{A}$ are the components of a
natural isomorphism. Let $\delta^{\prime},\gamma^{\prime}$ and $\rho^{\prime}$
be three other relations and let $f:\delta\longrightarrow\delta^{^{\prime}%
},g:\gamma\longrightarrow\gamma^{\prime}$ and $h:\rho\longrightarrow
\rho^{\prime}$ be three morphisms. For any three relations $\delta,\gamma$ and
$\rho$ define $\varphi_{\delta,\gamma,\rho}=(\pi_{\delta,\gamma}^{A}%
\otimes1_{\rho})\circ\pi_{\delta\otimes^{A}\gamma,\rho}^{A}$ and $\psi
_{\delta,\gamma,\rho}=(1_{\delta}\otimes^{A}\pi_{\gamma,\rho}^{A})\circ
\pi_{\delta,\gamma\otimes^{A}\rho}^{A}$. Then we have
\begin{align*}
& \varphi_{\delta^{\prime},\gamma^{\prime},\rho^{\prime}}\circ((f\otimes
^{A}g)\otimes^{A}h)\circ m_{\delta,\gamma,\rho}^{A}\\
& =((f\boxtimes^{A}g)\boxtimes^{A}h)\circ\varphi_{\delta,\gamma,\rho}\circ
m_{\delta,\gamma,\rho}^{A}\\
& =((f\boxtimes^{A}g)\boxtimes^{A}h)\circ M_{\delta,\gamma,\rho}^{A}\circ
\psi_{\delta,\gamma,\rho}\\
& =M_{\delta^{\prime},\gamma^{\prime},\rho^{\prime}}^{A}\circ(f\boxtimes
^{A}(g\boxtimes^{A}h))\circ\psi_{\delta,\gamma,\rho}\\
& =M_{\delta^{\prime},\gamma^{\prime},\rho^{\prime}}^{A}\circ\psi
_{\delta^{\prime},\gamma^{\prime},\rho^{\prime}}\circ(f\otimes^{A}%
(g\otimes^{A}h))\\
& =\varphi_{\delta^{\prime},\gamma^{\prime},\rho^{\prime}}\circ m_{\delta
^{\prime},\gamma^{\prime},\rho^{\prime}}^{A}\circ(f\otimes^{A}(g\otimes^{A}h))
\end{align*}

>From this the naturality of $m_{\delta,\gamma,\rho}^{A}$ follows because
$\varphi_{\delta^{\prime},\gamma^{\prime},\rho^{\prime}}$ is mono. We have
thus far proved that we have a natural isomorphism
\[
m_{\delta,\gamma,\rho}^{A}:\delta\otimes^{A}(\gamma\otimes^{A}\rho
)\longrightarrow(\delta\otimes^{A}\gamma)\otimes^{A}\rho
\]

A induced left unit constraint is a natural isomorphism in $\mathcal{S}%
^{A}(C)$ that satisfy the equation
\[
\delta^{l}\circ\varphi=\pi_{a,\delta}^{A}%
\]

By definition $a\otimes^{A}\delta$ is a universal cone on the diagram
$\mathcal{P}_{a,\delta}^{A}$. But $\delta^{l}$ is also a universal cone on
this diagram so there exists a isomorphism $\varphi:a\otimes^{A}%
\delta\longrightarrow\delta$ such that $\delta^{l}\circ\varphi=\pi_{a,\delta
}^{A}$. We have seen that the only solution of this equation in $\mathcal{S}%
^{A}(C)$ is given by $l_{\delta}^{A}=\beta_{B}\circ(\epsilon_{A}\otimes
1_{B})\circ\pi_{a,\delta}^{A}$ where the underlying object for $\delta$ is
$B$. We can therefore conclude that $l_{\delta}^{A}:a\otimes^{A}%
\delta\longrightarrow\delta$ is a isomorphism. This isomorphism is natural
because if $\delta$ and $\delta^{\prime}$ are objects in $\mathcal{S}^{A}(C)$
with underlying objects $B$ and $B^{\prime}$ and $f:\delta\longrightarrow
\delta^{\prime}$ is any morphism the naturality of $\beta$ and $\pi^{A}$ give
\begin{align*}
& f\circ l_{\delta}^{A}\\
& =f\circ\beta_{B}\circ(\epsilon_{A}\otimes1_{B})\circ\pi_{a,\delta}^{A}\\
& =\beta_{B^{\prime}}\circ(1_{e}\otimes f)\circ(\epsilon_{A}\otimes1_{B}%
)\circ\pi_{a,\delta}^{A}\\
& =\beta_{B^{\prime}}\circ(\epsilon_{A}\otimes f)\circ\pi_{a,\delta}^{A}\\
& =\beta_{B^{\prime}}\circ(\epsilon_{A}\otimes1_{B^{\prime}})\circ
(1_{A}\otimes f)\circ\pi_{a,\delta}^{A}\\
& =\beta_{B^{\prime}}\circ(\epsilon_{A}\otimes1_{B^{\prime}})\circ
\pi_{a,\delta^{\prime}}^{A}\circ f\\
& =l_{\delta^{\prime}}^{A}\circ f
\end{align*}

In a similar way we find a natural isomorphism $r_{\delta}^{A}=\gamma_{B}%
\circ(1_{B}\otimes\epsilon_{A})\circ\pi_{\delta,a}^{A}$. The
proposition is proved if we can show that these three natural
isomorphisms satisfy the MacLane coherence conditions for a
monoidal category. We clearly have
\[
a^{l}\circ l_{a}^{A}=\pi_{a,a}^{A}=a^{r}\circ r_{a}^{A}%
\]

But $a^{l}=\delta_{A}=a^{r}$ and $\delta_{A}$ is a monomorphism so we have
$l_{a}^{A}=r_{a}^{A}$. This is the third MacLane condition. Let us now
consider the second MacLane condition. Let $\delta$ and $\gamma$ have
underlying objects $B$ and $E$. Using the formulas for $l_{\delta}%
^{A},r_{\delta}^{A}$ and the definition of $m_{\delta,\gamma,\rho}^{A}$ we
find
\begin{align*}
& \pi_{\delta,\gamma}^{A}\circ(r_{\delta}^{A}\otimes^{A}1_{\gamma})\circ
m_{\delta,a,\gamma}^{A}\\
& =(r_{\delta}^{A}\otimes1_{E})\circ\pi_{\delta\otimes^{A}a,\gamma}^{A}\circ
m_{\delta,a,\gamma}^{A}\\
& =(\gamma_{B}\circ(1_{B}\otimes\epsilon_{A})\circ\pi_{\delta,a}^{A}%
\otimes1_{E})\circ\pi_{\delta\otimes^{A}a,\gamma}^{A}\circ m_{\delta,a,\gamma
}^{A}\\
& =(\gamma_{B}\otimes1_{E})\circ((1_{B}\otimes\epsilon_{A})\otimes1_{E}%
)\circ(\pi_{\delta,a}^{A}\boxtimes^{A}1_{\gamma})\circ\pi_{\delta\otimes
^{A}a,\gamma}^{A}\circ m_{\delta,a,\gamma}^{A}\\
& =(\gamma_{B}\otimes1_{E})\circ((1_{B}\otimes\epsilon_{A})\otimes1_{E}%
)\circ\alpha_{B,A,E}\circ(1_{\delta}\boxtimes^{A}\pi_{a,\gamma}^{A})\circ
\pi_{\delta,a\otimes^{A}\gamma}^{A}\\
& =(\gamma_{B}\otimes1_{E})\circ\alpha_{B,e,E}\circ(1_{B}\otimes(\epsilon
_{A}\otimes1_{E}))\circ(1_{\delta}\boxtimes^{A}\pi_{a,\gamma}^{A})\circ
\pi_{\delta,a\otimes^{A}\gamma}^{A}\\
& =(1_{B}\otimes\beta_{E})\circ(1_{B}\otimes(\epsilon_{A}\otimes1_{E}%
))\circ(1_{B}\otimes\pi_{a,\gamma}^{A})\circ\pi_{\delta,a\otimes^{A}\gamma
}^{A}\\
& =(1_{B}\otimes\beta_{E}\circ(\epsilon_{A}\otimes1_{E})\circ\pi_{a,\gamma
}^{A})\circ\pi_{\delta,a\otimes^{A}\gamma}^{A}\\
& =(1_{\delta}\boxtimes^{A}l_{\gamma}^{A})\circ\pi_{\delta,a\otimes^{A}\gamma
}^{A}\\
& =\pi_{\delta,\gamma}^{A}\circ(1_{\delta}\otimes^{A}l_{\gamma}^{A})
\end{align*}

and $\pi_{\delta,\gamma}^{A}$ is mono so we have
\[
(r_{\delta}^{A}\otimes^{A}1_{\gamma})\circ m_{\delta,a,\gamma}^{A}=(1_{\delta
}\otimes^{A}l_{\gamma}^{A})
\]

and this is the second MacLane condition. The first MacLane condition follow
from the assumptions in the Theorem and the fact that $\boxtimes^{A}$ is a
semimonoidal structure on $\mathcal{S}^{A}(C)$ with associativity constraint
$M^{A}$.
\end{proof}

Since $\mathcal{S}^{A}(C)$ is isomorphic to the category of relations a
monoidal structure on $\mathcal{S}^{A}(C)$ will induce one on the category of
relations. Let the product in $\mathcal{R}^{A}(C)$ corresponding to
$\otimes^{A}$ be $\odot^{A}:\mathcal{R}^{A}(C)\times\mathcal{R}^{A}%
(C)\longrightarrow\mathcal{R}^{A}(C)$. We thus have
\[
\odot^{A}=\Psi\circ\otimes^{A}\circ(\Phi\times\Phi)
\]

We have the following explicit expression for the product

\begin{proposition}
For any pair of objects $r$ and $s$ in$\mathcal{R}^{A}(C)$ we have
\[
r\odot^{A}s=(\gamma_{A}\otimes\beta_{A})\circ(1_{A}\otimes\epsilon_{A}%
\otimes\epsilon_{A}\otimes1_{A})\circ(r\otimes s)\circ\pi_{\Phi(r),\Phi
(s)}^{A}%
\]
\end{proposition}

\begin{proof}
We have a natural monomorphism
\[
\pi^{A}:\otimes^{A}\longrightarrow\boxtimes^{A}%
\]

Since by definition $\boxdot^{A}=\Psi\circ\boxtimes^{A}\circ(\Phi\times\Phi)$
and $\odot^{A}=\Psi\circ\otimes^{A}\circ(\Phi\times\Phi)$, horizontal
composition of natural transformations give us a natural transformation
\[
(1_{\Psi}\circ\pi^{A}\circ(1_{\Phi}\times1_{\Phi})):\odot^{A}\longrightarrow
\boxdot^{A}%
\]

If we evaluate this natural transformation at a pair of objects $r$ and $s$ in
$\mathcal{R}^{A}(C)$ we get the following morphism in $\mathcal{R}^{A}(C)$%
\[
\pi_{\Phi(r),\Phi(s)}^{A}:r\odot^{A}s\longrightarrow r\boxdot^{A}s
\]

But this means that
\[
r\odot^{A}s=(r\boxdot^{A}s)\circ\pi_{\Phi(r),\Phi(s)}^{A}%
\]

and this is the formula in the proposition if we take into account the formula
for $r\boxdot^{A}s$ that we have derived earlier.
\end{proof}

We will now consider a few examples of the tensor product. Let us
first assume that the underlying category $C$ is $Sets$ with its
unique choice of natural $C$-category $C$. We have seen that all
possible $A-A$ bicomodule structures
$\delta=\{\delta^{l},\delta^{r}\}$on a set $B$ are of the form
$\delta ^{l}(x)=(f(x),x)$ and $\delta^{r}(x)=(x,g(x))$ for some
functions $f,g:B\longrightarrow A$. The relation on $A$
corresponding to $\delta$ is clearly given by $r(x)=(f(x),g(x))$.
Let now $r(x)=(f(x),g(x))$ and $s(y)=(h(y),k(y))$ be two relations
with domains $B$ and $E$ and let the corresponding $A-A$
bicomodules be $\delta$ and $\gamma$. The two maps
$(\delta^{r}\times1_{E})$ and $(1_{B}\times\gamma^{l})$ are given
by
\begin{align*}
(\delta^{r}\times1_{E})(x,y)  & =(x,g(x),y)\\
(1_{B}\times\gamma^{l})(x,y)  & =(x,h(y),y)
\end{align*}

In $Sets$ the underlying object $X$ for the $A-A$ bicomodule
$\delta \otimes^{A}\gamma$ is the equalizer of the two given maps.
We therefore find that
\[
X=\{(x,y)|\;g(x)=h(y)\}
\]

and
\begin{align*}
(\delta\otimes^{A}\gamma)^{l}(x,y)  & =(f(x),x,y)\\
(\delta\otimes^{A}\gamma)^{r}(x,y)  & =(x,y,x,k(y))
\end{align*}

The map $\pi_{\delta,\gamma}^{A}:X\longrightarrow B\times E$ is the inclusion
map. The relation $r\odot^{A}s$ corresponding to $\delta\otimes^{A}\gamma$ is
then given by $(r\odot^{A}s)(x,y)=(f(x),k(y))$.

We have seen that each relation $r$ and $s$ is in fact a directed
labelled graph. Each element in $B$ can be thought of as a arrow
that has a source and a target in the vertex set $A$ for the graph
and similarly for elements in $E$. Let us define two arrows to be
composable if the target of the first is the same as the source of
the second. The set $X$ then consists of all composable pairs of
arrows from $B$ and $E$. Two relations on $A$, $B\subset A\times
A$ and $E\subset A\times A$, in the usual sense corresponds to
relations $r$ and $s$ in our sense if we let
$r(x,x^{\prime})=(x,x^{\prime})$ and
$s(y,y^{\prime})=(y,y^{\prime})$ be the inclusion maps. If we use
the same notation as above we find
$f(x,x^{\prime})=x,g(x,x^{\prime})=x^{\prime },h(y,y^{\prime})=y$
and $k(y,y^{\prime})=y^{\prime}$.

For this special case we find
\begin{align*}
& X=\{((x,y),(y,y^{\prime}))\}\\
& (r\odot^{A}s)((x,y),(y,y^{\prime}))=(x,y^{\prime})
\end{align*}

We then observe that
\[
(r\odot^{A}s)(X)=B\circ E
\]

where $B\circ E$ is the usual composition of relations.

Let $\ t:D\longrightarrow A\times A$ be a third relation with
$t(z)=(p(z),q(z))$ and let $\rho$ be the $A-A$ bicomodule corresponding to $t
$. Let $X$ be the underlying object for $(\delta\otimes^{A}\gamma)\otimes
^{A}\rho$ and $Y$ the underlying object for $\delta\otimes^{A}(\gamma
\otimes^{A}\rho)$. Direct calculation show that
\begin{align*}
X  & =\{((x,y),z)|\;g(x)=h(y),k(y)=p(z)\}\\
Y  & =\{(x,(y,z))|\;g(x)=h(y),k(y)=p(z)\}
\end{align*}

Define
\begin{align*}
m_{\delta,\gamma,\rho}^{A}(x,(y,z))  & =((x,y),z)\\
l_{\delta}^{A}(a,x)  & =x\\
r_{\delta}^{A}(x,a)  & =x
\end{align*}

It is easy to see that $m_{\delta,\gamma,\rho}^{A}$ is a morphism of
relations. The underlying object for $a\otimes^{A}\delta$ is easily seen to
given by
\[
Z=\{(f(x),x)|\;x\in B\}
\]

Therefore $l_{\delta}^{A}$ is clearly a isomorphism and%

\begin{align*}
& (\delta^{l}\circ l_{\delta}^{A})(f(x),x)\\
& =\delta^{l}(x)\\
& =(f(x),x)\\
& =\pi_{a,\delta}^{A}(f(x),x)
\end{align*}

In a similar way we show that $r_{\delta}^{A}$ is a isomorphism
that satisfy $\delta^{r}\circ r_{\delta}^{A}=\pi_{\delta,a}^{A}$.
\ It is easy to see that $\omega_{r,s,t}^{A},\beta_{r}^{A}$ and
$\gamma_{r}^{A}$ satisfy the MacLane coherence conditions. They
are therefore the associativity and unit constraints for a
monoidal structure $\otimes^{A}$on $\mathcal{S}^{A}(C)$. A simple
calculation show that they are induced.

In a similar way we can define products of any number of
relations. It is evident that the product of $n$ relations
consists of strings of composable arrows of length $n$, one arrow
from each relation. Note that $\delta_{A}$ is a relation on $A$.
Let us assume that there exists a morphism of relations
$f:\delta_{A}\longrightarrow r$. \ This means that for each
element $a\in A$ there exists a element $b_{a}=f(a)$ in $B$ such
that $a$ is both the source and target of $b_{a}$. If we now take
all possible products of the relation $r$ we observe that the
result is in fact the (internal) category generated by the graph
defined by the relation.

Let us next consider $Vect_{k}$ with $\oplus$ as monoidal structure. Let $A$
be a linear space and let $r:B\longrightarrow A\oplus A$ and
$s:E\longrightarrow A\oplus A$ be relations on $A$. The domain for the product
$r\odot^{A}s$ is a linear subspace of $B\oplus E$
\begin{align*}
V  & =\{(u,v)|\;g(u)=h(v)\}\\
(r\odot^{A}s)(u,v)  & =(f(u),k(v))
\end{align*}

where $r(u)=(f(u),g(u))$ and $s(v)=(h(v),k(v))$. Let $L:A\longrightarrow A$
and $S:A\longrightarrow A$ be two endomorphism of $A$ and let $B=E=A$ and
$r:B\longrightarrow A\oplus A$ $s:E\longrightarrow A\oplus A$ be relations
where $r(a)=(a,L(a))$ and $s(a)=(a,S(a)\dot{)}$. We then have
$f(a)=a,g(a)=L(a),h(a)=a$ and $k(a)=S(a)$. Therefore \ the underlying object
for $r\odot^{A}s$ is $X=\{(a,L(a))\}$ and
\[
(r\odot^{A}s)(a,L(a))=(a,(S\circ L)(a))
\]

so the image of $X$ in $A\oplus A$ is the graph of the composition of $L$ and
$S$. More generally let $L\subset A\oplus A$ and $S\subset A\oplus A$ be two
linear subspaces and let $r:L\longrightarrow A\oplus A$ and
$s:S\longrightarrow A\oplus A$ be the corresponding relations with $r$ and $s$
the inclusion maps. Then the image of the product relation of $L$ and $S$ in
$A\oplus A$ is formed by selecting vectors in $u\in L$ , decomposing them as
$u=a+b$ with $a,b\in A$, selecting vectors $v\in S$ with decomposition $v=b+c
$ and finally forming the vectors $w=a+c$.

A monoid in the category of relations is a relation $r$ and two morphisms of
relations
\begin{align*}
\mu_{r}  & :r\odot^{A}r\longrightarrow r\\
u_{r}  & :\delta_{A}\longrightarrow r
\end{align*}

such that the associativity and unit diagrams commute. Let us consider the
case when the basic category is $Sets$. Then we have seen that a relation is a
graph with vertex set $A$ and arrow set $B$. The source and target for any
arrow $x$ is given by $f(x),g(x)\in A$ where $r(x)=(f(x),g(x))$. The domain
$X$ for the relation $r\odot^{A}r$ consists of all composable pairs of arrows
from the graph $B$.
\[
X=\{(x,x^{\prime})|\;g(x)=f(x^{\prime})\;\}
\]

The map $\mu_{r}$ will define a associative rule of composition for composable
pair of arrows in $B$. Furthermore the unitmap will provide for each vertex
$a\in A$ an arrow $u(a)$ that has $a$ as both source and target and that acts
as left and right unit for composition. The structure we have described is
clearly a internal category in $Sets$ with objects $A$ and arrows $B$. Let
$B\subset A\times A$ be a transitive and reflexive relation in the usual
sense. If we define $r((x,y))=(x,y)$ then $r$ is clearly a relation in our
sense. We have seen that the domain of the relation $r\otimes^{A}r$ is of the
form
\[
X=\{((x,y),(y,z))\;|\;(x,y)\in B,(y,z)\in B\}
\]

and $(r\otimes^{A}r)(((x,y),(y,z)))=(x,z)$. Define a map of sets $\mu
_{r}:X\longrightarrow A\times A$ by $\mu_{r}(((x,y),(y,z)))=(x,z)$. But both
$(x,y)\in B$ and $(y,z)\in B$ and since $B$ is transitive we have $(x,z)\in B$
and so we have in fact $\mu_{r}:X\longrightarrow B$. But we also have
\begin{align*}
& (r\circ\mu_{r})(((x,y),(y,z))\\
& =r((x,z))\\
& =(x,z)\\
& =(r\otimes^{A}r)(((x,y),(y,z)))
\end{align*}

so $\mu_{r}:r\otimes^{A}r\longrightarrow r$. The map $\mu_{r}$ is clearly
associative. Define a map of sets $u_{r}:A\longrightarrow A\times A$ by
$u_{r}(x)=(x,x)$. Since the relation $B$ is reflexive we have $(x,x)\in B$ for
all $x\in A$ and therefore we have $u_{r}:A\longrightarrow B$. This map is
clearly a morphism of relations and acts as a left and right unit for the rule
of composition $\mu_{r}$. We therefore have proved that the relation in our
sense,corresponding to a reflexive and transitive relation in the usual sense,
is in fact a monoid in the category of relations. We can thus think of monoids
in $\mathcal{R}^{A}(C)$ as generalized reflexive and transitive relations or
generalized categories.

\subsection{Symmetries for the category of relations}

>From an algebraic point of view we know that commutative monoids is an
important and interesting subclass of all monoids. From a categorical point of
view the notion of commutativity can not be formulated unless there is a
symmetry defined on the category.

Let us therefore consider the notion of a symmetry for the
category of relations. In $Sets$ a relation $r$ with domain $B$ is
a labelled and directed graph with vertex set $A$ and arrow set
$B$. We have seen that the tensor product of two relations $r$ and
$s$ with domains $B$ and $E$ is a new graph on the vertex set $A$
where the set of arrows consists of all composable pairs of arrows
from $B$ and $E$. If $(x,y)$ with $x\in B$ and $y\in E$ is a
composable pair of arrows it is clear that in general the pair
$(y,x)$ is not a composable pair. It is thus evident that for
relations the simple transposition $(x,y)\longrightarrow(y,x)$ is
not the right notion for a symmetry. We will develop our theory
for the category $\mathcal{S}^{A}(C)$ and use the isomorphism
whenever we need the corresponding structures in the
category of relations. Dual properties holds for the categories $\mathcal{S}%
_{A}(C)$ and the category of $C$-corelations.

Before we give the right definition of symmetry for the category of relations
we need to introduce a new structure. Let $\delta=\{\delta^{l},\delta^{r}\}$
be a object in $\mathcal{S}^{A}(C)$. Let $r=\Psi(\delta)$ be the corresponding
relation and define $r^{\ast}=\sigma_{A,A}\circ r$ and
\[
\delta^{\ast}=\Phi(r^{\ast})
\]

An explicit expression for the new object $\delta^{\ast}$ is given by the following.

\begin{proposition}
Let $\delta$ be a object in $\mathcal{S}^{A}(C)$ with underlying object $B$.
Then we have
\begin{align*}
(\delta^{\ast})^{l}  & =\sigma_{B,A}\circ\delta^{r}\\
(\delta^{\ast})^{r}  & =\sigma_{A,B}\circ\delta^{l}%
\end{align*}
\end{proposition}

\begin{proof}
Let $r=\Psi(\delta)$. Then we have
\begin{align*}
& (\delta^{\ast})^{l}\\
& =(\gamma_{A}\otimes1_{B})\circ((1_{A}\otimes\epsilon_{A})\otimes1_{B}%
)\circ(r^{\ast}\otimes1_{B})\circ\delta_{B}\\
& =(\gamma_{A}\otimes1_{B})\circ((1_{A}\otimes\epsilon_{A})\otimes1_{B}%
)\circ(\sigma_{A,A}\otimes1_{B})\circ(r\otimes1_{B})\circ\delta_{B}\\
& =(\gamma_{A}\otimes1_{B})\circ(\sigma_{e,A}\otimes1_{B})\circ((\epsilon
_{A}\otimes1_{A})\otimes1_{B})\circ(r\otimes1_{B})\circ\delta_{B}\\
& =(\beta_{A}\otimes1_{B})\circ((\epsilon_{A}\otimes1_{A})\otimes1_{B}%
)\circ(r\otimes1_{B})\circ\sigma_{B,B}\circ\delta_{B}\\
& =\sigma_{B,A}\circ(1_{B}\otimes\beta_{A})\circ(1_{B}\otimes(\epsilon
_{A}\otimes1_{A}))\circ(1_{B}\otimes r)\circ\delta_{B}\\
& =\sigma_{B,A}\circ\delta^{r}%
\end{align*}

where we have used the commutativity of $\delta_{B}$. In a similar way we show
that $(\delta^{\ast})^{r}=\sigma_{A,B}\circ\delta^{l}$.
\end{proof}

Note that $\delta$ and $\delta^{\ast}$ both have the same underlying object.
In order to extend the new operation $^{\ast}$ to morphisms we need the
following lemma

\begin{lemma}
Let $\delta$ and $\gamma$ be two objects in $\mathcal{S}^{A}(C)$ with
underlying objects $B$ and $E$ and let $f:\delta\longrightarrow\gamma$ be a
morphism. Then the corresponding arrow in $C$ define a morphism in
$\mathcal{S}^{A}(C)$ with domain $\delta^{\ast}$ and codomain $\gamma^{\ast}$.
\end{lemma}

\begin{proof}
We have
\begin{align*}
& (1_{A}\otimes f)\circ(\delta^{\ast})^{l}\\
& =(1_{A}\otimes f)\circ\sigma_{B,A}\circ\delta^{r}\\
& =\sigma_{E,A}\circ(f\otimes1_{A})\circ\delta^{r}\\
& =\sigma_{E,A}\circ\gamma^{r}\circ f\\
& =(\gamma^{\ast})^{l}\circ f
\end{align*}

and in a similar way we prove that $(f\otimes1_{A})\circ(\delta^{\ast}%
)^{r}=(\gamma^{\ast})^{r}\circ f$.
\end{proof}

We define $f^{\ast}:\delta^{\ast}\longrightarrow\gamma^{\ast}$ to be the
morphism described in the previous lemma. It is evident that the map
$T:\mathcal{S}^{A}(C)\longrightarrow\mathcal{S}^{A}(C)$ defined on objects and
arrows by $T(\delta)=\delta^{\ast}$ and $\ T(f)=f^{\ast}$ is a endofunctor on
$\mathcal{S}^{A}(C)$ and since $\sigma$ is a symmetry we have $T\circ
T=1_{\mathcal{S}^{A}(C)}$. This show that the category $\mathcal{S}^{A}(C)$
has a nontrivial action by the group $S_{2}=\langle t\;|\;t^{2}=1\rangle$.

The action have at least one fixpoint

\begin{proposition}
\label{fixunit}Let $a=\{\delta_{A},\delta_{A}\}$ be unit the object for the
monoidal structure $\otimes^{A}$on $\mathcal{S}^{A}(C)$ Then we have $a^{\ast}=a$
\end{proposition}

\begin{proof}
Recall that $\delta_{A}:A\longrightarrow A\otimes A$ defines a commutative
coalgebra structure on $A$. But then we have
\begin{align*}
& a^{\ast}\\
& =\{\delta_{A},\delta_{A}\}^{\ast}\\
& =\{\sigma_{A,A}\circ\delta_{A},\sigma_{A,A}\circ\delta_{A}\}\\
& =\{\delta_{A},\delta_{A}\}\\
& =a
\end{align*}
\end{proof}

The nontrivial group of symmetries must be taken into account when
the notion of a symmetry for the product structures
$\boxtimes^{A}$and $\otimes^{A}$ in $\mathcal{S}^{A}(C)$ are
defined. We have previously shown that for a symmetric monoidal
category in the usual sense we have an interpretation of the
Yang-Baxter equation and the unit symmetry conditions in terms of
invariance with respect to the group
$H=\{T_{t}(\sigma),1_{\otimes}\}$. This whole construction was
based on a certain choice of action by the group $S_{2}=\{1,t\}$
on the category $C$, $C^{2}$ and $C^{3}$ generated by the functors
$T_{1}=1_{C}$,$T_{2}=\tau$ and $T_{3}=(1_{C}\times\tau)\circ
(\tau\times1_{C})\circ(1_{C}\times\tau)$. What is new for the
category of relations is that we have a nontrivial action of
$S_{2}$. We will generalize this and consider monoidal categories
$\langle C,\otimes,K_{e},\alpha ,\beta,\gamma\rangle$ where we
have a nontrivial action of $S_{2}$ generated by a functor
$T_{1}:C\longrightarrow C$. We use this action together with
$\tau$ to define actions of $S_{2}$ on the categories $C^{2}$ and
$C^{3}$ generated by $T_{2}=\tau\circ(T_{1}\times T_{1})$ and
$T_{3}=(T_{1}\times T_{2})\circ(T_{2}\times
T_{1})\circ(T_{1}\times T_{2})$. It is easy to see that
$T_{2}\circ T_{2}=1_{C^{2}}$ and $T_{3}\circ T_{3}=1_{C^{3}}$ so
that these functors really defines an action of $S_{2}$. Note that
if $T_{1}=1_{C}$ we get the action we discussed previously in the
section on symmetries and group action. We now lift this action to
the functor categories $[C^{2},C]$ and $[C^{3},C]$ in the usual
way. From this point we proceed in a way that is exactly parallel
to what we did in the section on symmetries and group action. In
general one could imagine that the functor $T_{1}$ does not fix
the unit so that $T_{1}(e)\neq e$. In this general situation we
would assume the existence of a natural isomorphism
$\theta:K_{e}\longrightarrow tK_{e}$ in addition to the
isomorphism $\sigma:\otimes\longrightarrow t\otimes$. We would
thus allow the constant functor $K_{e}$ to be fixed only up to
natural isomorphism. In this paper we will not consider such a
possibility. This is because the unit is fixed both for the usual
case with trivial action and for the case of the category of
relations $S^{A}(C)$ as proved in proposition \ref{fixunit}.
Allowing the unit to move would also make all formulas and
derivations more complicated. With this out of the way we can now
state that all results derived in the section on symmetries and
group actions,up to and including corollary \ref{symcor} ,also
holds for the current situation if we substitute the $S_{2}$
action from this section in all statements. The proofs of these
results are of course different since they must take into account
the more general action considered in this section. We do not
reproduce these proofs here since they are long and rather similar
to the ones already given in the section on symmetries and group
action.

We now reverse proposition \ref{symdefprop} that characterized symmetric
monoidal categories in terms of invariance and is lead to the following
definition of symmetries for monoidal categories with a action of the group
$S_{2}$.

\begin{definition}
\label{mysym}Let $C$ be a category where there is defined an action of the
group $S_{2}$. Then $\langle C,\otimes,K_{e},\alpha,\beta,\gamma,\sigma
\rangle$ is a symmetric monoidal category if $\langle C,\otimes,K_{e}%
,\alpha,\beta,\gamma\rangle$ is a monoidal category and $\sigma:\otimes
\longrightarrow t\otimes$ is a natural isomorphism such that the following
identities holds.
\begin{align*}
\sigma\circ(1_{1_{C}}\times\sigma)  & =(t\alpha)\cdot(\sigma\circ(\sigma
\times1_{1_{C}}))\cdot\alpha\\
\beta & =(t\gamma)\cdot(\sigma\circ1_{K_{e}\times1_{C}})\\
\gamma & =(t\beta)\cdot(\sigma\circ1_{1_{C}\times K_{e}})\\
t\sigma & =\sigma^{-1}%
\end{align*}

We say that $\sigma$ is the symmetry for the monoidal category $\langle
C,\otimes,K_{e},\alpha,\beta,\gamma\rangle$.
\end{definition}

The first condition is equivalent to the Yang-Baxter equation if
we consider symmetric monoidal categories in the usual sense with
trivial action of $S_{2}$. We will in all cases call the first
condition for the Yang-Baxter equation. Note that even for the
case of trivial action our notion of symmetric monoidal category
is more general than the standard one. The standard definition of
symmetry for a monoidal category implies that the Yang-Baxter
equation holds but the fact that the Yang-Baxter equation holds
for $\sigma$ does not necessarily imply that $\sigma$ is a
symmetry in the usual sense.

We say that $\langle C,\otimes,\alpha,\sigma\rangle$ is a symmetric
semimonoidal category if $\langle C,\otimes,\alpha\rangle$ is a semimonoidal
category and the first and last of the above conditions hold.

We will now apply the definition of symmetry for the case of the category of
relations. For this case we denoted the map $T_{1}$ by $^{\ast}$. Let us first
consider the structure $\boxtimes^{A}$. \ In terms of objects, the definition
of a symmetry for the semimonoidal category $\langle\mathcal{S}^{A}%
(C),\boxtimes^{A},M^{A}\rangle$ is as follows.

\begin{definition}
A symmetry for the semimonoidal category $\langle\mathcal{S}^{A}%
(C),\boxtimes^{A},M^{A}\rangle$ is a isomorphism $S_{\delta,\gamma}^{A}%
:\delta\boxtimes^{A}\gamma\longrightarrow(\gamma^{\ast}\boxtimes^{A}%
\delta^{\ast})^{\ast}$ that is natural in $\delta$ and $\gamma$ and such that
the following identities are satisfied for all $\delta,\gamma$ and $\rho$.
\begin{align*}
(M_{\rho^{\ast},\gamma^{\ast},\delta^{\ast}}^{A})^{\ast}\circ(1_{\rho}^{\ast
}\boxtimes^{A}(S_{\delta,\gamma}^{A})^{\ast})^{\ast}\circ S_{\delta
\boxtimes^{A}\gamma,\rho}^{A}\circ M_{\delta,\gamma,\rho}^{A}  &
=((S_{\gamma,\rho}^{A})^{\ast}\boxtimes^{A}1_{\delta}^{\ast})^{\ast}\circ
S_{\delta,\gamma\boxtimes^{A}\rho}^{A}\\
(S_{\gamma^{\ast},\delta^{\ast}}^{A})^{\ast}\circ S_{\delta,\gamma}^{A}  &
=1_{\delta\boxtimes^{A}\gamma}%
\end{align*}
\end{definition}

In general many symmetric semimonoidal structures may exist for $\mathcal{S}%
^{A}(C)$ with the product $\boxtimes^{A}$. We will now show there is always at
least one.

\begin{proposition}
Let $\delta,\gamma$ and $\rho$ be objects in $\mathcal{S}^{A}(C)$ with
underlying objects $B$,$E$ and $D$. Define
\begin{align*}
M_{\delta,\gamma,\rho}^{A}  & =\alpha_{B,E,D}\\
S_{\delta,\gamma}^{A}  & =\sigma_{B,E}%
\end{align*}

where $\alpha$ and $\sigma$ are the associativity constraint and symmetry for
the category $C$. Then $\langle\mathcal{S}^{A}(C),\boxtimes^{A},M^{A}%
,S^{A}\rangle$ is a symmetric semimonoidal category.
\end{proposition}

\begin{proof}
We already know that $\langle\mathcal{S}^{A}(C),\boxtimes^{A},M^{A}\rangle$ is
a semimonoidal category. First we need to prove that $S_{\delta,\gamma}^{A}$
is a morphism in $\mathcal{S}^{A}(C)$. Note that the underlying object for
$(\gamma^{\ast}\boxtimes^{A}\delta^{\ast})^{\ast}$ is $E\otimes B$. If we use
the fact that $\sigma$ is a symmetry in $C$ we have
\begin{align*}
& ((\gamma^{\ast}\boxtimes^{A}\delta^{\ast})^{\ast})^{l}\circ S_{\delta
,\gamma}^{A}\\
& =\sigma_{E\otimes B,A}\circ(\gamma^{\ast}\boxtimes^{A}\delta^{\ast}%
)^{r}\circ\sigma_{B,E}\\
& =\sigma_{E\otimes B,A}\circ\alpha_{E,B,A}\circ(1_{E}\otimes(\delta^{\ast
})^{r})\circ\sigma_{B,E}\\
& =\sigma_{E\otimes B,A}\circ\alpha_{E,B,A}\circ(1_{E}\otimes\sigma_{A,B}%
\circ\delta^{l})\circ\sigma_{B,E}\\
& =\sigma_{E\otimes B,A}\circ\alpha_{E,B,A}\circ(1_{E}\otimes\sigma
_{A,B})\circ(1_{E}\otimes\delta^{l})\circ\sigma_{B,E}\\
& =\sigma_{E\otimes B,A}\circ\alpha_{E,B,A}\circ(1_{E}\otimes\sigma
_{A,B})\circ\sigma_{A\otimes B,E}\circ(\delta^{l}\otimes1_{E})\\
& =(1_{A}\otimes\sigma_{B,E})\circ\alpha_{A,B,E}^{-1}\circ(\delta^{l}%
\otimes1_{E})\\
& =(1_{A}\otimes S_{\delta,\gamma}^{A})\circ(\delta\boxtimes^{A}\gamma)^{l}%
\end{align*}

and this proves that $S_{\delta,\gamma}^{A}$ is a morphism in $\mathcal{S}%
^{A}(C)$. It is clearly a isomorphism and naturality is evident.
The condition for $S^{A}$ to be a symmetry in $\mathcal{S}^{A}(C)$
is satisfied since it turns into the condition for $\sigma$ being
is a symmetry in the category $C$.
\end{proof}

The previous proposition leads us to make the following definition.

\begin{definition}
A symmetric semimonoidal structure $\langle\mathcal{S}^{A}(C),\boxtimes
^{A},M^{A},S^{A}\rangle$ on the category $\mathcal{S}^{A}(C)$ is external if
for all objects $\delta,\gamma$ and $\rho$ we have
\begin{align*}
M_{\delta,\gamma,\rho}^{A}  & =\alpha_{B,E,D}\\
S_{\delta,\gamma}^{A}  & =\sigma_{B,E}%
\end{align*}

where $B$,$E$ and $D$ are the underlying objects for $\delta,\gamma$ and
$\rho$ and where $\alpha$ and $\sigma$ are the associativity constraint and
symmetry for the category $C$.
\end{definition}

The previous proposition then proves that an external symmetries on
$\mathcal{S}^{A}(C)$ with product $\boxtimes^{A}$always exists.

We now turn to the definition of symmetries for $\mathcal{S}^{A}(C)$ with the
product $\otimes^{A}$. In terms of objects the general definition now takes
the form

\begin{definition}
A symmetry for the monoidal category$\langle\mathcal{S}^{A}(C),\otimes
^{A},a,m^{A},l^{A},r^{A}\rangle$ is a isomorphism $s_{\delta,\gamma}%
^{A}:\delta\otimes^{A}\gamma\longrightarrow(\gamma^{\ast}\otimes^{A}%
\delta^{\ast})^{\ast}$ that is natural in $\delta$ and $\gamma$ and such that
the following identities are satisfied for all $\delta,\gamma$ and $\rho$.
\begin{align*}
(m_{\rho^{\ast},\gamma^{\ast},\delta^{\ast}}^{A})^{\ast}\circ(1_{\rho}^{\ast
}\otimes^{A}(s_{\delta,\gamma}^{A})^{\ast})^{\ast}\circ s_{\delta\otimes
^{A}\gamma,\rho}^{A}\circ m_{\delta,\gamma,\rho}^{A}  & =((s_{\gamma,\rho}%
^{A})^{\ast}\otimes^{A}1_{\delta}^{\ast})^{\ast}\circ s_{\delta,\gamma
\otimes^{A}\rho}^{A}\\
(l_{\delta^{\ast}}^{A})^{\ast}\circ s_{\delta,a}^{A}  & =r_{\delta}^{A}\\
(r_{\delta^{\ast}}^{A})^{\ast}\circ s_{a,\delta}^{A}  & =l_{\delta}^{A}\\
(s_{\delta^{\ast},\gamma^{\ast}}^{A})^{\ast}\circ s_{\delta,\gamma}^{A}  &
=1_{\delta\otimes^{A}\gamma}%
\end{align*}
\end{definition}

Note that identity two and three are not independent. One can be derived from
the other by using identity four and the fact that the neutral object $a $ is
fixed by the action of $S_{2}$. There are several equivalent formulations of
the first symmetry condition

\begin{proposition}
Let $s^{A}$ be a natural isomorphism $s_{\delta,\gamma}^{A}:\delta\otimes
^{A}\gamma\longrightarrow(\gamma^{\ast}\otimes^{A}\delta^{\ast})$ such that
the following identities hold
\begin{align*}
(l_{\delta^{\ast}}^{A})^{\ast}\circ s_{\delta,a}^{A}  & =r_{\delta}^{A}\\
(s_{\delta^{\ast},\gamma^{\ast}}^{A})^{\ast}\circ s_{\delta,\gamma}^{A}  &
=1_{\delta\otimes^{A}\gamma}%
\end{align*}

Then the following statements are equivalent

\begin{enumerate}
\item $s^{A}$ is a symmetry

\item $(m_{\rho^{\ast},\gamma^{\ast},\delta^{\ast}}^{A})^{\ast}\circ(1_{\rho
}^{\ast}\otimes^{A}(s_{\delta,\gamma}^{A})^{\ast})^{\ast}\circ s_{\delta
\otimes^{A}\gamma,\rho}^{A}\circ m_{\delta,\gamma,\rho}^{A}=s_{\delta
,(\rho^{\ast}\otimes^{A}\gamma^{\ast})^{\ast}}^{A}\circ(1_{\delta}\otimes
^{A}s_{\gamma,\rho}^{A})$

\item $(m_{\rho^{\ast},\gamma^{\ast},\delta^{\ast}}^{A})^{\ast}\circ
s_{(\gamma^{\ast}\otimes^{A}\delta^{\ast})^{\ast},\rho}^{A}\circ
(s_{\delta,\gamma}^{A}\otimes^{A}1_{\rho})\circ m_{\delta,\gamma,\rho}%
^{A}=((s_{\gamma,\rho}^{A})^{\ast}\otimes^{A}1_{\delta}^{\ast})^{\ast}\circ
s_{\delta,\gamma\otimes^{A}\rho}^{A}$

\item $(m_{\rho^{\ast},\gamma^{\ast},\delta^{\ast}}^{A})^{\ast}\circ
s_{(\gamma^{\ast}\otimes^{A}\delta^{\ast})^{\ast},\rho}^{A}\circ
(s_{\delta,\gamma}^{A}\otimes^{A}1_{\rho})\circ m_{\delta,\gamma,\rho}%
^{A}=s_{\delta,(\rho^{\ast}\otimes^{A}\gamma^{\ast})^{\ast}}^{A}%
\circ(1_{\delta}\otimes^{A}s_{\gamma,\rho}^{A})$
\end{enumerate}
\end{proposition}

\begin{proof}
By naturality of $s^{A}$ we have the following two identities
\begin{align*}
((s_{\gamma,\rho}^{A})^{\ast}\otimes^{A}1_{\delta}^{\ast})^{\ast}\circ
s_{\delta,\gamma\otimes^{A}\rho}^{A}  & =s_{\delta,(\rho^{\ast}\otimes
^{A}\gamma^{\ast})^{\ast}}^{A}\circ(1_{\delta}\otimes^{A}s_{\gamma,\rho}%
^{A})\\
(1_{\rho}^{\ast}\otimes^{A}(s_{\delta,\gamma}^{A})^{\ast})^{\ast}\circ
s_{\delta\otimes^{A}\gamma,\rho}^{A}  & =s_{(\gamma^{\ast}\otimes^{A}%
\delta^{\ast})^{\ast},\rho}^{A}\circ(s_{\delta,\gamma}^{A}\otimes^{A}1_{\rho})
\end{align*}

The proposition now follows directly from these identities.
\end{proof}

Let us now consider the existence of symmetries.

\begin{definition}
A symmetry for the monoidal category $\langle\mathcal{S}^{A}(C),\otimes
^{A},a,m^{A},l^{A},r^{A}\rangle$ is induced by a symmetry $S^{A}$ of the
semimonoidal category $\langle\mathcal{S}^{A}(C),\boxtimes^{A},M^{A}\rangle$
if for all objects $\delta,\gamma$ in $\mathcal{S}^{A}(C)$ the following
diagram commute%

{\tiny
\begin{diagram}
\delta\boxtimes^A \gamma&\rTo^{S_{\delta,\gamma}^A}%
&                       &(\gamma^*\boxtimes^A\delta^*)^*        \\
\uTo^{\pi_{\delta,\gamma}^A}%
&                                  &                         & \uTo
_{(\pi_{\gamma^*,\delta^*}^A)^*}   \\
\delta\otimes^A \gamma& \rTo_{s_{\delta,\gamma}^A}%
&                         &(\gamma^*\otimes^A \delta^*)^*            \\
\end{diagram}
}%
\end{definition}

We will in the following show that an induced symmetry exists in the external
case and is uniquely determined by the symmetry $S^{A}$.

Recall that for any pair of objects $\delta$ and $\gamma$ in $\mathcal{S}%
^{A}(C)$, the diagram $\mathcal{P}_{\delta,\gamma}^{A}$ was given by%

{\tiny
\begin{diagram}
\delta\boxtimes^A(a\boxtimes^A\gamma
)  &                                                      &   \rTo
^{M^A_{\delta,a,\gamma}}%
&                                                    &    (\delta
\boxtimes^A a)\boxtimes\gamma\\
& \luTo_{1_{\delta}\boxtimes^A\gamma^l}%
&                                         & \ruTo_{\delta^r\boxtimes
^A 1_{\gamma}}  &                                                     \\
&                                                      &       \delta
\boxtimes\gamma
&                                                     &                                                     \\
\end{diagram}
}%

>From the general theory of categories it is well known that isomorphisms of
categories preserve universal cones. By definition $\langle\gamma^{\ast
}\otimes^{A}\delta^{\ast},\pi_{\gamma^{\ast},\delta^{\ast}}^{A}\rangle$ is a
universal cone on the diagram $\mathcal{P}_{\gamma^{\ast},\delta^{\ast}}^{A}%
$and therefore $\langle(\gamma^{\ast}\otimes^{A}\delta^{\ast})^{\ast}%
,(\pi_{\gamma^{\ast},\delta^{\ast}}^{A})^{\ast}\rangle$ is a universal cone on
the diagram $(\mathcal{P}_{\gamma^{\ast},\delta^{\ast}}^{A})^{\ast}$. But we
have the following result

\begin{lemma}
Let the symmetric semimonoidal category $\langle\mathcal{S}^{A}(C),\boxtimes
^{A},M^{A},S^{A}\rangle$ be external, then $\langle\delta\otimes^{A}%
\gamma,S_{\delta,\gamma}^{A}\circ\pi_{\delta,\gamma}^{A}\rangle$ is a
universal cone on $(\mathcal{P}_{\gamma^{\ast},\delta^{\ast}}^{A})^{\ast}$.
\end{lemma}

\begin{proof}
Let us first prove that it is a cone. For this we must prove that the
following identity holds
\[
(M_{\gamma^{\ast},a,\delta^{\ast}}^{A})^{\ast}\circ(1_{\gamma^{\ast}}%
\boxtimes^{A}(\delta^{\ast})^{l})^{\ast}\circ S_{\delta,\gamma}^{A}\circ
\pi_{\delta,\gamma}^{A}=((\gamma^{\ast})^{r}\boxtimes^{A}1_{\delta^{\ast}%
})^{\ast}%
\]

holds in $\mathcal{S}^{A}(C)$. Since the semimonoidal structure on
$\mathcal{S}^{A}(C)$ is external, the previous identity is for the strict case
equivalent to the following identity in $C$%
\[
\alpha_{E,A,B}\circ(1_{E}\otimes(\delta^{\ast})^{l})\circ\sigma_{B,E}\circ
\pi_{\delta,\gamma}^{A}=((\gamma^{\ast})^{r}\otimes1_{B})\circ\sigma
_{B,E}\circ\pi_{\delta,\gamma}^{A}%
\]

But this identity follows from the Yang baxter equation and the fact that
$\langle\delta\otimes^{A}\gamma,\pi_{\delta,\gamma}^{A}\rangle$ is a cone on
$\mathcal{P}_{\delta,\gamma}^{A}$.
\begin{align*}
& =\alpha_{E,A,B}\circ(1_{\gamma^{\ast}}\boxtimes^{A}(\delta^{\ast}%
)^{l})^{\ast}\circ\sigma_{B,E}\circ\pi_{\delta,\gamma}^{A}\\
& =\alpha_{E,A,B}\circ(1_{E}\otimes\sigma_{B,A})\circ(1_{E}\otimes\delta
^{r})\circ\sigma_{B,E}\circ\pi_{\delta,\gamma}^{A}\\
& =\alpha_{E,A,B}\circ(1_{E}\otimes\sigma_{B,A})\circ\sigma_{B\otimes
A,E}\circ(\delta^{r}\otimes1_{E})\circ\pi_{\delta,\gamma}^{A}\\
& =\alpha_{E,A,B}\circ(1_{E}\otimes\sigma_{B,A})\circ\sigma_{B\otimes
A,E}\circ\alpha_{B,A,E}\circ(1_{B}\otimes\gamma^{l})\circ\pi_{\delta,\gamma
}^{A}\\
& =\alpha_{E,A,B}\circ(1_{E}\otimes\sigma_{B,A})\circ\sigma_{B\otimes
A,E}\circ\alpha_{B,A,E}\circ\sigma_{A\otimes E,B}\circ(\gamma^{l}\otimes
1_{B})\circ\sigma_{B,E}\circ\pi_{\delta,\gamma}^{A}\\
& =\alpha_{E,A,B}\circ(1_{E}\otimes\sigma_{B,A})\circ\sigma_{B\otimes
A,E}\circ\alpha_{B,A,E}\circ\sigma_{A\otimes E,B}\circ(\sigma_{E,A}%
\otimes1_{B})\\
& \circ((\gamma^{\ast})^{r}\otimes1_{B})\circ\sigma_{B,E}\circ\pi
_{\delta,\gamma}^{A}\\
& =((\gamma^{\ast})^{r}\otimes1_{B})\circ\sigma_{B,E}\circ\pi_{\delta,\gamma
}^{A}%
\end{align*}

Let now $\langle\theta,u\rangle$ be any cone on $(\mathcal{P}_{\gamma^{\ast
},\delta^{\ast}}^{A})^{\ast}$. The proposition is proved if we can show that
the following equations has a unique solution $\varphi:\theta\longrightarrow
\delta\otimes^{A}\gamma$%
\[
S_{\delta,\gamma}^{A}\circ\pi_{\delta,\gamma}^{A}\circ\varphi=u
\]

The equation has at most one solution since $S_{\delta,\gamma}^{A}\circ
\pi_{\delta,\gamma}^{A}$ is a monomorphism. In a calculation very similar to
previous one we can prove that $\langle\theta,(S_{\delta,\gamma}^{A}%
)^{-1}\circ u\rangle$ is a cone on $\mathcal{P}_{\delta,\gamma}^{A}$. But
$\langle\delta\otimes^{A}\gamma,\pi_{\delta,\gamma}^{A}\rangle$ is a universal
cone on $\mathcal{P}_{\delta,\gamma}^{A}$ and therefore there exists a
morphism $h:\theta\longrightarrow\delta\otimes^{A}\gamma$ in $\mathcal{S}%
^{A}(C)$ such that $\pi_{\delta,\gamma}^{A}\circ h=(S_{\delta,\gamma}%
^{A})^{-1}\circ u$. Composing on both sides with $S_{\delta,\gamma}^{A}$ show
that $S_{\delta,\gamma}^{A}\circ\pi_{\delta,\gamma}^{A}\circ\varphi=u$ and the
proposition is proved.
\end{proof}

We can now prove the existence of induced symmetries in the external case.

\begin{theorem}
Let the symmetric semimonoidal category $\langle\mathcal{S}^{A}(C),\boxtimes
^{A},M^{A},S^{A}\rangle$ be external, then there exists a induced symmetry for
the monoidal category $\langle\mathcal{S}^{A}(C),\otimes^{A},a,m^{A}%
,l^{A},r^{A}\rangle$.
\end{theorem}

\begin{proof}
The previous lemma show that both $\langle\delta\otimes^{A}\gamma
,S_{\delta,\gamma}^{A}\circ\pi_{\delta,\gamma}^{A}\rangle$ and $\langle
(\gamma^{\ast}\otimes^{A}\delta^{\ast})^{\ast},(\pi_{\gamma^{\ast}%
,\delta^{\ast}}^{A})^{\ast}\rangle$ are universal cones on $(\mathcal{P}%
_{\gamma^{\ast},\delta^{\ast}}^{A})^{\ast}$. We can therefore conclude that
there exists a unique morphism $s_{\delta,\gamma}^{A}:\delta\otimes^{A}%
\gamma\longrightarrow(\gamma^{\ast}\otimes^{A}\delta^{\ast})^{\ast}$ such that
$(\pi_{\gamma^{\ast},\delta^{\ast}}^{A})^{\ast}\circ s_{\delta,\gamma}%
^{A}=S_{\delta,\gamma}^{A}\circ\pi_{\delta,\gamma}^{A}$. We will show that
$s_{\delta,\gamma}^{A}$ is a symmetry for the monoidal category $\langle
\mathcal{S}^{A}(C),\otimes^{A},a,m^{A},l^{A},r^{A}\rangle$ on $\mathcal{S}%
^{A}(C)$. The first symmetry condition for $s^{A}$ follows from the first
symmetry condition for $S^{A}$, the identity $(\pi_{\gamma^{\ast},\delta
^{\ast}}^{A})^{\ast}\circ s_{\delta,\gamma}^{A}=S_{\delta,\gamma}^{A}\circ
\pi_{\delta,\gamma}^{A}$ and the fact that $m^{A}$ is induced by $M^{A}$. For
the second symmetry condition we have for the strict case
\begin{align*}
& (l_{\delta^{\ast}}^{A})^{\ast}\circ s_{\delta,a}^{A}\\
& =\beta_{B}\circ(\epsilon_{A}\otimes1_{B})\circ\pi_{a,\delta^{\ast}}^{A}\circ
s_{\delta,a}^{A}\\
& =\beta_{B}\circ(\epsilon_{A}\otimes1_{B})\circ S_{\delta,a}^{A}\circ
\pi_{\delta,a}^{A}\\
& =\beta_{B}\circ(\epsilon_{A}\otimes1_{B})\circ\sigma_{B,A}\circ\pi
_{\delta,a}^{A}\\
& =\beta_{B}\circ\sigma_{B,e}\circ(1_{B}\otimes\epsilon_{A})\circ\pi
_{\delta,a}^{A}\\
& =\gamma_{B}\circ(1_{B}\otimes\epsilon_{A})\circ\pi_{\delta,a}^{A}\\
& =r_{\delta}^{A}%
\end{align*}

where we have used the fact that the symmetry $S^{A}$ is external. The last
symmetry condition follows easily from the commutative diagram defining
$s^{A}$ in terms of $S^{A}$and from the fact that $S^{A}$ is a symmetry.
\end{proof}

\subsection{Commutative monoids in the category of relation}

We will define the notion of a commutative monoid for categories with an
action of $S_{2}$ and then apply this definition to the case of relations. Let
now $\langle C,\otimes,K_{e},\alpha,\beta,\gamma,\sigma\rangle$ be a symmetric
monoidal category with an action of $S_{2}$ generated by the functor
$T_{1}:C\longrightarrow C$. \ The conditions from definition \ref{symcat} thus
holds for $\alpha,\beta,\gamma$ and $\sigma$.

Our definition of a commutative monoid is a natural extension and
categorization of the notion of a commutative monoids in algebra. Let $\langle
M,\cdot,e\rangle$ be a monoid in the usual algebraic sense, so that $M$ is a
set and $\cdot$ is a associative product on $M$ with unit element $e$. Define
a new associative product on $M$ by $x\ast y=y\cdot x$. Then $\langle
M,\ast,e\rangle$ is a new monoid on the same underlying set. The monoid $M$ is
said to be commutative if $\langle M,\ast,e\rangle$ is the same monoid as
$\langle M,\cdot,e\rangle$ and this is equivalent to the condition $x\cdot
y=y\cdot x$ for all $x$ and $y$ in $M$. The previous condition is really too
strict since in algebra we consider isomorphic monoids to be essentially the
same. Thus it would be more natural to require that the two monoids $\langle
M,\cdot,e\rangle$ and $\langle M,\ast,e\rangle$ are isomorphic. From a
categorical point of view the last condition is the only one that really makes
sense since the relation of equality exists only between arrows and not
between objects. If we now recall that the symmetry $\sigma$ is the
categorization of the idea of changing order in the category $C$ we arrive at
our definition of commutativity.

Let $X$ be a monoid in the category $C$ with product $\mu:X\otimes
X\longrightarrow X$ and unit $u:e\longrightarrow X$. Define morphisms
$\mu^{\sigma}:T_{1}(X)\otimes T_{1}(X)\longrightarrow T_{1}(X)$ and
$u^{s}:e\longrightarrow T_{1}(X)$ by
\begin{align*}
\mu^{\sigma}  & =T_{1}(\mu)\circ\sigma_{T_{1}(X),T_{1}(X)}\\
u^{s}  & =T_{1}(u)
\end{align*}

\begin{proposition}
$\langle T_{1}(X),\mu^{\sigma},u^{s}\rangle$ is a monoid in $C$
\end{proposition}

\begin{proof}
The Yang-Baxter equation and the naturality of $\sigma$ implies when evaluated
on $(T_{1}(X),T_{1}(X),T_{1}(X))$ the following relation
\begin{align*}
\sigma_{T_{1}(X),T_{1}(X\otimes X)}\circ(1_{T_{1}(X)}\otimes\sigma
_{T_{1}(X),T_{1}(X)})  & =T_{1}(\alpha_{X,X,X})\circ\sigma_{T_{1}(X\otimes
X),T_{1}(X)}\\
& \circ(\sigma_{T_{1}(X),T_{1}(X)}\otimes1_{T_{1}(X)})
\end{align*}

Using this relation we have
\begin{align*}
& \mu^{\sigma}\circ(1_{T_{1}(X)}\otimes\mu^{\sigma})\\
& =T_{1}(\mu)\circ\sigma_{T_{1}(X),T_{1}(X)}\circ(1_{T_{1}(X)}\otimes
(T_{1}(\mu)\circ\sigma_{T_{1}(X),T_{1}(X)}))\\
& =T_{1}(\mu)\circ\sigma_{T_{1}(X),T_{1}(X)}\circ(1_{T_{1}(X)}\otimes
T_{1}(\mu))\circ(1_{T_{1}(X)}\otimes\sigma_{T_{1}(X),T_{1}(X)})\\
& =T_{1}(\mu)\circ T_{1}(\mu\otimes1_{X})\circ\sigma_{T_{1}(X),T_{1}(X\otimes
X)}\circ(1_{T_{1}(X)}\otimes\sigma_{T_{1}(X),T_{1}(X)})\\
& =T_{1}(\mu\circ(\mu\otimes1_{X}))\circ T_{1}(\alpha_{X,X,X})\circ
\sigma_{T_{1}(X\otimes X),T_{1}(X)}\circ(\sigma_{T_{1}(X),T_{1}(X)}%
\otimes1_{T_{1}(X)})\\
& =T_{1}(\mu)\circ T_{1}(1_{X}\otimes\mu)\circ\sigma_{T_{1}(X\otimes
X),T_{1}(X)}\circ(\sigma_{T_{1}(X),T_{1}(X)}\otimes1_{T_{1}(X)})\\
& =T_{1}(\mu)\circ\sigma_{T_{1}(X),T_{1}(X)}\circ((T_{1}(\mu)\circ
\sigma_{T_{1}(X),T_{1}(X)})\otimes1_{T_{1}(X)})\\
& =\mu^{\sigma}\circ(\mu^{\sigma}\otimes1_{T_{1}(X)})
\end{align*}

so the morphism $\mu^{\sigma}$ is associative. The first unit
condition evaluated at the pair of objects $(e,T_{1}(X))$ given
the identity
\[
\beta_{e,T_{1}(X)}=T_{1}(\gamma_{X,e})\circ\sigma_{e,T_{1}(X)}%
\]

>From the naturality of $\sigma$ and the fact that $\langle X,\mu,u\rangle$ is
a monoid we have
\begin{align*}
& \mu^{\sigma}\circ(u^{s}\otimes1_{T_{1}(X)})\\
& =T_{1}(\mu)\circ\sigma_{T_{1}(X),T_{1}(X)}\circ(T_{1}(u)\otimes1_{T_{1}%
(X)})\\
& =T_{1}(\mu)\circ T_{1}(1_{X}\otimes u)\circ\sigma_{e,T_{1}(X)}\\
& =T_{1}(\mu\circ(1_{X}\otimes u))\circ\sigma_{e,T_{1}(X)}\\
& =T_{1}(\gamma_{X,e})\circ\sigma_{e,T_{1}(X)}\\
& =\beta_{e,T_{1}(X)}%
\end{align*}

and this is the left condition on the unit. The proof for the right condition
is similar.
\end{proof}

Recall that $\varphi:X\longrightarrow Y$ is a morphism of monoids
$\varphi:\langle X,\mu,u\rangle\longrightarrow\langle Y,\mu^{\prime}%
,u^{\prime}\rangle$ if the following two diagrams commute%

{\tiny
\begin{equation*}
\begin{aligned}
\begin{diagram}
X\otimes X   &                               &\rTo^{\phi\otimes\phi}%
&                          &Y\otimes Y         \\
\dTo^{\mu}%
&                                &                                &                           &  \dTo
_{\mu'}   \\
X                          &                                &\rTo_{\phi}%
&                           &  Y
\end{diagram}
\end{aligned}
\;
\;
\;
\begin{aligned}
\begin{diagram}
X &                                                      &   \rTo^{\phi}%
&                                                    &    Y             \\
& \luTo_{u}                                      &                 & \ruTo
_{u'}                                   &                              \\
&                                                      &       e       &                                                     &                              \\
\end{diagram}
\end{aligned}
\end{equation*}
}%

We are now ready to define the notion of a commutative monoid in the symmetric
monoidal category $C$.

\begin{definition}
\label{commutativity}Let $\langle C,\otimes,K_{e},\alpha,\beta,\gamma
,\sigma\rangle$ be a symmetric monoidal category. A monoid $\langle
X,\mu,u\rangle$ in $C$ is commutative if there exists an isomorphism of
monoids
\[
\varphi:\langle X,\mu,u\rangle\longrightarrow\langle T_{1}(X),\mu^{\sigma
},u^{s}\rangle
\]
\end{definition}

We will now apply these definitions to the $Sets$. For this case
there is only one possible choice that makes $Sets$ into a
C-category. Let $r(x)=(f(x),g(x))$ and $s(y)=(h(y),k(y))$ be two
relations with domains $B$ and $E$. Then $r^{\ast}(x)=(g(x),f(x))$
and $s^{\ast}(y)=(k(y),h(y))$. If $X$ and $Y$ are the underlying
sets for $r\otimes^{A}s$ and $(s^{\ast}\otimes
^{A}r^{\ast})^{\ast}$ we have
\begin{align*}
X  & =\{(x,y)\;|\;g(x)=h(y)\}\\
Y  & =\{(y,x)\;|\;h(y)=g(x)\}
\end{align*}

and the relations are given by
\begin{align*}
(r\otimes^{A}s)(x,y)  & =(f(x),k(y))\\
(s^{\ast}\otimes^{A}r^{\ast})^{\ast}(y,x)  & =(f(x),k(y))
\end{align*}

Define $s_{r,s}^{A}(x,y)=(y,x)$. Then clearly we have $s_{r,s}^{A}%
:X\longrightarrow Y$ and also
\begin{align*}
& ((s^{\ast}\otimes^{A}r^{\ast})^{\ast}\circ s_{r,s}^{A})(x,y)\\
& =(s^{\ast}\otimes^{A}r^{\ast})^{\ast}(y,x)\\
& =(f(x),k(y))\\
& =(r\otimes^{A}s)(x,y)
\end{align*}

so that we have a morphism in
$s_{r,s}^{A}:r\otimes^{A}s\longrightarrow
(s^{\ast}\otimes^{A}r^{\ast})^{\ast}$. It is straight forward to
prove that $s^{A}$ is a symmetry on the category of relations. It
is in fact induced by the symmetry of the external category
$Sets$. Since a relation $r:B\longrightarrow A\times A$ is a
directed labelled graph it is clear that we get the relation
$r^{\ast}:B\longrightarrow A\times A$ by reversing all arrows in
the relation $r$. We have seen that $r$ is a monoid if there
exists a associative rule of composition for composable arrows in
$r$ such that for each object $x\in A$ there exists an arrow with
source and target given by $x$ and that acts as right and left
unit for the composition. Let $b $ and $b^{\prime}$ be two objects
in $B$. Then the rule of composition for the relation $r^{\ast}$is
defined by first reversing both arrows, then composing them as
arrows in $r$ and then reversing the result to get a arrow in
$r^{\ast}$. Now a isomorphism $\varphi:r\longrightarrow r^{\ast}$
is a bijective map with domain and codomain given by $B$ and such
that
\begin{align*}
g(\varphi(x))  & =f(x)\\
f(\varphi(x))  & =g(x)
\end{align*}

for all $x\in B$. If $\varphi$ is also a isomorphism of the monoids $\langle
r,\mu,u\rangle$ and $\langle r^{\ast},\mu^{s},u^{s}\rangle$ we must have
\[
\varphi(\mu((x,y)))=\mu((\varphi(y),\varphi(x)))
\]

for all objects $(x,y)\in A\times A$ such that $g(x)=f(y)$. These
conditions are in general impossible to satisfy for the identity
map $\varphi=1_{B}$. Let $B\subset A\times A$ be a relation in the
usual sense. Then $B$ corresponds to a relation in our sense if we
define $r:B\longrightarrow A\times A$ by $r((x,y))=(x,y)$ so that
$f((x,y))=x$ and $g((x,y))=y$. We know that if $r$ is a monoid in
the category of relations then $B$ is a reflexive and transitive
relation in the usual sense. \ Assume that $B$ is also symmetric
so that it is in fact a equivalence relation. Thus we have
$(x,y)\in B$ if and only if $(y,x)\in B$. Then we can define a map
$\varphi :B\longrightarrow B$ by $\varphi((x,y))=(y,x)$. For this
map we have
\begin{align*}
g(\varphi((x,y)))  & =g((y,x))=x=f((x,y))\\
f(\varphi((x,y)))  & =f((y,x))=y=g((x,y))
\end{align*}

so that $\varphi:r\longrightarrow r^{\ast}$. We have seen that the rule of
composition and unit maps for $r$ are given by
\begin{align*}
\mu_{r}(((x,y),(y,z)))  & =(x,z)\\
u_{r}(x)  & =(x,x)
\end{align*}

but then the composition and unit maps for $r^{\ast}$ must be given by
\begin{align*}
\mu^{s}(((y,x),(z,y)))  & =(z,x)\\
u^{s}(x)  & =(x,x)
\end{align*}

It is evident that $\varphi$ preserve that unit and the following computation
show that $\varphi:\langle r,\mu,u\rangle\longrightarrow\langle r^{\ast}%
,\mu^{s},u^{s}\rangle$ also preserve the product
\begin{align*}
& \varphi(\mu(((x,y),(y,z))))\\
& =\varphi((x,z))\\
& =(z,x)\\
& =\mu^{s}(((y,x),(z,y)))
\end{align*}

We have thus proved the following result

\begin{proposition}
Let $B\subset A\times A$ be a equivalence relation. Define $r:B\longrightarrow
A\times A$ by $r((x,y))=(x,y)$. Then $r$ is a commutative monoid in the
category of relations with respect to the symmetry in $\mathcal{R}^{A}(C)$
induced by the symmetry $\sigma(x,y)=(y,x)$ in $Sets$.
\end{proposition}

Note that this result show that relations that are not equivalence relations
in the usual sense might correspond to commutative monoids with respect to a
different symmetry than the standard one used in the proposition. Such a class
of relations would corresponds to a extension of the notion of equivalence
that might be of interest.

\section{Quantization of relations}

In this section we apply our ideas of quantization as properties
of functors in categories of representations of constraints. The
constraints here are the system of functors and natural
transformations defining a symmetric monoidal category where we
have a action of the group $S_{2}$. Morphisms in this category of
representations are what we call quantized functors. These are
determined by a functor and a triple of natural isomorphisms that
satisfy certain conditions that ensure that the functors behave in
a natural way with respect to the representations. \ Properties of
relations are coded in terms of commutative diagrams of arrows in
the category of relations. Equivalence relations appears as
commutative associative algebras with unit. In the last section we
show how we can quantize relations by mapping them with quantized
functors.

\subsection{Quantized functors}

Quantization has in our view its most natural formulation as a
property of functors between categories. We will define
quantization in the context of symmetric monoidal categories with
an action of the group $S_{2}$. The symmetries are supposed to be
symmetries in our modifies sense, they are natural isomorphisms
that satisfy the conditions given in definition \ref{mysym} .

Let now $\langle C_{i},\otimes_{i},P_{e_{i}},\alpha_{i},\beta_{i},\gamma
_{i},\sigma_{i}\rangle$ for $i=1,2$ be two symmetric monoidal categories and
let $F:C_{1}\longrightarrow C_{2}$ be a functor.

\begin{definition}
\bigskip\label{quantdef}A quantization of the functor $F$ is a triple of
natural isomorphisms $\langle\lambda,\mu,\eta\rangle$
\begin{align*}
\lambda & :\otimes_{2}\circ(F\times F)\longrightarrow F\circ\otimes_{1}\\
\mu & :F\longrightarrow tF\\
\eta & :K_{e_{2}}\longrightarrow F\circ K_{e_{1}}%
\end{align*}

such that the following relations hold
\begin{align*}
\alpha_{2}\circ1_{F\times F\times F}  & =(1_{\otimes_{2}}\circ(\lambda
^{-1}\times1_{F}))\cdot(\lambda^{-1}\circ1_{\otimes_{1}\times1_{C_{1}}}%
)\cdot(1_{F}\circ\alpha_{1})\\
& \cdot(\lambda\circ1_{1_{C_{1}}\times\otimes_{1}})\cdot(1_{\otimes_{2}}%
\circ(1_{F}\times\lambda))\\
\beta_{2}\circ1_{F\times F}  & =(1_{F}\circ\beta_{1})\cdot(\lambda
\circ1_{K_{e_{1}}\times1_{C_{1}}})\cdot(1_{\otimes_{2}}\circ(\eta\times
1_{F}))\\
\gamma_{2}\circ1_{F\times F}  & =(1_{F}\circ\gamma_{1})\cdot(\lambda
\circ1_{1_{C_{1}}\times K_{e_{1}}})\cdot(1_{\otimes_{2}}\circ(1_{F}\times
\eta))\\
\sigma_{2}\circ1_{F\times F}  & =(1_{t\otimes_{2}}\circ(\mu^{-1}\times\mu
^{-1}))\cdot(t\lambda^{-1})\cdot(\mu\circ\sigma_{1})\cdot\lambda\\
t\mu & =\mu^{-1}%
\end{align*}
\end{definition}

The only true justification of this definition, as for any mathematical
definition, lies in the importance and depth of its consequences. We will now
start investigating some of those consequences. We will first show that
quantized functors are composable.

\begin{proposition}
Let $F:C_{1}\longrightarrow C_{2}$ and $G:C_{2}\longrightarrow C_{3}$ be
quantized functors with quantizations $\langle\lambda_{F},\mu_{F},\eta
_{F}\rangle$ and $\langle\lambda_{G},\mu_{G},\eta_{G}\rangle$. Then $G\circ F$
is a quantized functor with quantization $\langle\lambda_{G\circ F}%
,\mu_{G\circ F},\eta_{G\circ F}\rangle$ where
\begin{align*}
\lambda_{G\circ F}  & =(1_{G}\circ\lambda_{F})\cdot(\lambda_{G}\circ1_{F\times
F})\\
\mu_{G\circ F}  & =\mu_{G}\circ\mu_{F}\\
\eta_{G\circ F}  & =(1_{G}\circ\eta_{F})\cdot\eta_{G}%
\end{align*}
\end{proposition}

\begin{proof}
For the first condition we have
\begin{align*}
& \alpha_{3}\circ1_{G\circ F\times G\circ F\times G\circ F}\\
& =\alpha_{3}\circ1_{G\times G\times G}\circ1_{F\times F\times F}\\
& =[(1_{\otimes_{3}}\circ(\lambda_{G}^{-1}\times1_{G})\circ1_{F\times F\times
F})\cdot(\lambda_{G}^{-1}\circ1_{\otimes_{2}\times1_{C_{2}}}\circ1_{F\times
F\times F})\\
& \cdot(1_{G}\circ\alpha_{2}\circ1_{F\times F\times F})\cdot(\lambda_{G}%
\circ1_{C_{2}\times\otimes_{2}\circ1_{F\times F\times F}})\\
& \cdot(1_{\otimes_{3}}\circ(1_{G}\times\lambda_{G})\circ1_{F\times F\times
F})]\\
& =(1_{\otimes_{3}}\circ(\lambda_{G}^{-1}\circ1_{F\times F}\times1_{G\circ
F}))\cdot(\lambda_{G}^{-1}\circ(1_{\otimes_{2}\circ(F\times F)}\times
1_{F}))\\
& \cdot(1_{G}\circ1_{\otimes_{2}}\circ(\lambda_{F}^{-1}\times1_{F}%
))\cdot(1_{G}\circ\lambda_{F}^{-1}\circ1_{\otimes_{1}\times1_{C_{1}}})\\
& \cdot(1_{G}\circ1_{F}\circ\alpha_{1})\cdot(1_{G}\circ\lambda_{F}%
\circ1_{C_{1}\times\otimes_{1}})\\
& \cdot(1_{G}\circ1_{\otimes_{2}}\circ(1_{F}\times\lambda_{F}))\cdot
(\lambda_{G}\circ(1_{F}\times1_{\otimes_{2}\circ(F\times F)}))\\
& \cdot(1_{\otimes_{3}}\circ(1_{G\circ F}\times(\lambda_{G}\circ1_{F\times
F})))\\
& =(1_{\otimes_{3}}\circ((\lambda_{G}^{-1}\circ1_{F\times F})\times1_{G\circ
F}))\cdot(\lambda_{G}^{-1}\circ(1_{\otimes_{2}\circ(F\times F)}\times
1_{F}))\\
& \cdot(1_{G\circ\otimes_{2}}\circ(\lambda_{F}^{-1}\times1_{F}))\cdot
(1_{G}\circ\lambda_{F}^{-1}\circ(1_{\otimes_{1}}\times1_{1_{C_{1}}}))\\
& \cdot(1_{G\circ F}\circ\alpha_{1})\cdot(1_{G}\circ\lambda_{F}\circ
(1_{1_{C_{1}}}\times1_{\otimes_{1}}))\cdot(1_{G\circ\otimes_{2}}\circ
(1_{F}\times\lambda_{F}))\\
& \cdot(\lambda_{G}\circ(1_{F}\times1_{\otimes_{2}\circ(F\times F)}%
))\cdot(1_{\otimes_{3}}\circ(1_{G\circ F}\times(\lambda_{G}\circ1_{F\times
F})))\\
& =(1_{\otimes_{3}}\circ((\lambda_{G}^{-1}\circ1_{F\times F})\times1_{G\circ
F}))\cdot(\lambda_{G}^{-1}\circ(\lambda_{F}^{-1}\times1_{F}))\\
& \cdot(1_{G}\circ\lambda_{F}^{-1}\circ(1_{\otimes_{1}}\times1_{1_{C_{1}}%
}))\cdot(1_{G\circ F}\circ\alpha_{1})\\
& \cdot(1_{G}\circ\lambda_{F}\circ(1_{C_{1}}\times1_{\otimes_{1}}%
))\cdot(\lambda_{G}\circ(1_{F}\times\lambda_{F}))\\
& \cdot(1_{\otimes_{3}}\circ(1_{G\circ F}\times(\lambda_{G}\circ1_{F\times
F}))\\
& =(1_{\otimes_{3}}\circ((\lambda_{G}^{-1}\circ1_{F\times F})\times1_{G\circ
F}))\cdot(1_{\otimes_{3}}\circ1_{G\times G}\circ(\lambda_{F}^{-1}\times
1_{F}))\\
& \cdot(\lambda_{G}^{-1}\circ(1_{F\circ\otimes_{1}}\times1_{F}))\cdot
(1_{G}\circ\lambda_{F}^{-1}\circ(1_{\otimes_{1}}\times1_{1_{C_{1}}}))\\
& \cdot(1_{G\circ F}\circ\alpha_{1})\cdot(1_{G}\circ\lambda_{F}\circ
(1_{1_{C_{1}}}\times1_{\otimes_{1}}))\\
& \cdot(\lambda_{G}\circ(1_{F}\times1_{F\circ\otimes_{1}}))\cdot
(1_{\otimes_{3}}\circ1_{G\times G}\circ(1_{F}\times\lambda_{F}))\\
& \cdot(1_{\otimes_{3}}\circ(1_{G\circ F}\times(\lambda_{G}\circ1_{F\times
F})))\\
& =(1_{\otimes_{3}}\circ([(\lambda_{G}^{-1}\circ1_{F\times F})\cdot(1_{G}%
\circ\lambda_{F}^{-1})]\times1_{G\circ F})\\
& \cdot([(\lambda_{G}^{-1}\circ1_{F\times F})\cdot(1_{G}\circ\lambda_{F}%
^{-1})]\circ(1_{\otimes_{1}}\times1_{1_{C_{1}}}))\\
& \cdot(1_{G\circ F}\circ\alpha_{1})\cdot([(1_{G}\circ\lambda_{F}%
)\cdot(\lambda_{G}\circ1_{F\times F})]\circ(1_{1_{C_{1}}}\times1_{\otimes_{1}%
}))\\
& \cdot(1_{\otimes_{3}}\circ(1_{G\circ F}\times\lbrack(1_{G}\circ\lambda
_{F})\cdot(\lambda_{G}\circ1_{F\times F})])\\
& =(1_{\otimes_{3}}\circ(\lambda_{G\circ F}^{-1}\times1_{G\circ F}%
))\cdot(\lambda_{G\circ F}^{-1}\circ1_{\otimes_{1}\times1_{C_{1}}})\\
& \cdot(1_{G\circ F}\circ\alpha_{1})\cdot(\lambda_{G\circ F}\circ1_{1_{C_{1}%
}\times\otimes_{1}})\cdot(1_{\otimes_{3}}\circ(1_{G\circ F}\times
\lambda_{G\circ F}))
\end{align*}

The proof for the second and third conditions are the similar and we only show
the proof for the third condition
\begin{align*}
& \gamma_{3}\circ1_{G\circ F\times1_{T}}\\
& =\gamma_{3}\circ1_{G\times1_{T}}\circ1_{F\times1_{T}}\\
& =[(1_{G}\circ\gamma_{2})\cdot(\lambda_{G}\circ1_{1_{C_{2}}\times K_{e_{2}}%
})\cdot(1_{\otimes_{3}}\circ(1_{G}\times\eta_{G}))]\circ1_{F\times1_{T}}\\
& =(1_{G}\circ\gamma_{2}\circ1_{F\times1_{T}})\cdot(\lambda_{G}\circ
1_{1_{C_{2}}\times K_{e_{2}}}\circ1_{F\times1_{T}})\cdot(1_{\otimes_{3}}%
\circ(1_{G}\times\eta_{G})\circ1_{F\times1_{T}})\\
& =(1_{G}\circ\lbrack(1_{F}\circ\gamma_{1})\cdot(\lambda_{F}\circ1_{1_{C_{1}%
}\times K_{e_{1}}})\cdot(1_{\otimes_{2}}\circ(1_{F}\times\eta_{F}))])\\
& \cdot(\lambda_{G}\circ1_{1_{C_{2}}\times K_{e_{2}}}\circ1_{F\times1_{T}%
})\cdot(1_{\otimes_{3}}\circ(1_{G}\times\eta_{G})\circ1_{F\times1_{T}})\\
& =(1_{G}\circ1_{F}\circ\gamma_{1})\cdot(1_{G}\circ\lambda_{F}\circ
1_{1_{C_{1}}\times K_{e_{1}}})\cdot(1_{G}\circ1_{\otimes_{2}}\circ(1_{F}%
\times\eta_{F}))\\
& \cdot(\lambda_{G}\circ1_{1_{C_{2}}\times K_{e_{2}}}\circ1_{F\times1_{T}%
})\cdot(1_{\otimes_{3}}\circ(1_{G}\times\eta_{G})\circ1_{F\times1_{T}})\\
& =(1_{G\circ F}\circ\gamma_{1})\cdot(1_{G}\circ\lambda_{F}\circ1_{1_{C_{1}%
}\times K_{e_{1}}})\\
& \cdot(1_{G\circ\otimes_{2}}\circ(1_{F}\times\eta_{F}))\circ(\lambda_{G}%
\circ(1_{F}\times1_{K_{e_{2}}}))\cdot(1_{\otimes_{3}}\circ(1_{G\circ F}%
\times\eta_{G}))\\
& =(1_{G\circ F}\circ\gamma_{1})\cdot(1_{G}\circ\lambda_{F}\circ1_{1_{C_{1}%
}\times K_{e_{1}}})\cdot(\lambda_{G}\circ(1_{F}\times\eta_{F}))\cdot
(1_{\otimes_{3}}\circ(1_{G\circ F}\times\eta_{G}))\\
& =(1_{G\circ F}\circ\gamma_{1})\cdot(1_{G}\circ\lambda_{F}\circ1_{1_{C_{1}%
}\times K_{e_{1}}})\cdot(\lambda_{G}\circ(1_{F}\times\eta_{F}))\\
& \cdot(1_{\otimes_{3}}\circ(1_{G\circ F}\times1_{G\circ K_{e_{2}}}%
))\cdot(1_{\otimes_{3}}\circ(1_{G\circ F}\times\eta_{G}))\\
& =(1_{G\circ F}\circ\gamma_{1})\cdot(1_{G}\circ\lambda_{F}\circ1_{1_{C_{1}%
}\times K_{e_{1}}})\cdot(\lambda_{G}\circ(1_{F}\times\eta_{F}))\\
& \cdot(1_{\otimes_{3}}\circ(1_{G\circ F}\times\eta_{G}))\\
& =(1_{G\circ F}\circ\gamma_{1})\cdot(1_{G}\circ\lambda_{F}\circ1_{1_{C_{1}%
}\times K_{e_{1}}})\cdot(\lambda_{G}\circ(1_{F}\times1_{F\circ K_{e_{1}}}))\\
& \cdot(1_{\otimes_{3}}\circ(1_{G}\times1_{G})\circ(1_{F}\times\eta_{F}%
))\cdot(1_{\otimes_{3}\circ(1_{G\circ F}}\times\eta_{G}))\\
& =(1_{G\circ F}\circ\gamma_{1})\cdot(1_{G}\circ\lambda_{F}\circ1_{1_{C_{1}%
}\times K_{e_{1}}})\cdot(\lambda_{G}\circ1_{F\times F}\circ(1_{1_{C_{1}}%
}\times1_{K_{e_{1}}}))\\
& \cdot(1_{\otimes_{3}}\circ(1_{G\circ F}\times(1_{G}\circ\eta_{F}%
)))\cdot(1_{\otimes_{3}}\circ(1_{G\circ F}\times\eta_{G}))\\
& =(1_{G\circ F}\circ\gamma_{1})\cdot([(1_{G}\circ\lambda_{F})\cdot
(\lambda_{G}\circ1_{F\times F})]\circ1_{1_{C_{1}}\times K_{e_{1}}})\\
& \cdot(1_{\otimes_{3}}\circ(1_{G\circ F}\times\lbrack(1_{G}\circ\eta
_{F})\cdot\eta_{G}]))\\
& =(1_{G\circ F}\circ\gamma_{1})\cdot(\lambda_{G\circ F}\circ1_{1_{C_{1}%
}\times K_{e_{1}}})\cdot(1_{\otimes_{3}}\circ(1_{G\circ F}\times\eta_{G\circ
F}))
\end{align*}

For the fifth condition we have
\begin{align*}
\sigma_{3}\circ1_{G\circ F\times G\circ F}\\
=\sigma_{3}\circ1_{G\times G}\circ1_{F\times F}\\
=[(1_{t\otimes_{3}}\circ(\mu_{G}^{-1}\times\mu_{G}^{-1}))\cdot(t\lambda
_{G}^{-1})\cdot(\mu_{G}\circ\sigma_{2})\cdot\lambda_{G}]\circ1_{F\times F}\\
=(1_{t\otimes_{3}}\circ(\mu_{G}^{-1}\times\mu_{G}^{-1})\circ1_{F\times
F})\cdot(t\lambda_{G}^{-1}\circ1_{F\times F})\cdot(\mu_{G}\circ\sigma_{2}%
\circ1_{F\times F})\\
\cdot(\lambda_{G}\circ1_{F\times F})\\
=(1_{t\otimes_{3}}\circ(\mu_{G}^{-1}\times\mu_{G}^{-1})\circ1_{F\times
F})\cdot(t\lambda_{G}^{-1}\circ1_{F\times F})\\
\cdot(\mu_{G}\circ\lbrack(1_{t\otimes_{2}}\circ(\mu_{F}^{-1}\times\mu_{F}%
^{-1}))\cdot(t\lambda_{F}^{-1})\cdot(\mu_{F}\circ\sigma_{1})\cdot\lambda
_{F}])\cdot(\lambda_{G}\circ1_{F\times F})\\
=(1_{t\otimes_{3}}\circ(\mu_{G}^{-1}\times\mu_{G}^{-1})\circ1_{F\times
F})\cdot(t\lambda_{G}^{-1}\circ1_{F\times F})\cdot(1_{t(G\circ\otimes_{2}%
)}\circ(\mu_{F}^{-1}\times\mu_{F}^{-1}))\\
\cdot(1_{tG}\circ(t\lambda_{F}^{-1}))\cdot(\mu_{G}\circ\mu_{F}\circ\sigma
_{1})\cdot(1_{G}\circ\lambda_{F})\cdot(\lambda_{G}\circ1_{F\times F})\\
=(1_{t\otimes_{3}}\circ(\mu_{G}^{-1}\times\mu_{G}^{-1})\circ1_{F\times
F})\cdot((t\lambda_{G}^{-1})\circ(\mu_{F}^{-1}\times\mu_{F}^{-1}))\cdot
(1_{tG}\circ(t\lambda_{F}^{-1}))\\
\cdot(\mu_{G}\circ\mu_{F}\circ\sigma_{1})\cdot(1_{G}\circ\lambda_{F}%
)\cdot(\lambda_{G}\circ1_{F\times F})\\
=(1_{t\otimes_{3}}\circ(\mu_{G}^{-1}\times\mu_{G}^{-1})\circ1_{F\times
F})\cdot(1_{t\otimes_{3}}\circ(1_{tG}\times1_{tG})\circ(\mu_{F}^{-1}\times
\mu_{F}^{-1}))\\
\cdot((t\lambda_{G}^{-1})\circ(1_{tF}\times1_{tF}))\cdot(1_{tG}\circ
(t\lambda_{F}^{-1}))\cdot(\mu_{G}\circ\mu_{F}\circ\sigma_{1})\\
\cdot(1_{G}\circ\lambda_{F})\cdot(\lambda_{G}\circ1_{F\times F})\\
=(1_{t\otimes_{3}}\circ((\mu_{G}\circ\mu_{F})^{-1}\times(\mu_{G}\circ\mu
_{F})^{-1}))\cdot(t[(1_{G}\circ\lambda_{F})\cdot(\lambda_{G}\circ1_{F\times
F})]^{-1})\\
\cdot(\mu_{G}\circ\mu_{F}\circ\sigma_{1})\cdot\lbrack(1_{G}\circ\lambda
_{F})\cdot(\lambda_{G}\circ1_{F\times F})]\\
=(1_{t\otimes_{3}}\circ(\mu_{G\circ F}^{-1}\times\mu_{G\circ F}^{-1}%
))\cdot(t\lambda_{G\circ F}^{-1})\cdot(\mu_{G\circ F}\circ\sigma_{1}%
)\circ\lambda_{G\circ F}%
\end{align*}

The last condition is clearly satisfied because action by $t$ pass through
horizontal composition.
\end{proof}

As a consequence of this proposition the class of symmetric monoidal
categories form a category where arrows are four tuples $\langle F,\lambda
_{F},\mu_{F},\eta_{F}\rangle$ and where composition of four tuples is defined
using the previous proposition.
\[
\langle G,\lambda_{G},\mu_{G},\eta_{G}\rangle\circ\langle F,\lambda_{F}%
,\mu_{F},\eta_{F}\rangle=\langle G\circ F,\lambda_{G\circ F},\mu_{G\circ
F},\eta_{G\circ F}\rangle
\]

A given category $C$ with a product bifunctor $\otimes$ and unit functor
$K_{e}$ is a symmetric monoidal category if the conditions on $\alpha
,\beta,\gamma,\sigma$ and $\theta$ stated in definition \ref{mysym} are
satisfied. These conditions are equations that may have none or many solutions
depending on the category $C$ and the choice of functors $\otimes$ and $K_{e}%
$. We thus in general have a set of solutions. Let this set be denoted by $S$.
We will now show that there is a group acting on $S$. The definition of this
group action is derived from the formulas defining a quantized functor. Let
$G$ be the following group of natural isomorphisms
\[
G=\{(\lambda,\mu,\eta)\;|\;\lambda:\otimes\longrightarrow\otimes,\mu
:1_{C}\longrightarrow1_{C},\eta:K_{e}\longrightarrow K_{e}\}
\]

where the product is taken componentwise.The size of this group depends on the
category $C$ and functors $\otimes$ and $P_{e}$. Let now $(\lambda,\mu,\eta)$
be any element of the group $G$ and define a mapping $F_{\lambda,\mu,\eta}$ on
$S$ by
\[
F_{\lambda,\mu,\eta}(\alpha,\beta,\gamma,\sigma)=(\widehat{\alpha}%
,\widehat{\beta},\widehat{\gamma},\widehat{\sigma})
\]

where
\begin{align*}
\widehat{\alpha}  & =(\lambda^{-1}\circ(\lambda^{-1}\times1_{1_{C}}%
))\cdot\alpha\cdot(\lambda\circ(1_{1_{C}}\times\lambda))\\
\widehat{\beta}  & =\beta\cdot(\lambda\circ(\eta\times1_{1_{C}})\\
\widehat{\gamma}  & =\gamma\cdot(\lambda\circ(1_{1_{C}}\times\eta))\\
\widehat{\sigma}  & =(1_{t\otimes}\circ(\mu^{-1}\times\mu^{-1}))\cdot
(t\lambda^{-1})\cdot(\mu\circ\sigma)\cdot\lambda
\end{align*}

Let $H$ be the subgroup of $G$ defined by the relations
\begin{align*}
t\mu & =\mu^{-1}\\
\mu\circ1_{\otimes}\circ(\mu^{-1}\times1_{K_{e}})  & =1_{\otimes\circ
(1_{C}\times K_{e})}\\
\mu\circ1_{\otimes}\circ(1_{K_{e}}\times\mu^{-1})  & =1_{\otimes\circ
(K_{e}\times1_{C})}\\
t\eta & =(\mu\circ1_{K_{e}})\cdot\eta
\end{align*}

Then we have the following important result.

\begin{theorem}
$F_{\lambda,\mu,\eta}:S\longrightarrow S$ and defines a action of the group
$H$ on the set $S$.
\end{theorem}

\begin{proof}
In order to prove that $(\widehat{\alpha},\widehat{\beta},\widehat{\gamma
},\widehat{\sigma})\in S$ we must show that $(\widehat{\alpha},\widehat{\beta
},\widehat{\gamma},\widehat{\sigma})$ defines a symmetric monoidal structure.
There are eight such conditions. For the first condition we have (this is also
a proof that the map $T_{1}(\sigma)$ from section \ref{symgroupaction} maps
associativity constraints to associativity constraints)%

\begin{align*}
(\widehat{\alpha}\circ1_{\otimes\times1_{C}\times1_{C}})\cdot(\widehat{\alpha
}\circ1_{1_{C}\times1_{C}\times\otimes})\\
=(\lambda^{-1}\circ(\lambda^{-1}\times1_{1_{C}})\circ1_{\otimes\times
1_{C}\times1_{C}})\cdot(\alpha\circ1_{\otimes\times1_{C}\times1_{C}})\\
\cdot(\lambda\circ(1_{C}\times\lambda)\circ1_{\otimes\times1_{C}\times1_{C}%
})\cdot(\lambda^{-1}\circ(\lambda^{-1}\times1_{C})\circ1_{1_{C}\times
1_{C}\times\otimes})\\
\cdot(\alpha\circ1_{1_{C}\times1_{C}\times\otimes})\cdot(\lambda\circ
(1_{C}\times\lambda)\circ1_{1_{C}\times1_{C}\times\otimes})\\
=(\lambda^{-1}\circ(\lambda^{-1}\times1_{1_{C}})\circ1_{\otimes\times
1_{C}\times1_{C}})\cdot(\lambda\circ(1_{\otimes}\times\lambda))\cdot
(\lambda^{-1}\circ(\lambda^{-1}\times1_{\otimes}))\\
\cdot(\alpha\circ1_{C}\times1_{C}\times\otimes)\cdot(\lambda\circ(1_{C}%
\times\lambda)\circ1_{1_{C}\times1_{C}\times\otimes})\\
=(\lambda^{-1}\circ(\lambda^{-1}\times1_{1_{C}})\circ1_{\otimes\times
1_{C}\times1_{C}})\cdot(1_{\otimes}\circ(\lambda^{-1}\times\lambda
))\cdot(\alpha\circ1_{1_{C}\times1_{C}\times\otimes})\\
\cdot(\lambda\circ(1_{C}\times\lambda)\circ1_{1_{C}\times1_{C}\times\otimes
})\\
=(\lambda^{-1}\circ((\lambda^{-1}\circ(1_{\otimes}\times1_{1_{C}}%
))\times1_{1_{C}}))\cdot(\alpha\circ1_{\otimes\times1_{C}\times1_{C}}%
)\cdot(1_{\otimes}\circ(\lambda^{-1}\times\lambda))\\
\cdot(\alpha\circ1_{1_{C}\times1_{C}\times\otimes})\cdot(\lambda\circ
(1_{1_{C}}\times(\lambda\circ(1_{1_{C}}\times1_{\otimes}))))\\
=([(\lambda^{-1}\circ(\lambda^{-1}\times1_{1_{C}}))\cdot\alpha]\circ
(1_{\otimes}\times1_{1_{C}}\times1_{1_{C}}))\\
\cdot((1_{\otimes}\circ(1_{1_{C}}\times1_{\otimes}))\circ(\lambda^{-1}%
\times1_{1_{C}}\times1_{1_{C}}))\\
\cdot((1_{\otimes}\circ(1_{\otimes}\times1_{1_{C}}))\circ(1_{1_{C}}%
\times1_{1_{C}}\times\lambda))\\
\cdot([\alpha\cdot(\lambda\circ(1_{1_{C}}\times\lambda))]\circ(1_{1_{C}}%
\times1_{1_{C}}\times1_{\otimes}))\\
=([(\lambda^{-1}\circ(\lambda^{-1}\times1_{1_{C}}))\cdot\alpha]\circ
(\lambda^{-1}\times1_{1_{C}}\times1_{1_{C}}))\\
\cdot([\alpha\cdot(\lambda\circ(1_{1_{C}}\times\lambda))]\circ(1_{1_{C}}%
\times1_{1_{C}}\times\lambda))\\
=(\lambda^{-1}\circ(\lambda^{-1}\times1_{1_{C}})\circ(\lambda^{-1}%
\times1_{1_{C}}\times1_{1_{C}}))\cdot(\alpha\circ(1_{\otimes}\times1_{1_{C}%
}\times1_{1_{C}}))\\
\cdot(\alpha\circ(1_{1_{C}}\times1_{1_{C}}\times1_{\otimes}))\cdot
(\lambda\circ(1_{1_{C}}\times\lambda)\circ(1_{1_{C}}\times1_{1_{C}}%
\times\lambda))\\
=(\lambda^{-1}\circ((\lambda^{-1}\circ(\lambda^{-1}\times1_{1_{C}}%
))\times1_{1_{C}}))\cdot(\alpha\circ1_{\otimes\times1_{C}\times1_{C}})\\
\cdot(\alpha\circ1_{1_{C}\times1_{C}\times\otimes})\cdot(\lambda\circ
(1_{1_{C}}\times(\lambda\circ(1_{1_{C}}\times\lambda))))\\
=(\lambda^{-1}\circ((\lambda^{-1}\circ(\lambda^{-1}\times1_{1_{C}}%
))\times1_{1_{C}}))\cdot(1_{\otimes}\circ(\alpha\times1_{1_{C}}))\cdot
(\alpha\circ1_{1_{C}\times\otimes\times1_{C}})\\
\cdot(1_{\otimes}\circ(1_{1_{C}}\times\alpha))\cdot(\lambda\circ(1_{1_{C}%
}\times(\lambda\circ(1_{1_{C}}\times\lambda))))\\
=(\lambda^{-1}\circ((\lambda^{-1}\circ(\lambda^{-1}\times1_{1_{C}}%
))\times1_{1_{C}}))\cdot(1_{\otimes}\circ(\alpha\times1_{1_{C}}))\\
\cdot(1_{\otimes}\circ(1_{\otimes}\times1_{1_{C}})\circ(1_{1_{C}}\times
\lambda\times1_{1_{C}}))\cdot(\alpha\circ(1_{1_{C}}\times1_{\otimes}%
\times1_{1_{C}}))\\
\cdot(1_{\otimes}\circ(1_{1_{C}}\times1_{\otimes})\circ(1_{1_{C}}\times
\lambda^{-1}\times1_{1_{C}}))\cdot(1_{\otimes}\circ(1_{1_{C}}\times\alpha))\\
\cdot(\lambda\circ(1_{1_{C}}\times(\lambda\circ(1_{1_{C}}\times\lambda))))\\
=(1_{\otimes}\circ((\lambda^{-1}\circ(\lambda^{-1}\times1_{1_{C}}%
))\times1_{1_{C}}))\cdot(\lambda^{-1}\circ((1_{\otimes}\circ(1_{\otimes}%
\times1_{1_{C}}))\times1_{1_{C}}))\\
\cdot(1_{\otimes}\circ(\alpha\times1_{1_{C}}))\cdot(1_{\otimes}\circ
((1_{\otimes}\circ(1_{1_{C}}\times\lambda))\times1_{1_{C}}))\cdot(\alpha
\circ1_{1_{C}\times\otimes\times1_{C}})\\
\cdot(\lambda\circ(1_{1_{C}}\times(1_{\otimes}\circ(\lambda^{-1}\times
1_{1_{C}}))))\cdot(1_{\otimes}\circ(1_{1_{C}}\times\alpha))\\
\cdot(1_{\otimes}\circ(1_{1_{C}}\times(\lambda\circ(1_{1_{C}}\times
\lambda))))\\
=(1_{\otimes}\circ((\lambda^{-1}\circ(\lambda^{-1}\times1_{1_{C}}%
))\times1_{1_{C}}))\cdot(1_{\otimes}\circ(\alpha\times1_{1_{C}}))\\
\cdot(\lambda^{-1}\circ(1_{\otimes}\times1_{1_{C}})\circ(1_{1_{C}}%
\times\lambda\times1_{1_{C}}))\cdot(\alpha\circ1_{1_{C}\times\otimes
\times1_{C}})\\
\cdot(\lambda\circ(1_{1_{C}}\times1_{\otimes})\circ(1_{1_{C}}\times
\lambda^{-1}\times1_{1_{C}}))\cdot(1_{\otimes}\circ(1_{1_{C}}\times\alpha))\\
\cdot(1_{\otimes}\circ(1_{1_{C}}\times(\lambda\circ(1_{1_{C}}\times\lambda))))
\end{align*}%

\begin{align*}
=(1_{\otimes}\circ((\lambda^{-1}\circ(\lambda^{-1}\times1_{1_{C}}%
))\times1_{1_{C}}))\cdot(1_{\otimes}\circ(\alpha\times1_{1_{C}}))\\
\cdot(\lambda^{-1}\circ((1_{\otimes}\circ(1_{1_{C}}\times\lambda
))\times1_{1_{C}}))\cdot(\alpha\circ1_{1_{C}\times\otimes\times1_{C}})\\
\cdot(\lambda\circ(1_{1_{C}}\times(1_{\otimes}\circ(\lambda^{-1}\times
1_{1_{C}}))))\cdot(1_{\otimes}\circ(1_{1_{C}}\times\alpha))\\
\cdot(1_{\otimes}\circ(1_{1_{C}}\times(\lambda\circ(1_{1_{C}}\times
\lambda))))\\
=(1_{\otimes}\circ((\lambda^{-1}\circ(\lambda^{-1}\times1_{1_{C}}%
))\times1_{1_{C}}))\cdot(1_{\otimes}\circ(\alpha\times1_{1_{C}}))\\
\cdot(1_{\otimes}\circ((\lambda\circ(1_{1_{C}}\times\lambda))\times1_{1_{C}%
}))\cdot(\lambda^{-1}\circ(\lambda^{-1}\times1_{1_{C}})\circ1_{1_{C}%
\times\otimes\times1_{C}})\\
\cdot(\alpha\circ1_{1_{C}\times\otimes\times1_{1_{C}}})\cdot(\lambda
\circ(1_{1_{C}}\times\lambda)\circ1_{1_{C}\times\otimes1_{C}})\\
\cdot(1_{\otimes}\circ(1_{1_{C}}\times(\lambda^{-1}\circ(\lambda^{-1}%
\times1_{1_{C}}))))\cdot(1_{\otimes}\circ(1_{1_{C}}\times\alpha))\\
\cdot(1_{\otimes}\circ(1_{1_{C}}\times(\lambda\circ(1_{1_{C}}\times
\lambda))))\\
=(1_{\otimes}\circ(\widehat{\alpha}\times1_{1_{C}}))\cdot(\widehat{\alpha
}\circ1_{1_{C}\times\otimes\times1_{C}})\cdot(1_{\otimes}\circ(1_{1_{C}}%
\times\widehat{\alpha}))
\end{align*}

This proves the first condition. For the second condition we have
\begin{align*}
(1_{\otimes}\circ(\widehat{\gamma}\times1_{1_{C}}))\cdot(\widehat{\alpha}%
\circ1_{1_{C}\times K_{e}\times1_{C}})\\
=(1_{\otimes}\circ([\gamma\cdot(\lambda\circ1_{1_{C}\times K_{e}}%
)\cdot(1_{\otimes}\circ(1_{1_{C}}\times\eta))]\times1_{1_{C}}))\cdot
([1_{\otimes}\circ(\lambda^{-1}\times1_{1_{C}}))\\
\cdot(\lambda^{-1}\circ1_{\otimes\times1_{C}})\cdot\alpha\cdot(\lambda
\circ1_{1_{C}\times\otimes})\cdot(1_{\otimes}\circ(1_{1_{C}}\times
\lambda))]\circ1_{1_{C}\times K_{e}\times1_{C}})\\
=(1_{\otimes}\circ((\gamma\times1_{1_{C}})\cdot((\lambda\circ1_{1_{C}\times
K_{e}})\times1_{1_{C}})\cdot((1_{\otimes}\circ(1_{1_{C}}\times\eta
))\times1_{1_{C}})))\cdot\\
(1_{\otimes}\circ(\lambda^{-1}\times1_{1_{C}})\circ1_{1_{C}\times K_{e}%
\times1_{C}})\cdot(\lambda^{-1}\circ1_{\otimes\times1_{C}}\circ1_{1_{C}\times
K_{e}\times1_{C}})\\
\cdot(\alpha\circ1_{1_{C}\times K_{e}\times1_{C}})\cdot(\lambda\circ
1_{1_{C}\times\otimes}\circ1_{1_{C}\times K_{e}\times1_{C}})\cdot(1_{\otimes
}\circ(1_{1_{C}}\times\lambda)\circ1_{1_{C}\times K_{e}\times1_{C}})\\
=(1_{\otimes}\circ(\gamma\times1_{1_{C}}))\cdot(1_{\otimes}\circ((\lambda
\circ1_{1_{C}\times K_{e}})\times1_{1_{C}}))\cdot(1_{\otimes}\circ
((1_{\otimes}\circ(1_{1_{C}}\times\eta))\times1_{1_{C}}))\\
\cdot(1_{\otimes}\circ(\lambda^{-1}\times1_{1_{C}})\circ1_{1_{C}\times
K_{e}\times1_{C}})\cdot(\lambda^{-1}\circ1_{\otimes\times1_{C}}\circ
1_{1_{C}\times K_{e}\times1_{C}})\\
\cdot(\alpha\circ1_{1_{C}\times K_{e}\times1_{C}})\cdot(\lambda\circ
1_{1_{C}\times\otimes}\circ1_{1_{C}\times K_{e}\times1_{C}})\cdot(1_{\otimes
}\circ(1_{1_{C}}\times\lambda)\circ1_{1_{C}\times K_{e}\times1_{C}})\\
=(1_{\otimes}\circ(\gamma\times1_{1_{C}}))\cdot(1_{\otimes}\circ((\lambda
\circ(1_{1_{C}}\times1_{K_{e}}))\times1_{1_{C}}))\\
\cdot(1_{\otimes}\circ((\lambda^{-1}\circ(1_{1_{C}}\times\eta))\times1_{1_{C}%
}))\cdot(\lambda^{-1}\circ1_{\otimes\times1_{C}}\circ1_{1_{C}\times
K_{e}\times1_{C}})\\
\cdot(\alpha\circ1_{1_{C}\times K_{e}\times1_{C}})\cdot(\lambda\circ
1_{1_{C}\times\otimes}\circ1_{1_{C}\times K_{e}\times1_{C}})\cdot(1_{\otimes
}\circ(1_{1_{C}}\times\lambda)\circ1_{1_{C}\times K_{e}\times1_{C}})\\
=(1_{\otimes}\circ(\gamma\times1_{1_{C}}))\cdot(1_{\otimes}\circ((1_{\otimes
}\circ(1_{1_{C}}\times\eta))\times1_{1_{C}}))\\
\cdot(\lambda^{-1}\circ((1_{\otimes}\circ(1_{1_{C}}\times1_{K_{e}}%
))\times1_{1_{C}}))\cdot(\lambda\circ1_{1_{C}\times\otimes}\circ1_{1_{C}\times
K_{e}\times1_{C}})\\
\cdot(1_{\otimes}\circ(1_{1_{C}}\times\lambda)\circ1_{1_{C}\times K_{e}%
\times1_{C}})\\
=(1_{\otimes}\circ(\gamma\times1_{1_{C}}))\cdot(1_{\otimes}\circ((1_{\otimes
}\circ(1_{1_{C}}\times\eta))\times1_{1_{C}}))\\
\cdot(\lambda^{-1}\circ((1_{\otimes}\circ(1_{1_{C}}\times1_{K_{e}}%
))\times1_{1_{C}}))\cdot(\lambda\circ(1_{1_{C}}\times(\lambda\circ
1_{K_{e}\times1_{C}})))\\
=(\lambda^{-1}\circ(\gamma\times1_{1_{C}}))\cdot(1_{\otimes}\circ((1_{\otimes
}\circ(1_{1_{C}}\times\eta))\times1_{1_{C}}))\\
\cdot(\alpha\circ1_{1_{C}\times K_{e}\times1_{C}})\cdot(\lambda\circ(1_{1_{C}%
}\times(\lambda\circ1_{K_{e}\times1_{C}})))\\
=(\lambda^{-1}\circ(1_{G}\times1_{1_{C}}))\cdot(1_{\otimes}\circ(\gamma
\times1_{1_{C}}))\\
\cdot(1_{\otimes}\circ(1_{\otimes}\times1_{1_{C}})\circ(1_{1_{C}}\times
\eta\times1_{1_{C}}))\cdot(\alpha\circ1_{1_{C}\times K_{e}\times1_{C}})\\
\cdot(\lambda\circ(1_{1_{C}}\times(\lambda\circ1_{K_{e}\times1_{C}})))\\
=(\lambda^{-1}\circ(1_{G}\times1_{1_{C}}))\circ(1_{\otimes}\circ(\gamma
\times1_{1_{C}}))\\
\cdot(1_{\otimes}\circ(1_{\otimes}\times1_{1_{C}})\circ(1_{1_{C}}\times
\eta\times1_{1_{C}}))\cdot(\alpha\circ(1_{1_{C}}\times1_{K_{e}}\times1_{1_{C}%
}))\\
\cdot(\lambda\circ(1_{1_{C}}\times(\lambda\circ1_{K_{e}\times1_{C}})))\\
=(\lambda^{-1}\circ(1_{G}\times1_{1_{C}}))\cdot(1_{\otimes}\circ(\gamma
\times1_{1_{C}}))\\
\cdot(\alpha\circ(1_{1_{C}}\times1_{K_{e}}\times1_{1_{C}})\circ(1_{1_{C}%
}\times\eta\times1_{1_{C}}))\\
\cdot(\lambda\circ(1_{1_{C}}\times(\lambda\circ1_{K_{e}}c)))\\
=(\lambda^{-1}\circ(1_{G}\times1_{1_{C}}))\cdot(1_{\otimes}\circ(\gamma
\times1_{1_{C}}))\\
\cdot(\alpha\circ1_{1_{C}\times K_{e}\times1_{C}}\circ(1_{1_{C}}\times
\eta1_{1_{C}}))\\
\cdot(\lambda\circ(1_{1_{C}}\times\lambda)\circ(1_{1_{C}}\times1_{K_{e}}%
\times1_{1_{C}}))\\
=(\lambda^{-1}\circ(1_{G}\times1_{1_{C}}))\cdot(1_{\otimes}\circ(\gamma
\times1_{1_{C}}))\cdot([(\alpha\circ1_{1_{C}\times K_{e}\times1_{C}})\\
\cdot(\lambda\circ(1_{1_{C}}\times\lambda))]\circ(1_{1_{C}}\times\eta
\times1_{1_{C}}))\\
=(\lambda^{-1}\circ(1_{1_{C}}\times1_{B}))\cdot(1_{\otimes}\circ(\gamma
\times1_{1_{C}}))\cdot(\alpha\circ1_{1_{C}\times K_{e}\times1_{C}})\\
\cdot((\lambda\circ(1_{1_{C}}\times\lambda))\circ(1_{1_{C}}\times\eta
\times1_{1_{C}}))
\end{align*}%

\begin{align*}
=(\lambda^{-1}\circ(1_{1_{C}}\times1_{B}))\cdot(1_{\otimes}\circ(1_{1_{C}%
}\times\beta))\\
\cdot(\lambda\circ(1_{1_{C}}\times\lambda)\circ(1_{1_{C}}\times\eta
\times1_{1_{C}}))\\
=(\lambda^{-1}\circ(1_{1_{C}}\times\beta))\cdot(\lambda\circ(1_{1_{C}}%
\times\lambda)\circ(1_{1_{C}}\times\eta\times1_{1_{C}}))\\
=(\lambda^{-1}\circ(1_{1_{C}}\times\beta))\cdot(\lambda\circ(1_{1_{C}}%
\times\lambda\circ(\eta\times1_{1_{C}})))\\
=(1_{\otimes}\circ(1_{1_{C}}\times(\beta\cdot(\lambda\circ(\eta\times1_{1_{C}%
})))))\\
=1_{\otimes}\circ(1_{1_{C}}\times(\beta\cdot(\lambda\circ1_{K_{e}\times1_{C}%
})\cdot(1_{\otimes}\circ(\eta\times1_{1_{C}}))))\\
=1_{\otimes}\circ(1_{1_{C}}\times\widehat{\beta})
\end{align*}

This proves the second condition. For the third condition we have
\begin{align*}
\widehat{\beta}\circ1_{1_{C}\times K_{e}}\\
=[\beta\cdot(\lambda\circ1_{K_{e}\times1_{C}}\cdot(1_{\otimes}\circ(\eta
\times1_{1_{C}}))]\circ1_{1_{C}\times K}\\
=(\beta\circ1_{1_{C}\times K})\cdot(\lambda\circ(1_{K_{e}}\times1_{1_{C}%
})\circ1_{1_{C}\times K})\cdot(1_{\otimes}\circ(\eta\times1_{1_{C}}%
)\circ1_{1_{C}\times K})\\
=(\gamma\circ1_{K\times1_{C}})\cdot(\lambda\circ(1_{K_{e}}\times1_{K_{e}%
}))\cdot(1_{\otimes}\circ(\eta\times1_{K_{e}}))\\
=(\gamma\circ1_{K\times1_{C}})\cdot(\lambda\circ(1_{1_{C}}\times1_{K_{e}%
})\circ(1_{K_{e}}\times1_{1_{C}}))\cdot(1_{\otimes}\circ(1_{K_{e}}\times
\eta))\\
=(\gamma\circ1_{K\times1_{C}})\cdot(\lambda\circ(1_{1_{C}}\times1_{K_{e}%
})\circ1_{K\times1_{C}})\cdot(1_{\otimes}\circ(1_{1_{C}}\times\eta
)\circ1_{K\times1_{C}})\\
=[(\gamma\cdot(\lambda\circ1_{1_{C}\times K_{e}})\cdot(1_{\otimes}%
\circ(1_{1_{C}}\times\eta))]\circ1_{K\times1_{C}}\\
=\widehat{\gamma}\circ1_{K_{e}\times1_{C}}%
\end{align*}

This proves the third condition. The proof of the fourth condition is very
technical. In the proof we will use the following symbols
\begin{align*}
L_{1}=(1_{t\otimes}\circ(1_{1_{C}}\times t\lambda^{-1})\cdot(t\lambda
^{-1}\circ1_{t\otimes\times1_{C}})\cdot(t\alpha)\cdot(t\lambda\circ
1_{t\otimes\times1_{C}})\\
L_{2}=(1_{t\otimes}\circ(1_{1_{C}}\times t\lambda^{-1})\cdot(t\lambda
^{-1}\circ1_{t\otimes\times1_{C}})\cdot(t\alpha)\\
L_{3}=(1_{t\otimes}\circ(1_{1_{C}}\times t\lambda^{-1})\cdot(t\lambda
^{-1}\circ1_{t\otimes\times1_{C}})\\
L_{4}=(1_{t\otimes}\circ(1_{1_{C}}\times t\lambda^{-1})\\
R=((t\lambda^{-1})\circ(t\lambda^{-1}\times1_{1_{C}}))\cdot((\mu\circ
\sigma)\circ((\mu\circ\sigma)\times1_{1_{C}}))\cdot(\lambda\circ(\lambda
\times1_{1_{C}}))\\
R_{1}=((\mu\circ\sigma)\circ((\mu\circ\sigma)\times1_{1_{C}}))\cdot
(\lambda\circ(\lambda\times1_{1_{C}}))\\
R_{2}=(\lambda\circ(\lambda\times1_{1_{C}}))\\
S_{1}=\alpha^{-1}\cdot(\lambda\circ1_{\otimes\times1_{C}})\cdot(1_{\otimes
}\circ(\lambda\times1_{1_{C}}))\\
S_{2}=(\lambda^{-1}\circ1_{1_{C}\times\otimes})\cdot\alpha^{-1}\cdot
(\lambda\circ1_{\otimes\times1_{C}})\cdot(1_{\otimes}\circ(\lambda
\times1_{1_{C}}))
\end{align*}

Using these symbols we have for the fourth condition
\begin{align*}
(t\widehat{\alpha})\cdot(\widehat{\sigma}\circ(\widehat{\sigma}\times1_{1_{C}%
}))\\
=(t\widehat{\alpha})\cdot([(1_{t\otimes}\circ(\mu^{-1}\times\mu^{-1}%
))\cdot(t\lambda^{-1})\cdot(\mu\circ\sigma)\cdot\lambda]\circ([(1_{t\otimes
}\circ(\mu^{-1}\times\mu^{-1}))\\
\cdot(t\lambda^{-1})\cdot(\mu\circ\sigma)\cdot\lambda]\times1_{1_{C}}))\\
=L_{1}\cdot(1_{t\otimes}\circ(t\lambda\times1_{1_{C}}))\\
\cdot((1_{t\otimes}\circ(\mu^{-1}\times\mu^{-1}))\circ((1_{t\otimes}\circ
(\mu^{-1}\times\mu^{-1}))\times1_{1_{C}}))\cdot R\\
=L_{1}\cdot((1_{t\otimes}\circ(\mu^{-1}\times\mu^{-1}))\circ((t\lambda
\circ(\mu^{-1}\times\mu^{-1}))\times1_{1_{C}}))\cdot R\\
=L_{1}\cdot((1_{t\otimes}\circ(\mu^{-1}\times\mu^{-1}))\circ((t\lambda
\circ(\mu^{-1}\times\mu^{-1}))\times1_{1_{C}}))\\
\cdot(t\lambda^{-1}\circ(t\lambda^{-1}\times1_{1_{C}}))\cdot R_{1}\\
=L_{1}\cdot((1_{t\otimes}\circ(\mu^{-1}\times\mu^{-1}))\circ((1_{t\otimes
}\circ(\mu^{-1}\times\mu^{-1}))\times1_{1_{C}}))\cdot R_{1}\\
=L_{2}\cdot(t\lambda\circ(1_{t\otimes}\times1_{1_{C}}))\\
\cdot((1_{t\otimes}\circ(\mu^{-1}\times\mu^{-1}))\circ((1_{t\otimes}\circ
(\mu^{-1}\times\mu^{-1}))\times1_{1_{C}}))\cdot R_{1}\\
=L_{2}\cdot((1_{t\otimes}\circ(\mu^{-1}\times\mu^{-1}))\circ((1_{t\otimes
}\circ(\mu^{-1}\times\mu^{-1}))\times1_{1_{C}}))\cdot R_{1}\\
=L_{2}\cdot((1_{t\otimes}\circ(\mu^{-1}\times\mu^{-1}))\circ((1_{t\otimes
}\circ(\mu^{-1}\times\mu^{-1}))\times1_{1_{C}}))\\
\cdot((\mu\circ\sigma)\circ((\mu\circ\sigma)\times1_{1_{C}}))\cdot R_{2}\\
=L_{2}\cdot((1_{1_{C}}\circ1_{t\otimes}\circ(\mu^{-1}\times\mu^{-1}%
))\circ((1_{1_{C}}\circ1_{t\otimes}\circ(\mu^{-1}\times\mu^{-1}))\times
1_{1_{C}}))\\
\cdot((\mu\circ\sigma\circ(1_{1_{C}}\times1_{1_{C}}))\circ((\mu\circ
\sigma\circ(1_{1_{C}}\times1_{1_{C}}))\times1_{1_{C}}))\cdot R_{2}\\
=L_{2}\cdot((\mu\circ\sigma\circ(\mu^{-1}\times\mu^{-1}))\circ((\mu\circ
\sigma\circ(\mu^{-1}\times\mu^{-1}))\times1_{1_{C}}))\cdot R_{2}\\
=L_{2}\cdot((\mu\circ\sigma\circ(\mu^{-1}\times\mu^{-1}))\circ((\mu\circ
\sigma\circ(\mu^{-1}\times\mu^{-1}))\times1_{1_{C}}))\cdot R_{2}\\
=L_{2}\cdot((\mu\circ\sigma\circ(\mu^{-1}\times\mu^{-1}))\circ((\mu\circ
\sigma\circ(\mu^{-1}\times\mu^{-1}))\times(\mu\circ1_{1_{C}}\circ\mu
^{-1})))\cdot R_{2}\\
=L_{2}\cdot((1\mu\circ\sigma\circ(\mu^{-1}\times\mu^{-1}))\circ(\mu\times
\mu)\circ(\sigma\times1_{1_{C}})\circ(\mu^{-1}\times\mu^{-1}\times\mu
^{-1}))\cdot R_{2}\\
=L_{2}\cdot(\mu\circ\sigma\circ((\sigma\circ(\mu^{-1}\times\mu^{-1}))\times
\mu^{-1}))\cdot R_{2}\\
=L_{3}\cdot(1_{1_{C}}\circ t\alpha)\cdot(\mu\circ\sigma\circ(\sigma
\times1_{1_{C}}))\cdot(1_{\otimes}\circ((1_{\otimes}\circ(\mu^{-1}\times
\mu^{-1}))\times\mu^{-1}))\\
\cdot(\lambda\circ(\lambda\times1_{1_{C}}))\\
=L_{3}\cdot(\mu\circ(t\alpha\cdot(\sigma\circ(\sigma\times1_{1_{C}}%
))))\cdot(\lambda\circ((\lambda\circ(\mu^{-1}\times\mu^{-1}))\times\mu^{-1})\\
=L_{3}\cdot(\mu\circ((\sigma\circ(1_{1_{C}}\times\sigma))\cdot\alpha
^{-1}))\cdot(\lambda\circ((\lambda\circ(\mu^{-1}\times\mu^{-1}))\times\mu
^{-1})\\
=L_{3}\cdot(\mu\circ((\sigma\circ(1_{1_{C}}\times\sigma))\cdot\alpha
^{-1}))\cdot(\lambda\circ(\lambda\times1_{1_{C}})\circ(\mu^{-1}\times\mu
^{-1}\times\mu^{-1}))\\
=L_{3}\cdot(\mu\circ\sigma\circ(1_{1_{C}}\times\sigma))\cdot\alpha^{-1}%
\cdot(\lambda\circ(\lambda\times1_{1_{C}})\circ(\mu^{-1}\times\mu^{-1}%
\times\mu^{-1}))\\
=L_{3}\cdot(\mu\circ\sigma\circ(1_{1_{C}}\times\sigma))\cdot((\alpha^{-1}%
\cdot(\lambda\circ(\lambda\times1_{1_{C}})))\circ(\mu^{-1}\times\mu^{-1}%
\times\mu^{-1}))\\
=L_{3}\cdot(\mu\circ\sigma\circ(1_{1_{C}}\times\sigma))\cdot((1_{\otimes}%
\circ(1_{1_{C}}\times1_{\otimes}))\circ(\mu^{-1}\times\mu^{-1}\times\mu
^{-1}))\\
\cdot(\alpha^{-1}\cdot(\lambda\circ(\lambda\times1_{1_{C}})))\\
=L_{3}\cdot(\mu\circ\sigma\circ(1_{1_{C}}\times\sigma))\cdot(1_{\otimes}%
\circ(\mu^{-1}\times(1_{\otimes}\circ(\mu^{-1}\times\mu^{-1}))))\cdot S_{1}\\
=L_{3}\cdot(\mu\circ\sigma\circ(\mu^{-1}\times(\sigma\circ(\mu^{-1}\times
\mu^{-1}))))\cdot S_{1}\\
=L_{3}\cdot(((1_{1_{C}}\circ1_{\otimes})\cdot(\mu\circ\sigma))\circ((\mu
^{-1}\times(\sigma\circ(\mu^{-1}\times\mu^{-1})))\\
\cdot((1_{1_{C}}\times(1_{t\otimes}\circ(1_{1_{C}}\times1_{1_{C}}))))))\cdot
S_{1}%
\end{align*}%

\begin{align*}
=L_{3}\cdot((1_{1_{C}}\circ1_{t\otimes})\circ(\mu^{-1}\times(\sigma\circ
(\mu^{-1}\times\mu^{-1}))))\\
\cdot((\mu\circ\sigma)\circ(1_{1_{C}}\times(1_{t\otimes}\circ(1_{1_{C}}%
\times1_{1_{C}}))))\cdot S_{1}\\
=L_{3}\cdot(1_{t\otimes}\circ(\mu^{-1}\times\mu^{-1})\circ(1_{1_{C}}\times
(\mu\circ\sigma\circ(\mu^{-1}\times\mu^{-1}))))\\
\cdot((\mu\circ\sigma)\circ(1_{1_{C}}\times(1_{t\otimes}\circ(1_{1_{C}}%
\times1_{1_{C}}))))\cdot S_{1}\\
=L_{3}\cdot(1_{t\otimes}\circ(\mu^{-1}\times\mu^{-1})\circ(1_{1_{C}}\times
(\mu\circ\sigma\circ(\mu^{-1}\times\mu^{-1}))))\\
\cdot((\mu\circ\sigma)\circ(1_{1_{C}}\times1_{t\otimes}))\cdot S_{1}\\
=L_{3}\cdot(1_{t\otimes}\circ(\mu^{-1}\times\mu^{-1})\circ(1_{1_{C}}%
\times\lbrack(1_{t\otimes}\circ(\mu^{-1}\times\mu^{-1}))\cdot(\mu\circ
\sigma)]))\\
\cdot((\mu\circ\sigma)\circ1_{1_{C}\times\otimes})\cdot S_{1}\\
=L_{4}\cdot(t\lambda^{-1}\circ(1_{1_{C}}\times1_{t\otimes}))\cdot(1_{t\otimes
}\circ(\mu^{-1}\times\mu^{-1})\circ(1_{1_{C}}\times\lbrack(1_{t\otimes}%
\circ(\mu^{-1}\times\mu^{-1}))\\
\cdot(\mu\circ\sigma)]))\cdot((\mu\circ\sigma)\circ1_{1_{C}\times\otimes
})\cdot S_{1}\\
=L_{4}\cdot(t\lambda^{-1}\circ(\mu^{-1}\times\mu^{-1})\circ(1_{1_{C}}%
\times\lbrack(1_{t\otimes}\circ(\mu^{-1}\times\mu^{-1}))\cdot(\mu\circ
\sigma)]))\\
\cdot((\mu\circ\sigma)\circ1_{1_{C}\times\otimes})\cdot S_{1}\\
=L_{4}\cdot(t\lambda^{-1}\circ(\mu^{-1}\times\mu^{-1})\circ(1_{1_{C}}%
\times\lbrack(1_{t\otimes}\circ(\mu^{-1}\times\mu^{-1}))\cdot(\mu\circ
\sigma)]))\\
\cdot([(\mu\circ\sigma)\cdot\lambda]\circ1_{1_{C}\times\otimes})\cdot S_{2}\\
=(1_{t\otimes}\circ(1_{1_{C}}\times t\lambda^{-1}))\cdot(t\lambda^{-1}%
\circ(\mu^{-1}\times\mu^{-1})\circ(1_{1_{C}}\times\lbrack(1_{t\otimes}%
\circ(\mu^{-1}\times\mu^{-1}))\\
\cdot(\mu\circ\sigma)]))\cdot([(\mu\circ\sigma)\cdot\lambda]\circ
1_{1_{C}\times\otimes})\cdot S_{2}\\
=(t\lambda^{-1}\circ(\mu^{-1}\times\mu^{-1})\circ(1_{1_{C}}\times\lbrack
t\lambda^{-1}\cdot(1_{t\otimes}\circ(\mu^{-1}\times\mu^{-1}))\cdot(\mu
\circ\sigma)]))\\
\cdot([(\mu\circ\sigma)\cdot\lambda]\circ1_{1_{C}\times\otimes})\cdot S_{2}\\
=([(1_{t\otimes}\cdot t\lambda^{-1})\circ((\mu^{-1}\times\mu^{-1}%
)\cdot(1_{1_{C}}\times1_{1_{C}}))]\circ(1_{1_{C}}\times\lbrack(t\lambda
^{-1}\circ(\mu^{-1}\times\mu^{-1}))\\
\cdot(\mu\circ\sigma)]))\cdot([(\mu\circ\sigma)\cdot\lambda]\circ
1_{1_{C}\times\otimes})\cdot S_{2}\\
=([(1_{t\otimes}\circ(\mu^{-1}\times\mu^{-1}))\cdot t\lambda^{-1}%
]\circ(1_{1_{C}}\times\lbrack(1_{t\otimes}\circ(\mu^{-1}\times\mu^{-1}))\cdot
t\lambda^{-1}\cdot(\mu\circ\sigma)]))\\
\cdot([(\mu\circ\sigma)\cdot\lambda]\circ1_{1_{C}\times\otimes})\cdot S_{2}\\
=(\widehat{\sigma}\circ(1_{1_{C}}\times\lbrack(1_{t\otimes}\circ(\mu
^{-1}\times\mu^{-1}))\cdot t\lambda^{-1}\cdot(\mu\circ\sigma)]))\\
\cdot([(\mu\circ\sigma)\cdot\lambda]\circ1_{1_{C}\times\otimes})\cdot S_{2}\\
=(\widehat{\sigma}\circ(1_{1_{C}}\times\widehat{\sigma}))\cdot(1_{\otimes
}\circ(1_{1_{C}}\times\lambda^{-1}))\cdot S_{2}\\
=(\widehat{\sigma}\circ(1_{1_{C}}\times\widehat{\sigma}))\cdot\widehat{\alpha
}^{-1}%
\end{align*}

This proves the fourth condition. The fifth and sixth condition is proved in a
similar way and we only prove the sixth.
\begin{align*}
(t\widehat{\beta})\cdot(\widehat{\sigma}\circ1_{1_{C}\times K_{e}})\\
=t\beta\cdot(t\lambda\circ(1_{1_{C}}\times t\eta))\cdot(((1_{t\otimes}%
\circ(\mu^{-1}\times\mu^{-1}))\cdot(t\lambda^{-1})\cdot(\mu\circ\sigma
)\cdot\lambda)\circ1_{1_{C}\times K_{e}})\\
=\gamma\cdot(t\sigma\circ1_{1_{C}\times K_{e}})\cdot(t\lambda\circ(1_{1_{C}%
}\times t\eta))\cdot(1_{t\otimes}\circ(\mu^{-1}\times\mu^{-1})\circ
1_{1_{C}\times K_{e}})\\
\cdot(t\lambda^{-1}\circ1_{1_{C}\times K_{e}})\cdot(\mu\circ\sigma
\circ1_{1_{C}\times K_{e}})\cdot(\lambda\circ1_{1_{C}\times K_{e}})\\
=\gamma\cdot(\mu\circ\lambda\circ(\mu^{-1}\times(\mu^{-1}\circ t\eta)))\\
=\gamma\cdot(\mu\circ\lambda\circ(\mu^{-1}\times\eta))\\
=\gamma\cdot(\lambda\circ(1_{1_{C}}\times\eta))\cdot(\mu\circ1_{\otimes}%
\circ(\mu^{-1}\times1_{K_{e}}))\\
=\gamma\cdot(\lambda\circ(1_{1_{C}}\times\eta))\\
=\widehat{\gamma}%
\end{align*}

For the seventh condition we have
\begin{align*}
t\widehat{\sigma}\\
=t((1_{t\otimes}\circ(\mu^{-1}\times\mu^{-1}))\cdot t\lambda^{-1}\cdot
(\mu\circ\sigma)\cdot\lambda)\\
=(1_{\otimes}\circ(\mu\times\mu))\cdot\lambda^{-1}\cdot(\mu^{-1}\circ
\sigma^{-1})\cdot t\lambda\\
=(1_{\otimes}\circ(\mu\times\mu))\cdot(\lambda^{-1}\circ(1_{1_{C}}%
\times1_{1_{C}}))\cdot((\mu^{-1}\circ\sigma^{-1})\circ(1_{1_{C}}\times
1_{1_{C}}))\\
\cdot(t\lambda\circ(1_{1_{C}}\times1_{1_{C}}))\\
=(\lambda^{-1}\cdot(\mu^{-1}\circ\sigma^{-1})\cdot t\lambda\cdot1_{t\otimes
})\circ(1_{1_{C}\times1_{C}}\cdot(\mu\times\mu))\\
=\lambda^{-1}\cdot(\mu^{-1}\circ\sigma^{-1})\cdot t\lambda\cdot(1_{t\otimes
}\circ(\mu\times\mu))\\
=\widehat{\sigma}^{-1}%
\end{align*}
\end{proof}

>From this point of view the quantizations of the identity functor on a
symmetric monoidal category $\langle C,\otimes,K_{e},\alpha,\beta
,\gamma,\sigma\rangle$ is exactly equal to the subgroup of $H$ that fix the
point $(\alpha,\beta,\gamma,\sigma)$. \

\subsection{Quantization of algebraic structures in symmetric monoidal categories}

Let $\langle C_{i},\otimes_{i},P_{e_{i}},\alpha_{i},\beta_{i},\gamma
_{i},\sigma_{i}\rangle$ be symmetric monoidal categories for $i=1,2$ and let
$F:C_{1}\longrightarrow C_{2}$ be a quantized functor with quantization
$(\lambda,\mu,\eta)$. Let the $S_{2}$ action on $C_{1}$ and $C_{2}$ be
generated by the functors $T_{1}:C_{1}\longrightarrow C_{1}$ and $T_{2}%
:C_{2}\longrightarrow C_{2}$. In this section we will work with objects and
need the object formulation of the conditions defining a symmetric monoidal
category and quantized functors. We collect these conditions in the following
proposition whose proof consists of applying the definition of vertical
composition and horizontal composition.

\begin{proposition}
\label{objectform}
\begin{align*}
(\alpha_{2})_{F(X),F(Y),F(Z)}=(\lambda_{X,Y}^{-1}\otimes_{2}1_{F(Z)}%
)\circ\lambda_{X\otimes_{1}Y,Z}^{-1}\circ F((\alpha_{1})_{X,Y,Z})\\
\circ\lambda_{X,Y\otimes_{1}Z}\circ(1_{F(X)}\otimes_{2}\lambda_{Y,Z})\\
(\beta_{2})_{F(X),F(Y)}=F(\beta_{X;Y})\circ\lambda_{e_{1},Y}\circ(\eta
_{X}\otimes_{2}1_{F(Y)})\\
(\gamma_{2})_{F(X),F(Y)}=F((\gamma_{1})_{X,Y})\circ\lambda_{X,e_{1}}%
\circ(1_{F(X)}\otimes_{2}\eta_{Y})\\
(\sigma_{2})_{F(X),F(Y)}=T_{2}(T_{2}(\mu_{Y}^{-1})\otimes_{2}T_{2}(\mu
_{X}^{-1}))\circ T_{2}(\lambda_{T_{1}(Y),T_{1}(X)}^{-1})\\
\circ\mu_{T_{1}(T_{1}(Y)\otimes_{1}T_{1}(X))}\circ F(\sigma_{X;Y})\circ
\lambda_{X,Y}%
\end{align*}
\end{proposition}

Quantized functors preserve algebraic structures. Let $\langle X,\nu,u\rangle$
be a monoid in the symmetric monoidal category $C_{1}$ and Define arrows in
$C_{2}$%

\begin{align*}
\nu^{\lambda}  & :F(X)\otimes_{2}F(X)\longrightarrow F(X)\\
u^{\eta}  & :e_{2}\longrightarrow F(X)
\end{align*}

by $\nu^{\lambda}=F(\nu)\circ\lambda_{X,X}$ and $u^{\eta}=F(u)\circ\eta$. In

\begin{proposition}
$\langle F(X),\nu^{\lambda},u^{\eta}\rangle$ is a monoid in $C_{2}$.
\end{proposition}

\begin{proof}
Since $\langle X,\nu,u\rangle$ is a monoid in $C_{1}$ we have the identities
\begin{align*}
\nu\circ(1_{X}\otimes_{1}\nu)=\nu\circ(\nu\otimes_{1}1_{X})\circ(\alpha
_{1})_{X,X,X}\\
\nu\circ(u\otimes_{1}1_{X})=(\beta_{1})_{X,X}\\
\nu\circ(1_{X}\otimes_{1}\nu)=(\gamma_{1})_{X,X}%
\end{align*}

If we use these identities and the relations from proposition \ref{objectform}
we have
\begin{align*}
\nu^{\lambda}\circ(1_{F(X)}\otimes_{2}\nu^{\lambda})\\
=F(\nu)\circ\lambda_{X,X}\circ(1_{F(X)}\otimes_{2}F(\nu))\circ(1_{F(X)}%
\otimes_{2}\lambda_{X,X})\\
=F(\nu)\circ F(1_{X}\otimes_{1}\nu)\circ\lambda_{X,X\otimes_{1}X}%
\circ(1_{F(X)}\otimes_{2}\lambda_{X,X})\\
=F(\nu\circ(1_{X}\otimes_{1}\nu))\circ\lambda_{X,X\otimes_{1}X}\circ
(1_{F(X)}\otimes_{2}\lambda_{X,X})\\
=F(\nu)\circ F(\nu\otimes_{1}1_{X})\circ F((\alpha_{1})_{X,X,X})\circ
\lambda_{X,X\otimes_{1}X}\circ(1_{F(X)}\otimes_{2}\lambda_{X,X})\\
=F(\nu)\circ F(\nu\otimes_{1}1_{X})\circ\lambda_{X\otimes_{1}X,X}\circ
(\lambda_{X,X}\otimes_{2}1_{F(X)})\circ(\alpha_{2})_{F(X),F(X),F(X)}\\
=F(\nu)\circ\lambda_{X,X}\circ(F(\nu)\otimes_{2}1_{F(X)})\circ(\lambda
_{X,X}\otimes_{2}1_{F(X)})\circ(\alpha_{2})_{F(X),F(X),F(X)}\\
=\nu^{\lambda}\circ(\nu^{\lambda}\otimes_{2}1_{F(X)})\circ(\alpha
_{2})_{F(X),F(X),F(X)}%
\end{align*}

and
\begin{align*}
\nu^{\lambda}\circ(u^{\lambda}\otimes_{2}1_{F(X)})\\
=F(\nu)\circ\lambda_{X,X}\circ(F(u)\otimes_{2}F(1_{X}))\circ(\eta_{X}%
\otimes_{2}1_{F(X)})\\
=F(\nu)\circ F(u\otimes_{1}1_{X})\circ\lambda_{e_{1},X}\circ(\eta_{X}%
\otimes_{2}1_{F(X)})\\
=F(\nu\circ(u\otimes1_{X}))\circ\lambda_{e_{1},X}\circ(\eta_{X}\otimes
_{2}1_{F(X)})\\
=F((\beta_{1})_{X,X})\circ\lambda_{e_{1},X}\circ(\eta_{X}\otimes_{2}%
1_{F(X)})\\
=(\beta_{2})_{F(X),F(X)}%
\end{align*}
\end{proof}

We call the monoid $\langle F(X),\nu^{\lambda},u^{\eta}\rangle$ a quantization
of the monoid $\langle X,\nu,u\rangle$ in $C_{1}$. Quantization of comonoids
is defined by duality. Let us assume that the monoid $\langle X,\nu,u\rangle$
is commutative. This property is preserved by quantization.

\begin{proposition}
Let $\langle X,\nu,u\rangle$ be a commutative monoid in $C_{1}$. Then $\langle
F(X),\nu^{\lambda},u^{\eta}\rangle$ is a commutative monoid in $C_{2}$.
\end{proposition}

\begin{proof}
Using the exchange identity for horizontal and vertical composition of natural
transformations, the two last conditions in the definition of quantized
functors \ref{quantdef} and the symmetry conditions $t\sigma_{i}=\sigma
_{i}^{-1},i=1,2$ we get the following identity
\[
t(1_{F}\circ\sigma_{1})\cdot(t\lambda)\cdot(1_{t\otimes_{2}}\circ(\mu\times
\mu))=(\mu\circ1_{\otimes})\cdot\lambda\cdot(t\sigma_{2}\circ1_{F\times F})
\]

The $(X,X,X)$ component of this identity is gives after application of the
functor $T_{2}$ the following relation
\begin{align*}
F((\sigma_{1})_{T_{1}(X),T_{1}(X)})\circ\lambda_{T_{1}(X),T_{1}(X)}\circ
(T_{2}(\mu_{X})\otimes_{2}T_{2}(\mu_{X}))\\
=T_{2}(\mu_{X\otimes_{1}X})\circ T_{2}(\lambda_{X,X})\circ(\sigma_{2}%
)_{T_{2}(F(X)),T_{2}(F(X))}%
\end{align*}

But then we have
\begin{align*}
T_{2}(\mu_{X})\circ(\nu^{\lambda})^{\sigma_{2}}\\
=T_{2}(\mu_{X})\circ T_{2}(F(\nu))\circ T_{2}(\lambda_{X,X})\circ(\sigma
_{2})_{T_{2}(F(X)),T_{2}(F(X))}\\
=F(T_{1}(\nu))\circ T_{2}(\mu_{X\otimes_{1}X})\circ T_{2}(\lambda_{X,X}%
)\circ(\sigma_{2})_{T_{2}(F(X)),T_{2}(F(X))}\\
=F(T_{1}(\nu))\circ F((\sigma_{1})_{T_{1}(X),T_{1}(X)})\circ\lambda
_{T_{1}(X),T_{1}(X)}\circ(T_{2}(\mu_{X})\otimes_{2}T_{2}(\mu_{X}))\\
=(\nu^{\sigma_{1}})^{\lambda}\circ(T_{2}(\mu_{X})\otimes_{2}T_{2}(\mu_{X}))
\end{align*}

Since $\langle X,\nu,u\rangle$ is commutative in $C_{1}$ there exists a
isomorphism $\varphi:T_{1}(X)\longrightarrow X$ such that the following
identity holds
\[
\varphi\circ\nu^{\sigma_{1}}=\nu\circ(\varphi\otimes_{1}\varphi)
\]

Let \ the isomorphism $\widehat{\varphi}:T_{2}(F(X))\longrightarrow F(X)$ be
defined by $\widehat{\varphi}=F(\varphi)\circ T_{2}(\mu_{X})$. For this
isomorphism in $C_{2}$ we have
\begin{align*}
\widehat{\varphi}\circ(\nu^{\lambda})^{\sigma_{2}}\\
=F(\varphi)\circ T_{2}(\mu_{X})\circ(\nu^{\lambda})^{\sigma_{2}}\\
=F(\varphi)\circ(\nu^{\sigma_{1}})^{\lambda}\circ(T_{2}(\mu_{X})\otimes
_{2}T_{2}(\mu_{X}))\\
=F(\varphi)\circ F(\nu^{\sigma_{1}})\circ\lambda_{T_{1}(X),T_{1}(X)}%
\circ(T_{2}(\mu_{X})\otimes_{2}T_{2}(\mu_{X}))\\
=F(\nu)\circ F(\varphi\otimes_{1}\varphi)\circ\lambda_{T_{1}(X),T_{1}(X)}%
\circ(T_{2}(\mu_{X})\otimes_{2}T_{2}(\mu_{X}))\\
=F(\nu)\circ\lambda_{X,X}\circ(F(\varphi)\otimes_{2}F(\varphi))\circ(T_{2}%
(\mu_{X})\otimes_{2}T_{2}(\mu_{X}))\\
=\nu^{l}\circ(\widehat{\varphi}\otimes_{2}\widehat{\varphi})
\end{align*}

and this proves that $\langle F(X),\nu^{\lambda},u^{\lambda}\rangle$ is a
commutative monoid.
\end{proof}

Commutative comonoids will by duality also be preserved by
quantization. Similar results holds for other algebraic structures
like modules and condoles. As a special case of the above
constructions let $F=I_{C}$ and let $\langle X,\rho,u\rangle$ be a
commutative monoid in $C$. Then any element $(\lambda,\mu,\eta)$
in the group $H$ described in the previous section defines a
quantization $\langle X,\rho^{\lambda},u^{\eta}\rangle$ of the
given monoid. We thus get a whole family of quantized product and
unit structures on the object $X$. Each such quantized product and
unit does not define a commutative monoid with respect to the
original structure $\langle\alpha ,\beta,\gamma,\sigma\rangle$,
but with respect to the structure $\langle
\widehat{\alpha},\widehat{\beta},\widehat{\gamma},\widehat{\sigma}\rangle$.

\bigskip
\newif\ifabfull\abfulltrue

\end{document}